\documentclass[12pt, leqno]{amsart}

\usepackage[OT2,T1]{fontenc}
\usepackage{indentfirst}
\usepackage{amstext}
\usepackage{amsthm}
\usepackage{amsopn}
\usepackage{amsfonts}
\usepackage{amsmath}
\usepackage{latexsym}
\usepackage[francais,english]{babel}
\usepackage{amscd}
\usepackage{amssymb}
\usepackage{amsmath}
\usepackage[all,cmtip]{xy}
\usepackage{leftidx}
\usepackage{graphicx}
\usepackage{etoolbox}
\patchcmd{\section}{\normalfont\scshape\centering}{\normalfont\bfseries}{}{}
\patchcmd{\subsection}{-.5em}{.5em}{}{}

\newtheorem{theo}{{Theorem}}[section]
\newtheorem{coro}[theo]{{Corollary}}
\newtheorem{lemma}[theo]{{Lemma}}
\newtheorem{prop}[theo]{Proposition}

\theoremstyle{definition}
\newtheorem{remark}[theo]{\textbf{Remark}}
\newtheorem{defn}[theo]{Definition}
\newtheorem{example}[theo]{Example}
\numberwithin{equation}{section}

\newtheorem{notation}[theo]{Notation}
\DeclareFontEncoding{OT2}{}{} % to enable usage of cyrillic fonts
  \newcommand{\textcyr}[1]{%
    {\fontencoding{OT2}\fontfamily{wncyr}\fontseries{m}\fontshape{n}%
     \selectfont #1}}
\newcommand{\sha}{{\mbox{\textcyr{Sh}}}}

\begin{document}
\tolerance 400 \pretolerance 200 \selectlanguage{english}

\title{Isometries of lattices and Hasse principles}
\author{Eva  Bayer-Fluckiger}
\date{\today}
\maketitle

\begin{abstract}
We give necessary and sufficient conditions for an integral polynomial without linear factors to be the characteristic
polynomial of an isometry of some even, unimodular lattice of given signature. This gives rise to Hasse principle questions,
which we answer in a more general setting. As an application, we prove a Hasse principle for signatures of knots.

%\noindent {\em Keywords:} Lattices, automorphisms, characteristic polynomials, Tate-Shafarevich groups, local-global principles,
%signatures, knots, Alexander polynomials.

\medskip

%\noindent {\em MSC 2000:} 11G35, 12G05,14G05, 14G25.
\end{abstract}

\small{} \normalsize

\medskip

\selectlanguage{english}
\section{Introduction}

\bigskip
In \cite{GM}, Gross and McMullen give necessary conditions for a monic, irreducible polynomial to
be the characteristic polynomial of an isometry of some even, unimodular lattice of prescribed signature. They speculate
that these conditions may be sufficient; this is proved in \cite {BT}.  It turns out that the conditions of
Gross and McMullen are local conditions, and that (in the case of an irreducible polynomial) a local-global principle holds.

\medskip
More generally, if $F \in {\bf Z}[X]$ is a monic polynomial without linear factors, the conditions of
Gross and McMullen are still necessary. Moreover, they are also sufficient everywhere locally (see Theorem \ref{qrs}).
% this can
%be deduced from the results of \cite {BT} (see Theorem \ref {BT} of the present paper). 
However, when $F$ is reducible, the local-global principle no longer holds
in general, as shown by the following example

\medskip
\noindent
{\bf Example.} Let $F$ be the polynomial $$ (X^6 + X^5 + X^4 + X^3 + X^2 + X + 1)^2(X^6 - X^5 + X^4 - X^3 + X^2 - X + 1)^2,$$ and 
let $L$ be an even, unimodular, positive definite lattice of rank 24.
% [It is well-know that there exist exactly 24
%isomorphism classes of such lattices, and that they are isomorphic over ${\bf Z}_p$ for every
%prime number $p$]. 
The lattice $L$ has an isometry with characteristic polynomial $F$ everywhere locally, but
not globally (see Example \ref{intro 2}). 
In particular, the Leech lattice has no isometry with characteristic polynomial $F$, in spite of having such an
isometry everywhere locally. 

\medskip More generally, let $p$ be a prime number with $p \equiv 3 \ {\rm (mod \ 4})$ and set $F = \Phi_p^2 \Phi_{2p}^2$, where $\Phi_m$ denotes
the cyclotomic polynomial of the $m$-th roots of unity (the example above is the case $p = 7$); there exists an even, unimodular, positive definite lattice having
an isometry with characteristic polynomial $F$ if and only if $p \equiv  3 \ {\rm (mod \ 8})$. On the other hand, all
even, unimodular and positive definite lattices of rank ${\rm deg}(F)$ have such an isometry locally everywhere.

\medskip
Several other examples are given in \S \ref{Z}, both in the definite and indefinite cases.

\medskip
One of the aims of this paper is to give necessary and sufficient conditions for the local-global principle (also called Hasse principle) to hold, and hence for the existence of an even, unimodular lattice of given signature having an isometry
with characteristic polynomial $F$; this is done in Theorem \ref{final}. 

\medskip
We place ourselves in a more general setting : we consider isometries of lattices 
over rings of integers of global fields of characteristic $\not = 2$ with respect to a finite set of places.  The local conditions on the polynomial then also depend on
 the field and on the finite set of places  - in the case of unimodular, even lattices over the
integers, they are given by the Gross-McMullen conditions, and the same question
for rings of integers of arbitrary number fields is treated by Kirschmer in \cite{K}. We do not attempt to work out the local
conditions in general : we assume that they are satisfied, and give the obstruction to the local-global
principle.

\medskip
Even more generally, in \S \ref{obstruction} we construct an ``obstruction group" - a finite abelian group which gives rise to 
the obstruction group to the Hasse principle in several concrete situations. 

\medskip
To explain the results of the paper, let us come back to the original question~: we have a monic polynomial $F \in  {\bf Z}[X]$
without linear factors, we give ourselves a pair  of integers, $r,s \geqslant 0$.
 and we would like to know whether there
exists an even, unimodular lattice of signature $(r,s)$ having an isometry with characteristic polynomial $F$. 
%It is well-known that if $(r,s)$ is the signature of an even, unimodular lattice, then $r \equiv s \  ({\rm mod} \ 8)$.
 It
is easy to see that  all even, unimodular lattices with the same signature 
become isomorphic over ${\bf Q}$.
This leads to an easier question : let $F \in {\bf Q}[X]$ be a monic polynomial and let $q$ be a non-degenerate
quadratic form over ${\bf Q}$; under what conditions does $q$ have an isometry with characteristic polynomial $F$ ?
This question was raised, for arbitrary ground fields of characteristic not $2$, by Milnor in \cite{M}, and
an answer is given  in  \cite {B4} in the case of global fields; 
instead of just fixing the characteristic polynomial, the point of view of \cite{B4} is to fix a module over the group ring 
of the infinite cyclic group. As seen in \cite {B4},
the crucial case is the one of a {\it semi-simple} and {\it self-dual} module (the minimal polynomial is a product
of distinct, symmetric irreducible polynomials - recall that a monic polynomial $f \in K[X]$ of even degree is
said to be {\it symmetric} if $f(X) = X^{{\rm deg}(f)} f(X^{-1})$). 

\medskip
We start by studying this ``rational" question, and then come back to the ``integral" one.  Assume that $K$ is
a global field, and that the irreducible factors of $F$ are symmetric and of even degree.
%$f = \underset{i \in I} \prod f_i^{n_i}$,
%where $f_i \in K[X]$ are distinct symmetric irreducible polynomials of even degree, and set 
%$g = \underset{i \in I} \prod f_i$. 
If $q$ is a quadratic form over $K$, we
say that the {\it rational Hasse principle} holds for $q$ and $F$ if 
$q$ has a semi-simple  isometry with characteristic
polynomial $F$ if and only if such an isometry exists everywhere locally.
 Theorem \ref{rat} gives a necessary and sufficient condition for this to hold, in terms of an
obstruction group (see \S \ref {rational obstruction}). 

\medskip
Let $I$ be the set of irreducible factors of $F$. We show that the rational Hasse principle always holds if 
for all $f \in I$, the extensions $K[X]/(f)$ are
pairwise
independent over $K$ (see Corollary \ref{pairwise independent}).
%\medskip
%\noindent
%{\bf Theorem.} {\it If the extensions  $K[X]/(f_i)$ are pairwise independent over $K$, then a quadratic form has
%an isometry with characteristic polynomial $f$ and minimal polynomial $g$  if and only if such 
%an isometry exists everywhere locally.}
We also obtain an ``odd degree descent" result (see Theorem \ref{odd}) :

\medskip
\noindent
{\bf Theorem.} {\it A quadratic form has  an isometry with characteristic polynomial $F$ if
and only if such an isometry exists over a finite extension of $K$ of odd degree.}

\bigskip
We next come to the ``integral" questions, defining an integral obstruction group (see \S \ref{integral obstruction}), which
contains the rational one as a subgroup. The Hasse principle result is given in Theorem \ref{necessary and sufficient}.
Note that both the rational and integral Hasse principles are proved in a more general setting than the one
outlined in this introduction.

\bigskip The applications to the initial question are as follows.
Let $F$
be as above, with 
$F \in {\bf Z}[X]$. In  \S \ref{integral obstruction}, we define an integral obstruction group $\sha_F$. Recall that 
$I$ is the set of irreducible factors of $F$. For all $f,g \in I$, let $V_{f,g}$ be the set of prime numbers $p$ such that $f$ and $g$ have a common irreducible,
symmetric factor ${\rm (mod \ {\it p})}$. 
Consider the equivalence relation on $I$ generated by the elementary equivalence 
$$f \sim_e g \iff V_{f,g} \not = \varnothing;$$
let $\overline I$ be the set of equivalence classes, and $C(\overline I)$ the group of maps $\overline I \to {\bf Z}/2{\bf Z}$; the group $\sha_F$ is by definition the quotient of $C(\overline I)$ by
the constant maps.

\medskip
Set  ${\rm deg}(F) = 2n$, and let $(r,s)$ be a pair of integers, $r,s \geqslant 0$ such that $r+s = 2n$.
The following conditions (C 1) and (C 2) are necessary for the existence of an even, unimodular lattice with signature $(r,s)$ having an
isometry with characteristic polynomial $F$ (see Lemma \ref{c1c2})~:
%(see Gross-Mc Mullen, \cite {GM}) :

\medskip
(C 1) {\it The integers  $|F(1)|$, $|F(-1)|$ and $(-1)^nF(1)F(-1)$ are all squares.}

\medskip

Let $m(F)$ be the number of roots  $z$ of $F$ with $|z| > 1$ (counted with multiplicity). 

\medskip
(C 2)  {\it $r  \equiv s \ {\rm (mod \ 8)}$, 
$r \geqslant m(F)$, $s \geqslant m(F)$, and $m(F) \equiv r \equiv s  \ {\rm (mod \ 2)}.$}

\medskip
Assume that conditions (C 1) and (C 2) hold. We define a finite set of linear forms on $\sha_F$, and prove that there exists
an even, unimodular lattice with signature $(r,s)$ having a semi-simple isometry with characteristic polynomial $F$ if and only if one of these
linear forms is zero (see Theorem \ref {final}.) 
%Several examples are given to illustrate this - in many cases, only a small number of linear forms have
%to be determined. 

%\medskip Set  $g = \underset{i \in I}  \prod f_i$.

%\medskip Set $g = \underset{i \in I} \prod f_i$, $E = {\bf Q}[X]/(g)$ and $R = {\bf Z}[X]/(g)$. 

%\medskip
%\noindent
%{\bf Theorem.} 
%begin{coro}\label{sha 0} 
%{\it Assume that $f$ satisfies condition {\rm (C 1)}, and that $\sha_R = 0$. Then for each pair of integers $(r,s)$ such that condition {\rm (C 2)} holds,
%there exists an even, unimodular lattice with signature $(r,s)$ having an isometry with characteristic polynomial $f$.}
%and minimal polynomial $g$.}

%\medskip
%This is proved in Corollary \ref {sha 0}. When $\sha_R \not = 0$, a necessary and sufficient condition is given in Theorem \ref {final}.
%\medskip We give several examples (see Examples \ref {cyclotomic 2}, \ref {cyclotomic 3}, \ref {cyclotomic 4}, \ref {Salem 6}, \ref {Salem 10}, \ref {Salem 12}). For instance, we have

\medskip 
Several examples are given in  \S \ref{Z}. 
% involving cyclotomic polynomials and Salem polynomials. 
%For all integers $d \geqslant 1$, we denote by $\Phi_d$ the $d$-th cyclotomic polynomial. 

\medskip
\noindent
{\bf Example.} Let $f_1(X) = X^{10} + X^9 - X^7 - X^6 - X^ 5 - X^4 - X^3 + X +1$, and $f_2 = \Phi_{14}$. Set $F = f_1 f_2^2$.
Set  $L_{3,19} = (- E_8) \oplus (- E_8) \oplus H \oplus H \oplus H$, where $E_8$ is the $E_8$-lattice and
$H$ the hyperbolic lattice. It is shown in Examples \ref {Salem 10} and  \ref {Salem 10 bis} that $L_{3,19}$ has an isometry with characteristic polynomial $F$.

\medskip
\noindent
{\bf Example.} 
Let $p$ and $q$ be two distinct prime numbers, such that $p \equiv q \equiv 3 \ {\rm (mod \ 4})$. Let $n,m,t \in {\bf Z}$ with $n,m,t \geqslant 2$ and $m \not = t$, and set

$$f_1 = \Phi_{p^nq^m}, \ f_2=  \Phi_{p^nq^t}, \ {\rm and} \ \  F= f_1f_2.$$

\medskip
There exists a positive definite, even, unimodular lattice having an isometry with characteristic polynomial $F$ if and only if ${\rm ({p \over q})} = 1$
(see Examples \ref {cyclotomic 2} and \ref {cyclotomic 3}).

\medskip
With $F$ and $(r,s)$ as above, we ask a more precise question.

\medskip
\noindent
{\bf Question.} {\it Let $t \in {\rm SO}_{r,s}({\bf R})$ be an semi-simple isometry with characteristic polynomial $F$.
Does $t$ preserve an even, unimodular lattice ?}

\medskip
Still assuming that the ``local conditions" (C 1) and (C 2) hold, we define a linear form $\sha_F \to {\bf Z}/2{\bf}$ 
and we show  that the answer to the above question is affirmative if and only if this form is trivial (see Theorem \ref{final SO}).
In particular, we have (see Theorem \ref{all SO}) : 

\medskip
\noindent
{\bf Theorem.} {\it If $\sha_F = 0$, then all semi-simple elements of ${\rm SO}_{r,s}({\bf R})$ with characteristic polynomial $F$
preserve an even, unimodular lattice. }

\medskip
Finally, we give an application to {\it knot theory}. Let $\Delta \in {\bf Z}[X]$ be a symmetric polynomial with $\Delta(1) = 1$,
and let $r,s \geqslant 0$ be two integers with $r + s = {\rm deg}(\Delta)$.
The results of this paper can be used to 
 decide for which pair $(r,s)$ there exists
a knot with  Alexander polynomial $\Delta$
and signature $(r,s)$; for simplicity, we assume here that
$\Delta$ is monic, $\Delta(-1) = \pm 1$ and that $\Delta$ is a product of distinct irreducible and symmetric polynomials
(see \S \ref{knots and Seifert}). 
We have (see Corollary \ref {knot trivial sha})~:

\medskip
\noindent
{\bf Theorem.} {\it Assume that Conditions {\rm (C 1)} and {\rm (C2)} hold for $\Delta$ and $(r,s)$, and suppose
that $\sha_{\Delta} = 0$. Then there exists a knot with Alexander polynomial $\Delta$ and signature $(r,s)$.}

\medskip This is no longer the case in general if  $\sha_{\Delta} \not = 0$. In \S \ref{torus}, we discuss in
detail the case where $\Delta$ is the polynomial 

$$\Delta_{u,v} = {{(X^{uv} - 1)(X-1)} \over {(X^u-1)(X^v-1)}}$$
where $u,v > 1$ are two relatively prime odd integers. Using properties of cyclotomic fields,  the obstruction
group $\sha_{\Delta}$ and the homomorphisms associated to the local data can be determined explicitly. 
For instance, we have (see Example \ref{final pq}) :

\medskip
\noindent
{\bf Example.} Let $p$ and $q$ be two distinct prime numbers with $p \equiv q  \equiv 3 \ {\rm (mod \ 4)}$, let
$e \geqslant 1$ be an integer, and  set $\Delta = \Delta_{p^e,q}$. There exists a knot with 
Alexander polynomial $\Delta$  and signature $(r,s)$ if and only if  $r  \equiv s \ {\rm (mod \ 8)}$, and

\medskip $\bullet$  ${\rm ({p \over q})} =  - 1$, or 

\medskip
$\bullet$ ${\rm ({p \over q})} =  1$  and
$|r-s| \leqslant 2n - 4(e-1)$.

%\bigskip We also characterize the possible Milnor indices (see \S \ref {knots and Seifert}, \S \ref{torus}). 

\bigskip
\centerline {\bf Table of contents}

\bigskip

\noindent
1. The equivariant Witt group \dotfill 6

\medskip
\noindent
2. Lattices \dotfill 7

\medskip
\noindent
3. Bounded modules, semi-simplification and reduction mod $\pi$ \dotfill 7

\medskip
\noindent
4. The residue map  \dotfill 8

\medskip
\noindent
5. Lattices and discriminant forms  \dotfill 8

\medskip
\noindent
6, Local-global problems \dotfill 9

\medskip
\noindent
7. Stable factors, orthogonal decomposition and transfer  \dotfill 10

\medskip
\noindent
8. Signatures \dotfill 11

\medskip
\noindent
9. Simple modules and reduction  \dotfill 13

\medskip
\noindent
10. Twisting groups  \dotfill 13

\medskip
\noindent
11. Residue maps  \dotfill 14

\medskip
\noindent
12. Hermitian forms and Hasse-Witt invariants \dotfill 16

\medskip
\noindent
13. Obstruction \dotfill 18

\medskip
\noindent
14. Obstruction group - the rational case  \dotfill 23

\medskip
\noindent
15. Twisting groups and equivalence relations \dotfill 24

\medskip
\noindent
16. Local data - the rational case \dotfill 27

\medskip
\noindent
17. Local-global problem - the rational case  \dotfill 29

\medskip
\noindent
18. Independent extensions \dotfill 33

\medskip
\noindent
19. Local conditions \dotfill 34

\medskip
\noindent
20. Odd degree descent \dotfill 36

\medskip
\noindent
21. Obstruction group - the integral case  \dotfill 37

\medskip
\noindent
22. Local data and residue maps  \dotfill 38

\medskip
\noindent
23. Local-global problem - the integral case   \dotfill 40

\medskip
\noindent
24. The integral case - Hasse principle with additional conditions \dotfill 44

\medskip
\noindent
25. Lattices over $\bf Z$  \dotfill 46

\medskip
\noindent
26. Milnor signatures and Milnor indices \dotfill 54

\medskip
\noindent
27. Lattices and Milnor indices \dotfill 54

\medskip
\noindent
28. Knots \dotfill 58

\medskip
\noindent
29. Seifert forms \dotfill 58

\medskip
\noindent
30. Knots and Seifert forms \dotfill 59

\medskip
\noindent
31. Alexander polynomials of torus knots and indices \dotfill 61

\bigskip

\section{The equivariant Witt group}\label{witt}

\medskip
Let $K$ be a field, and let $A$ be a $K$-algebra with a $K$-linear involution 
$\sigma : A \to A$.  All $A$-modules are supposed to be finite dimensional over $K$.

\medskip

{\bf Bilinear forms compatible with a module structure}

\medskip
 Let $M$ be an $A$-module, and let $b : M \times M \to K$ be
a non-degenerate symmetric bilinear form. We  say that $(M,b)$ is 
an {\it $A$-bilinear form} if
$$b(ax,y) = b(x,\sigma(a)y)$$ for all $a \in A$ and
all $x,y \in M$. If $G$ is a group, $A = K[G]$, and if  $\sigma : K[G] \to K[G]$ is 
%the canonical involution of the group ring, namely 
the $K$-linear
involution 
%of $K[G]$ 
sending $g \in G$ to $g^{-1}$, then this becomes $b(gx,gy) = b(x,y)$ for all $g \in G$ and $x,y \in M$.
%We then say that $(M,b)$ is 
%{\it compatible} with the $A$-module structure of $M$, and $b$ is 
%called 
%an {\it $A$-bilinear form}.
%on $M$ [in the case $A = K[G]$, this is the usual notion of $G$-bilinear form on $M$]. 

\medskip
Let $V$ be a finite dimensional $K$-vector space, and let  $q : V \times V \to K$
be a non-degenerate symmetric bilinear form. We say that $M$ and $(V,q)$ are {\it compatible} if there exists
a $K$-linear isomorphism $\varphi : M \to V$ such that the bilinear form $b_{\varphi} : M \times M \to K$, defined by $b_{\varphi} = q(\varphi(x),\varphi(y))$, is an $A$-bilinear form.

\begin{example}\label{isometry 1} Take for $A$ the group ring $K[\Gamma]$, where $\Gamma$ is the infinite cyclic group, and for
$\sigma$ the canonical involution of $K[\Gamma]$; let $\gamma$ be a generator of $\Gamma$. If $b$ is a bilinear form compatible
with a module $M$, then $\gamma$ acts as an isometry of $b$. Conversely, an isometry of a bilinear form $b : V \times V \to K$ endows
$V$ with a structure of $A$-module compatible with $b$. 

\end{example}

\begin{theo}\label{Serre} Let $M$ and $q : V \times V \to K$ be as above, and assume that $M$ and $(V,q)$ are
compatible. Let $M^s$ be the semi-simplification of $M$. There exists a non-degenerate symmetric bilinear form $q$
on $M^s$ with the following two properties

\smallskip {\rm (a)} $q'$ is compatible with the $A$-module structure of $M^s$.

\smallskip {\rm (b)} $q'$ is isomorphic to $q$. 

\end{theo}

\noindent
{\bf Proof.} See \cite {S1}, Theorem 4.2.1.

\begin{example}\label{isometry semisimple} Let $A$, $\gamma$ and $(V,b)$  be as in Example \ref{isometry 1}. We
say that a polynomial in $k[X]$ is square-free if it is the product of distinct irreducible polynomials in $k[X]$. 
Let $f \in k[X]$. The
following properties are equivalent

\medskip {\rm (a)} $(V,b)$ has an isometry with characteristic polynomial $f$.

\smallskip {\rm (b)} $(V,b)$ has an isometry with characteristic polynomial $f$ and square-free minimal polynomial.

\medskip It is clear that (b) implies (a); the implication 
(a) $\implies$ (b) 
follows from Theorem \ref{Serre}.

\end{example}

\medskip

Let $A_M$ be the image of $A$ in ${\rm End}(M)$; if $M$ is compatible with a non-degenerate bilinear form, then the  kernel of the map $A \to A_M$ is stable by $\sigma$, and the algebra $A_M$ carries the induced involution $\sigma : A_M \to A_M$.

\begin{defn}
We say that an $A$-module $M$ is {\it self-dual} if  the kernel of the map $A \to A_M$
is stable by $\sigma$. 
\end{defn}

If ${\rm char}(K) \not = 2$, the notions of symmetric bilinear form and quadratic form coincide; we then use the
terminology of {\it $A$-quadratic form} instead of $A$-bilinear form.

\medskip

{\bf The equivariant Witt group}

\medskip
We denote by $W_A(K)$ the Witt group of $A$-bilinear forms (cf. \cite {BT},
Definition 3.3). 

\medskip
If $M$ is a simple $A$-module, let $W_A(K,M)$ be the subgroup of $W_A(K)$ generated by the classes of the $(A,\sigma)$-bilinear forms $(M,q)$. 
We have $$W_A(K) = \underset{M} \oplus  \ W_A(K,M),$$ where $M$ ranges over the isomorphism
classes of simple $A$-modules (see \cite {BT}, Theorem 3.12).

\medskip
More generally, if $M$ is a semi-simple $A$-module, set $W_A(K,M) = \underset{S} \oplus \ W_A(K,S)$, where $S$ ranges over the isomorphism
classes of simple $A$-modules arising in a direct sum decomposition of $M$. 
%Note that $W_A(K,M) = W_{A_M}(K)$. 

\begin{example}\label{cobordism}  Let $A = K[\Gamma]$ and $M$ as in Example \ref{isometry 1}.
If $(M,b)$ is an $A$-bilinear form, we associate to it an $A$-bilinear form $(M^s,b')$ as
in Theorem \ref{Serre}. Assuming that ${\rm char}(K) \not =2$, this is also done in Milnor's paper \cite{M}, \S 3. It follows from \cite {M}, Lemma 3.1
and Theorems 3.2 and 3.3 that the the classes of $(M,b)$ and $(M^s,b')$ in $W_A(K)$ are equal. 

\end{example}

\section{Lattices}\label{lattices}

Let $O$ be an integral domain, and let $K$ be its field of fractions; let $\Lambda$ be an $O$-algebra, and 
let $\sigma : \Lambda \to \Lambda$ be 
an $O$-linear involution. Set $\Lambda_K = \Lambda \otimes_O K$. Let $M$ be a $\Lambda_K$-module; a {\it $\Lambda$-lattice} is a 
sub-$\Lambda$-module $L$ of $M$ which is a projective $O$-module and
satisfies $KL = M$.

\medskip Let $(M,b)$ be a $\Lambda_K$-bilinear form;
if $L$ is an $\Lambda$-lattice, then so is its dual
$$L^{\sharp} = \{x \in M \ | \ b(x,L) \subset O \}.$$ We say that $L$ is {\it unimodular} if $L^{\sharp} = L$.

\section{Bounded modules, semi-simplification and reduction mod $\pi$}\label{bounded}

We keep the notation of the previous section, and assume that $O$ is a discrete valuation ring; let $\pi$ be a uniformizer, and 
let $k = O/\pi O$ be
the residue field. Set  $\Lambda_k = \Lambda \otimes_O k$.
We say that a $\Lambda_K$-module is {\it bounded} if it contains a $\Lambda$-lattice.

\medskip
Let $M$ be a bounded $\Lambda_K$-module, and let $L \subset M$ be a $\Lambda$-lattice; the quotient $L/\pi L$ is a $\Lambda_k$-module. 
The isomorphism classes of the simple $\Lambda_k$-modules occuring as quotients in a Jordan-H\"older filtration of $L/\pi L$ are
independent of the choice of the $\Lambda$-lattice $L$ in $M$; this is a generalization of the Brauer-Nesbitt theorem, see \cite {S1}, Theorem 2.2.1. 
The direct sum of these modules is called the {\it semi-simplification} of $L/\pi L$; by the above quoted result it is 
independent of the choice of $L$. It will be called the {\it reduction mod $\pi$} of
$M$, and will be denoted by $M(k)$; this $\Lambda_k$-module is defined up to a non-canonical isomorphism. 

\medskip The involution $\sigma : \Lambda \to \Lambda$ induces an involution $\Lambda_k \to \Lambda_k$, which we still denote by $\sigma$.
Recall that a $\Lambda_k$-module $N$ is said to be self-dual if the kernel of the homomorphism $\Lambda_k \to {\rm End}(N)$ is stable 
by the involution $\sigma$, and that the image of $\Lambda_k \to {\rm End}(N)$ is denoted by $ (\Lambda_{k_v})_N$. Set $\kappa(N) = (\Lambda_{k_v})_N$; if $N$ is self-dual, then $\sigma$ induces an involution $\sigma_N : \kappa(N) \to \kappa(N)$.

\section{The residue map}\label{res}

We keep the notation of the previous section.
%and assume that $O$ is a discrete valuation ring; let $k$ be the residue field. Set  $\Lambda_k = \Lambda \otimes_O k$.
We say that a $\Lambda_K$-bilinear form is {\it bounded} if it is defined on a bounded module, and we
denote by $W_{\Lambda_K}^b(K)$ the subgroup of $W_{\Lambda_K}(K)$ generated by the
classes of bounded forms.

\medskip
Let $(M,b)$ be a $\Lambda_K$-bilinear form, and let $L$ be a lattice in $M$; we say that $L$
is {\it almost unimodular} if $\pi L^{\sharp} \subset L \subset L^{\sharp}$, where $\pi$ is
a uniformizer. Every bounded $\Lambda_K$-bilinear form contains an almost unimodular
lattice (see \cite {BT}, Theorem 4.3. {\rm (i)}). Recall also from \cite {BT}, Theorem 4.3,  the following result

\begin{theo}\label{BT} The map $$\partial : W_{\Lambda_K}^b(K) \to W_{\Lambda_k}(k)$$ given by
$$[M,b] \mapsto [L^{\sharp}/L],$$ where $L$ is an almost unimodular lattice contained in $M$,  is
a homomorphism. Moreover, a $\Lambda_K$-bilinear form contains a unimodular lattice if and only if it is bounded
and the image of its Witt class by $\partial$ is zero.

\end{theo}

We call $\partial ([M,b])$ the {\it discriminant form} of the almost unimodular lattice $L$. 

\medskip Let $(V,q)$ be a quadratic form over $K$, and let $\delta \in W_{\Lambda_k}(k)$. 
We say that $(V,q)$ {\it contains an almost unimodular $\Lambda$-lattice with discriminant form $\delta$} if there exists  an isomorphism $\varphi : M \to V$  such that
$(M,b_{\varphi})$ contains an almost unimodular $\Lambda$-lattice with discriminant form $\delta$.

\section{Lattices and discriminant forms}\label{disc}

We keep the notation of \S \ref{lattices}; moreover we assume that $K$ is a global field, and that
 $O$ is the ring of integers of $K$ with respect to a finite non-empty set $\Sigma$ of places of $K$ (containing the infinite places when $K$ is
a number field). Let $V_{\Sigma}$ be the set of places of $K$ that are not in $\Sigma$. 

\medskip
We denote by $V_K$ the set of all places of $K$; if $v \in V_K$, let $K_v$ be
the completion of $K$ at $v$, let $O_v$ be the ring of integers of $K_v$,  let $k_v$ be the residue field, and set $\Lambda^v = \Lambda \otimes _O O_v$.

\medskip
Let $(M,b)$ be a $\Lambda_K$-bilinear form; if $L$ is a lattice of $M$ and 
$v \in V_K$, set $M^v = M \otimes_K K_v$  and $L^v = L \otimes _O O_v$. 

\begin{defn} We say that a $\Lambda$-lattice $L$ is {\it almost unimodular} if the lattice $L^v$ is almost unimodular for all  $v \in V_{\Sigma}$. The {\it discriminant
form} of an almost unimodular $\Lambda$-lattice $L$ is by definition the collection $\delta(L) = (\delta_v)$ of elements of
$W_{\Lambda_{k_v}}(k_v)$ where
$\delta_v$ is the discriminant form of $L^v$. 

\end{defn} 

Note that $\delta_v = 0$ for almost all $v \in V_K$. 

\begin{prop}\label{almostunimodular}  An almost unimodular lattice is unimodular if and only if its discriminant form is trivial.

\end{prop}

\noindent
{\bf Proof.} This follows from Theorem \ref{BT}. 

\section{Local-global problems}\label{localglobalproblem}

We keep the notation of the previous section. Let $(V,q)$ be a quadratic form over $K$; 
for all  places $v \in V_{\Sigma}$, let us fix $\delta = (\delta_v)$, with $\delta_v \in W_{\Gamma}(k_v)$ such that $\delta_v = 0$ for
almost all $v$. 
We say that $(V,q)$ {\it contains an almost unimodular $\Lambda$-lattice with discriminant form $\delta$} if there exists  an isomorphism $\varphi : M \to V$  such that
$(M,b_{\varphi})$ contains an almost unimodular $\Lambda$-lattice with discriminant form $\delta$.

\medskip
We consider the following local and global conditions :

\medskip
(L 1) {\it For all $v \in V_K$, the quadratic form $(V,q)\otimes_KK_v$ is compatible with the module $M \otimes_KK_v$.}

%\medskip If (L 1) holds, we obtain bounded  $A^v$-quadratic forms on the module $M^v$
% $(M^v,q^v)$ 
%over $K_v$ for all $v \in V_K$, 
%giving rise to an element of $W^b_{A}(K_v)$. 

\medskip
${\rm (L  \ 2)}_{\delta}$
 {\it For
all  $v \in V_{\Sigma}$, the quadratic form $(V,q) \otimes_K K_v$ contains an almost unimodular $\Lambda^v$-lattice with discriminant form $\delta_v$.}

%we have $$\partial_v [M^v,q^v] = \delta_v.$$}

%\bigskip

%We say that the {\it global conditions are satisfied} if conditions (G 1) and (G 2) below hold :

\bigskip
(G 1) {\it The quadratic form $(V,q)$ is compatible with the module $M$.}

%\medskip If (G 1) holds, we obtain $A$-quadratic form on the module $M$ over $K$, giving rise to bounded 
%$A^v$-quadratic form over $K_v$ for all $v \in V_{\Sigma}$.
%$(M,q)$ over $K$, giving rise to an element of $W_{A}(K)$
%such 
%that  $(M^v,q^v)$ is a bounded $A^v$-quadratic form over $K_v$ for all finite places $v \in V_K$.

\medskip
${\rm (G \ 2)}_{\delta}$
 {\it The quadratic form $(V,q)$ contains an almost unimodular $\Lambda$-lattice with discriminant form $\delta$.}

%here exists an $A$-quadratic form $(M,q)$ as above such that for all  $v \in V_{\Sigma}$, we have $$\partial_v [M^v,q^v] = \delta_v.$$

%\bigskip
%We say that the {\it Hasse principle holds}  if the local conditions imply the global conditions. 

\medskip

For all $v \in V_K$, set $M^v = M \otimes_K K_v$ and $V^v =  V \otimes _K K_v$.

\medskip

\begin{prop}\label{local disc} 
%Assume that condition  {\rm {(L 1)}} holds. 
The following are equivalent :

\medskip
{\rm {(i)}} For all $v \in V_{\Sigma}$, there exists an isomorphism $\varphi_v : M^v  \to V^v$ such that 
$$\partial_v ([M^v,b_{\varphi_v}])= \delta_v.$$

\medskip
{\rm {(ii)}} Conditions {\rm (L 1)} and   {\rm {${\rm (L  \ 2)}_{\delta}$
}} hold. 

\end{prop}

\medskip
\noindent
{\bf Proof.} This is clear. 

%\begin{coro}\label{localunimodular} Let $\delta = 0$, and assume that condition  {\rm {(L 1)}} holds. 
%The following are equivalent

%\medskip
%{\rm {(i)}} For all $v \in V_K$, the $A$-quadratic form $(M^v,q^v)$ contains a unimodular $\Lambda^v$-lattice.

%\medskip
%{\rm {(ii)}} Condition  {\rm {(L 2)}} holds. 

%\end{coro}

%\noindent
%{\bf Proof.} This is an immediate consequence of Proposition \ref{local disc}, and the fact that a lattice is unimodular if and only if
%it is almost unimodular with trivial discriminant form (see Proposition \ref{almostunimodular}). 

\begin{prop}\label{global disc} 
%Assume that {\rm (G 1)} holds. 
The following are equivalent

\medskip
{\rm {(i)}} There exists an isomorphism $\varphi : M \to V$ such that for all $v \in V_{\Sigma}$, we have 
$$\partial_v ([M,b_{\varphi}])= \delta_v.$$

\medskip
{\rm {(ii)}} Conditions {\rm (G 1)} and  {\rm {${\rm (G \ 2)}_{\delta}$
}} hold. 

\end{prop}

\noindent
{{\bf Proof.} We have {\rm {(ii)}} $\implies$ {\rm {(i)}} by definition.
%Assume now that {\rm {(i)}} holds, and let
%$N$ be a $\Lambda$-lattice in $M$. 
Condition   {\rm {(i)} implies that for all $v \in V_{\Sigma}$, there exists
a $\Lambda^v$-lattice $L^v$ in $M^v$ with discriminant form $\delta_v$. Set $$L = \{ x \in M \ | \ x \in L^v \ {\rm for \ all} \ v \in V_K \};$$
then $L$ is
an almost unimodular lattice of $(M,b_{\varphi})$ with discriminant form $\delta$, hence 
{\rm {(ii)}} holds.

%\begin{coro}\label{globalunimodular} Let $\delta = 0$, and assume that {\rm (G 1)} holds. The following are equivalent

%\medskip
%{\rm {(i)}} The $A$-quadratic form $(M,q)$ contains a  unimodular $\Lambda$-lattice.

%\medskip
%{\rm {(ii)}} Condition {\rm {(G 2)}} holds. 

%\end{coro}

%\noindent
%{\bf Proof.} This follows from Proposition \ref{global disc}, and the fact that a lattice is unimodular if and only if
%it is almost unimodular with trivial discriminant form. 

\section{Stable factors, orthogonal decomposition and transfer}\label{factors}

We keep the notation of Section \ref{witt}; in particular, $K$ is a field,  $A$ is a $K$-algebra with a $K$-linear involution $\sigma : A \to A$, and $M$ is an $A$-module. We also
assume that $M$ is {\it semi-simple}. Note that by  \cite {BT}, Theorem 3.12, a bilinear form compatible with $M$ is in the same Witt class as a
bilinear form on a semi-simple module.

\medskip
The module $M$ decomposes into the direct sum of isotypic submodules,
%$M \simeq \underset {i \in I} \oplus M_i$. 
$M \simeq  \underset {N} \oplus M_N$, where $N$ ranges over the simple factors of $M$. If $b$ is an $A$-bilinear form on $M$, and if  $N$ is a self-dual
simple factor of $M$, then the restriction of $b$ to $M_N$ is an orthogonal direct factor of $b$, and it is compatible with the module $M_N$. If $N$ a
simple factor of $M$ that is not self-dual, then $M$ has a simple factor $N'$ with $N'  \not = N$ such that the restriction of $b$ to $M_N \oplus M_{N'}$ is an orthogonal
direct factor of $b$, and is compatible with the module $M_N \oplus M_{N'}$. Summarizing, we get the following is well-known result :

\begin{prop}\label{decomposition} Let $b$ be an $A$-bilinear form on $M$. We have an orthogonal sum decomposition
$$b \simeq ( \underset{N} \oplus \ b_N ) \oplus ( \underset{N,N'} \oplus b_{N,N'}),$$
%$$b \simeq \underset{i \in I} \oplus \ b_i,$$ 
%where $E$ ranges over the set $\mathcal A_M$, and 
where in the first sum, $N$ ranges over the self-dual simple factors of $M$, and $b_N$ is compatible
with the module $M_N$; in the second sum, $N$ ranges over the simple factors of $M$ that are not self-dual, and $b_{N,N'}$ is compatible with
the module $M_N \oplus M_{N'}$. Moreover, the form  $b_{N,N'}$ is metabolic. 

\end{prop}

\medskip Recall that $A_M$ is the image of $A$ in ${\rm End}(M)$, and that we are assuming that the kernel of the map $A \to A_M$ is stable by $\sigma$;  the algebra $A_M$ carries the induced involution $\sigma : A_M \to A_M$. 

\medskip
Assume in addition that if $N$ is a simple factor of $M$, then $A_N$ is  a {\it commutative field}: therefore $A_M$
is a product of fields, finite extensions of  $K$.  Some of these are stable under $\sigma$, others come in pairs, exchanged by $\sigma$. 

%Note that a simple module $N$ is  self-dual if $A_N$ is stable by $\sigma$. 

%Some of the fields $A_N$ are stable under $\sigma$,
%others come in pairs, exchanged by $\sigma$. They are called {\it factors} of $A_M$. A factor of $A_M$ is called a {\it simple $\sigma$-stable factor}  if it is either a field stable by $\sigma$, or
%a product of two fields exchanged by $\sigma$.  The module $M$ decomposes into the direct sum of submodules $M_E$, indexed
%by the simple $\sigma$-stable factors $E$ of $A_M$. 

\medskip
We associate to the module $M$ a set $\mathcal A_M$ of $\sigma$-stable commutative $K$-algebras, as follows. The set $\mathcal A_M$
consists of the $\sigma$-stable fields $A_N$, where $N$ is a self-dual simple factor of $M$, and of the products of two fields interchanged by $\sigma$
associated to the simple factors of $M$ that are not self-dual. 

\medskip

The elements of $\mathcal A_M$ are called  simple $\sigma$-stable algebras associated to $M$. These can be of three types :

\medskip
Type (0) : A field  $E$ stable by $\sigma$, and the restriction of $\sigma$ to $E$ is the identity.

\smallskip
Type (1) : A  field  $E$ stable by $\sigma$, and the restriction of $\sigma$ to $E$ is not the identity; $E$ is then a quadratic extension of  the fixed field of $\sigma$ in $E$.

\smallskip
Type (2) : A product of two fields exchanged by $\sigma$. 

\medskip
Proposition \ref{decomposition} implies that every bilinear form $b$ compatible with $M$ decomposes as an orthogonal sum $b \simeq b(0) \oplus
b(1) \oplus b(2)$, corresponding to the factors of type 0, 1 and 2. Moreover, by Proposition \ref{decomposition}  the class of $b(2)$ in $W_A(K)$ is zero.
%without loss of generality, we can assume that $A_M$ only has factors of type 0 and 1. 

\begin{example}\label{isometry 2} Let $A = K[\Gamma]$ as in example \ref{isometry 1}, and let $\gamma$ be a generator of $\Gamma$. In this case, the simple $\sigma$-stable factors are
of the shape $K[X]/(f)$, where $f \in K[X]$ is as follows 

\smallskip
Type (0) $f(X) = X + 1$ or $X - 1$.

\smallskip
Type (1) $f \in K[X]$ is monic, irreducible, symmetric, of even degree [recall that an even degree polynomial $f$ is  symmetric if $f(X) = X^{{\rm deg}(f)}f(X^{-1})$].

\smallskip
Type (2) $f = g g^*$ for some irreducible, monic polynomial $g \in K[X]$ with non-zero constant term, and  $g^*(X) = g(0)^{-1} X^{{\rm deg}(g)} g(X^{-1})$
such that $g \not = g^*$.

\medskip We say that $f$ is a polynomial of type 0, 1 or 2. 

\medskip

Let  $E$ be a simple $\sigma$-stable factor as above, and let $N$ be an isotypic factor  of the direct sum decomposition of $M$ with  $A_N =  E$. We denote by $\gamma_N$ the image of the generator $\gamma$ by the map $A \to A_N$.
If $b$ is a bilinear form compatible with $N$, the endomorphism $\gamma_N$ is an isometry of
$b$ with minimal polynomial $f$ and characteristic polynomial $f^n$, where $n = {\rm dim}_E N$. 

\end{example}

\begin{example} 
Let $A$, $M$ and $(V,q)$ be as in Example \ref {isometry 1}. Let $f = h h'$, where $h \in K[X]$ is a product of
polynomials of type 0 and 1, and $h' \in K[X]$ is a product of polynomials of type 2. The quadratic form $(V,q)$ has
an isometry with characteristic polynomial $f$ if and only if $(V,q)$ is the orthogonal sum of a quadratic
form having an isometry with characteristic polynomial $h$, and of a hyperbolic form of dimension ${\rm deg}(h')$. 

\end{example}

\medskip

{\bf Transfer}

\medskip
Assume now that $\mathcal A_M$ has only one element, an involution invariant field $E$; note that $M$ is a finite dimensional $E$-vector space. Let $\ell : E \to K$ be a non-trivial $K$-linear map such that $\ell(\sigma (x)) = \ell(x)$ for all $x \in E$. The following is well-known :

\begin{prop}\label{transfer} A bilinear form $b$ is compatible with $M$ if and only if there exists a non-degenerate hermitian form 
$h : M \times M \to E$ on the $E$-vector space $M$ such that  $b = \ell \circ h$. 

\end{prop}

%\noindent
%{\bf Proof.} Let  $h : M \to M \to E$ be a non-degenerate hermitian form with respect to the involution $\sigma : E \to E$. Then
%$\ell \circ h : M \times M \to K$ is an $A$-bilinear form. Conversely, if $b : M \times M \to K$ is an $A$-bilinear form, then there exists a non-degenerate
%hermitian form $h : M \times M \to E$ such that $b = \ell \circ h$. 

\section{Signatures}\label{signatures}

We keep the notation of the previous section, with $K = {\bf R}$. If $(M,b)$ is an $A$-quadratic form, we define
a {\it signature} for each self-dual simple factor of $M$. Moreover, the signatures determine $(M,b)$
up to an isomorphism of $A$-quadratic forms.

\medskip
If $(M,b)$ is as above, and if $N$ is a simple self-dual factor of $M$, we set $A_N = E$; note that $E = {\bf C}$, and
that $\sigma : E \to E$ is the complex conjugation. By Proposition \ref {transfer}, there exists a non-degenerate
hermitian form $h : N \times N \to E$ such that $b = \ell \circ h$. Let $(r_N,s_N) = (2r'_N,2s'_N)$, where
$(r'_N,s'_N)$ is the signature of the hermitian form $h$. The following
is immediate

\begin{prop}\label{sum of signatures} Let $(M,b)$ be an $A$-quadratic form, and let $(r,s)$ be the signature of $b$; then

$$r - s = \underset{N} \sum r_N - s_N$$
where $N$ runs over the self-dual simple factors of $M$. Conversely, all pairs $(r_N,s_N)$ with $r_N, s_N \leqslant 2 {\rm dim}_E(M_N)$
and $r - s = \underset{N} \sum r_N - s_N$ are realized by an $A$-form 
$(M,b)$ such that the signature of $b$ is  $(r,s)$.

\end{prop}

\begin{notation}\label{index} The difference $r-s$ is called the {\it index} of $b$, denoted by $\tau(b)$. Similarly,
for each self-dual simple factor $N$ of $M$, we denote by $\tau_N(b) = r_N - s_N$ the index of $b$ at $N$. 

\end{notation}

Note that Proposition \ref{sum of signatures} implies that $\tau(b) = \underset{N} \sum \tau_N(b)$, where 
$N$ runs over the self-dual simple factors of $M$. 

\medskip
It is easy to see that the signatures characterize the $A$-quadratic form up to isomorphism. We have

\begin{prop}\label{R} Two $A$-quadratic forms on $M$ are isomorphic if and only if their
signatures coincide for each self-dual simple factor of $M$. 

\end{prop}

\noindent
{\bf Proof.} By Proposition \ref{decomposition} an $A$-quadratic form on $M$ is determined by its restriction to the
modules $M_N$, where $N$ runs over the self-dual simple factors of $M$. On the other hand, it is easy to
see that an $A$-quadratic
form  $(M_N,b_N)$ is determined by the hermitian form $h : M_N  \times M_N \to E$ with $b_N = \ell \circ h$.

\medskip
The following corollary is an immediate consequence of Proposition \ref{R} :

\begin{coro}\label{coro R} Let $(M,b)$ and $M,b')$ be two $A$-quadratic forms. Then
$$(M,b) \simeq (M,b') \iff \tau_N(b) = \tau_N(b')$$ for each self-dual simple factor $N$ of $M$. 

\end{coro}

\begin{example}\label{iso signature} We keep the notation of Example \ref{isometry 1}, with $A = {\bf R}[\Gamma]$. 
Let $(M,b)$ be an $A$-quadratic form; 
the action of a generator of $\Gamma$ on $M$ is an isometry of the quadratic form $b$.
The simple, self-dual factors  of $M$ are of the shape ${\bf R}[X]/\mathcal P$, where $\mathcal P \in {\bf R}[X]$ is
an irreducible, symmetric polynomial of degree $2$; 
hence each such polynomial $\mathcal P$ gives rise
to a signature of $(M,b)$ at $\mathcal P$.

\end{example}

\begin{example}\label {SO} The above results lead to a generalization of Theorem 2.4 of \cite{GM}. 
We keep the notation of Example \ref{isometry 1} and Example \ref{iso signature} with $A = {\bf R}[\Gamma]$. 
%%Let $f = f' f''$ 
%$f = \underset {i \in I} \prod f_i f'$ 
%be a symmetric polynomial, where $f'$ is a product of irreducible, symmetric polynomials of degree 2 (in other
%words, of polynomials of type (1), and $f''$ is a
%$f_i \in {\bf R}[X]$ is irreducible, symmetric of degree 2 (hence of type (1)) for all $i \in I$, and $f'$ is  a 
%product
%of irreducible polynomials of type (2). 
%Hence ${\rm deg}(f'')$ is even; let us write ${\rm deg}(f'') = 2m(f)$. Let $r, s \geqslant 0$ be integers with 
%$r \geqslant m(f)$, $s \geqslant m(f)$, and with $r \equiv s \equiv m(f) \ {\rm (mod \ 2)}$. 
Let $(r,s)$ be the signature of the quadratic form $b$. By the above construction, 
to all semi-simple elements of ${\rm SO}_{r,s}({\bf R})$ we associate a signature (and an index) at the irreducible, symmetric factors of its characteristic polynomial. Moreover,  Proposition
\ref{R} shows that two such elements of ${\rm SO}_{r,s}({\bf R})$ are conjugate if and only if these signatures coincide. 

\medskip If $t \in {\rm SO}_{r,s}({\bf R})$ is semi-simple and if $\mathcal P$ is an irreducible, symmetric factor of the
characteristic polynomial of $t$, we denote by $\tau_{\mathcal P}(t)$ the corresponding index. By Corollary \ref{coro R},
two semi-simple isometries $t, t' \in {\rm SO}_{r,s}({\bf R})$ with characteristic polynomial $f$ are conjugate if and only if $\tau_{\mathcal P}(t)
= \tau_{\mathcal P}(t')$ for each irreducible, symmetric factor $\mathcal P$ of $f$. 

\end{example}

\section{Simple modules and reduction}\label{simple}

We keep the notation of \S \ref{factors}; we assume moreover that $M$ is a simple $A$-module. Let $E = A_M$ be the image of
$A$ in ${\rm End}(M)$; by hypothesis, $E$ carries an involution $\sigma : E \to E$ induced by the involution of $A$. As in the previous section,
we assume that 
$E$ is a field, finite extension of $K$; let $F$ be the fixed field of $\sigma$ in $E$. 

\medskip
Assume moreover that $K$ is a global field, and let $O$ be a ring of integers with respect to a finite set $\Sigma$ of places of $K$, containing
the infinite places if $K$ is a number field.  Let $\Lambda$ be
an $O$-algebra; assume that $A = \Lambda \otimes _O K$, and that $\Lambda$ is stable by the involution $\sigma : A \to A$. 

\medskip
Recall that we denote by $V_K$ the set of places of $K$, by $V_{\Sigma}$ the set of places of $K$ that are not in $\Sigma$, and by $K_v$ the completion of $K$ at $v$; let $O_v$ be the ring
of integers of $K_v$, let $\pi_v \in O_v$ be a uniformizer, and let $k_v = O_v/\pi_v O_v$ be the residue field. Set $M^v = M \otimes _K  K_v$,
$E^v = E \otimes_K K_v$ 
and $\Lambda^v = \Lambda \otimes_O O_v$.
 As in \S \ref{bounded}, we denote by $M^v(k_v)$ the reduction mod $\pi_v$ of $M^v$. Let $S^v$ be the set of places $w$  of $F$ above $v$ such
 that $E_w = E \otimes F_w$ is a field. If $w \in S^v$, we denote by $O_w$ the ring of integers of $E_w$, and by $\kappa_w$
its residue field. 

\begin{prop}\label{reduction} Let $v \in V_{\Sigma}$. The self-dual, simple components of $M^v(k_v)$ are isomorphic to $\kappa_w$ for some $w \in S^v$. 

\end{prop}

\noindent
{\bf Proof.} Let $O_{E^v}$ be the maximal order of $E^v$, and note that $O_{E^v}$ is a $\Lambda^v$-lattice in $M^v = E^v$. The simple
components of $O_{E^v}/\pi_v O_{E^v}$ are the residue fields of the rings of integers of $E_w$ for the places $w \in V_E$ above $v$;
such a component is self-dual if and only if $w \in S^v$. This completes the proof of the proposition. 

\begin{example}\label{isometry 0}
Assume that $\Lambda = O[\Gamma]$,  hence  $A = \Lambda \otimes_O K = K[\Gamma]$. With the notation of Example \ref{isometry 2},
let $f \in O[X]$ of type 1, let $E = K[X]/(f)$ and let  $M = [ K[X]/(f)]^n$ for some 
integer $n \geqslant 1$. 
%Assume that all the simple $\sigma$-stable
%factors in $\mathcal A_M$ are of type (1). In other words, we have $\mathcal A_M = (E_i)_{i \in I}$
%\underset {i \in I} \prod E_i$ 
%with $E_i = K[X]/(f_i)$, where $f_i \in K[X]$
%are monic, irreducible,  symmetric polynomials of even degree, and $M = \underset{i \in I} \oplus M_i$, with $M_i = E_i^{n_i}$.
%Assume moreover that $f_i \in O[X]$; set $g = \underset {i \in I} \prod  f_i$, and $f = \underset {i \in I} \prod f^{n_i}$.
%If $v \in V_{\sigma}$, we denote by $O_v$ the ring of integers of $K_v$, and by $k_v$ the residue field. 
The homomorphism $O_v \to k_v$ induces $p_v : O_v[X] \to k_v[X]$. The self-dual, simple components of $M^v(k_v)$ are isomorphic to $k_v[X]/(\mathcal P)$, where $\mathcal P \in k_v[X]$ is an
irreducible, symmetric factor of the polynomial $p_v(f) \in k_v[X]$. Indeed, $[O_v[X]/(f)]^n$ is a lattice in $M^v = [K_v[X]/(f)]^n$,
hence the self-dual, simple components of the reduction mod $\pi_v$  of $M^v$ are of the shape  $k_v[X]/(\mathcal P)$ with $\mathcal P \in k_v[X]$ as above. 
\end{example}

\section{Twisting groups}\label{twist}

We start by recalling  some notions and facts from \cite{BT}, \S 5. Assume that $K$ is a field of characteristic
$\not= 2$.

\medskip
If $F$ is a commutative semi-simple  $K$-algebra of finite rank, and $E$ a $K$-algebra that is free
of rank 2 over $F$; we denote by $\sigma : E \to E$ the involution fixing $F$, and we set

$$T(E,\sigma) = F^{\times}/{\rm N}_{E/F}(E^{\times}).$$

\medskip \noindent
There is an exact sequence
\[
	1 \to T(E,\sigma) \overset\beta\to {\rm Br}(F) \overset{\rm res}\longrightarrow {\rm Br}(E),
\]
where  ${\rm res} :  {\rm Br}(F) \to {\rm Br}(E)$ is the base change map (cf. \cite{BT}, Lemma 5.3). 
Let $d \in F^{\times}$ such that $E = F(\sqrt d)$; then $\beta(\lambda)$ is the class of the
quaternion algebra $(\lambda,d)$ in ${\rm Br}(F)$.

\medskip

If $E$ is a local field, then we have a natural commutative diagram

\begin{equation}\label{e1}
\xymatrix{ 
   1  \ar[r]  & T(E,\sigma) \ar[d]_{\theta} \ar[r] & {\rm Br}(F) \ar[d]_{\rm inv} \ar[r]& {\rm Br}(E)  \ar[d]_{\rm inv}
   \\
 0 \ar[r]  & {\bf Z}/2{\bf Z} \ar[r]^-{1/2}&
   {\bf Q}/{\bf Z} \ar[r]^{2}&
    {\bf Q}/{\bf Z}}
  \end{equation} 
  in which the vertical map $\theta : T(E,\sigma) \to {\bf Z}/2 {\bf Z}$ is an isomorphism.

\medskip

Assume now that $E$ is a global field. Let $S$ be the set of places $w$ of $F$ such that
$E \otimes_F F_w$ is a field. Then we have (cf. \cite{BT}, Theorem 5.7)

\begin{theo} \label{reciprocity} The sequence

\[
	1 \to T(E,\sigma) \to \bigoplus_{w\in  S} T(E_w,\sigma) \overset{\sum \theta_w} \longrightarrow
	{\bf Z}/2 {\bf Z} \to 0
\]
is exact.

\end{theo}

\medskip
Let $v$ be a place of $K$, and let $w$ be a place of $F$ above $v$. Recall that the following diagram commutes

\begin{equation}\label{e2}
\xymatrix{
{\rm Br}(F_w) \ar[d]_{{\rm cor}} \ar[r]^{{\rm inv}_w} &  {\bf Q}/{\bf Z}\ar[d]_{\rm id}  \\
{\rm Br}(K_v) \ar[r]^{{\rm inv}_v}  &  {\bf Q}/{\bf Z}.}
\end{equation} 

\section{Residue maps}\label{residuemaps}

Assume that $K$ is a non-archimedean local field of characteristic $\not = 2$, that
$E$ is a finite field extension of $K$ with a non-trivial involution $\sigma : E \to E$ and
that $F$ is the fixed field of $\sigma$. 
Let $\ell : E \to K$ be a non-trivial $K$-linear map such
that $\ell(\sigma(x) )= \ell(x)$ for all $x \in E$. Let $O$ be the ring of integers of $K$, and let
$k$ be its residue field. Let $O_E$ be the ring of integers of $E$, and let $m_E$ be its maximal ideal; set $\kappa_E = O_E/m_E$.
%Let us write $E = K(\tau)$, and assume that the coefficients of the minimal polynomial of $\tau$ are in $O$. 

\medskip
If $n \geqslant 1$ is an integer, we define a map $$t_n: T(E,\sigma) \to W_E^b(K)$$ as follows. For
$\lambda \in T(E,\sigma)$, let us denote by $h_{n,\lambda}$ the $n$-dimensional diagonal
hermitian form $\langle \lambda, 1, \dots, 1 \rangle$ over $E$ with respect to the involution $\sigma$, and set
$q_{n,\lambda} = \ell(h_{n,\lambda})$;
% Note that the multiplication by $\tau$ is an isometry of 
%$q_{\lambda}$, and that $q_{\lambda}$ is bounded, hence 
we obtain an element
$[q_{n,\lambda}] \in W_{E}^b(K)$. 
Recall that we have a homomorphism $\partial : W_{E}^b(K) \to W_{\kappa_E}(k)$.

\begin{prop}\label{residue} For all integers $n \geqslant 1$, we have

\medskip

 {\rm (a)} If $E/F$ is inert, then 

$$\partial  \circ t_n : T(E,\sigma) \to W_{\kappa_E}(k)$$ is a bijection.

\medskip

 {\rm (b)} If $E/F$ is ramified and ${\rm char}(k) \not = 2$, then 

$$\partial  \circ t_n : T(E,\sigma) \to W_{\kappa_E}(k)$$ is injective, and its
image consists of the classes of forms of dimension $n[\kappa_E:k]$ ${\rm mod} \ 2$.

\medskip
{\rm (c)} If $E/F$ is ramified and ${\rm char}(k) = 2$, then 

$$\partial  \circ t_n : T(E,\sigma) \to W_{\kappa_E}(k)$$ is constant.

\end{prop}

\medskip
\noindent
{\bf Proof.} Let $\delta$ be the valuation of the different ideal $\mathcal {D}_{E/K}$; in other words, we
have $\mathcal {D}_{E/K} = m_E^{\delta}$. 

\medskip (a) Let $\overline \sigma : \kappa_E \to \kappa_E$ be the involution
induced by $\sigma$. Since $E/F$ is inert, this involution is non-trivial. The
 group $W_{\kappa_E}(k)$ is isomorphic to the Witt group $W(\kappa_E,\overline \sigma)$
of hermitian forms over $(\kappa_E,\overline \sigma)$. Note that since $k$ is a finite field, the group
$W(\kappa_E,\overline \sigma)$ is of order 2. 

\medskip
Let $\pi_E$ be a generator of $m_E$, and note that 
the class of $\pi_E$ is the unique non-trivial element of $T(E,\sigma)$. 
By \cite{BT}, Corollary 6.2 and Proposition 6.3, we see that
if $\delta$ is even, then $\partial  \circ t_n (1) = 0$ 
and $\partial  \circ t_n (\pi_E) \not = 0$;
if $\delta$ and $n$ are both odd, then $\partial  \circ t_n (\pi_E) = 0$ and $\partial  \circ t_n (1) \not = 0$;
if $\delta$ is odd and $n$ is even, then $\partial  \circ t_n (1) = 0$ and $\partial  \circ t_n (\pi_E) \not = 0$.
Hence $\partial \circ t_n$ is bijective.

\medskip (b) This follows from \cite{BT}, Proposition 6.6.

\medskip (c) Follows from \cite{BT}, Proposition 6.7. 

\bigskip
The following special case will be useful in the sequel.
\medskip
\begin{lemma}\label{unramified}  Assume that $E/K$ is unramified. Then $\partial \circ t_n(1) = 0$ in $W_{\kappa_E}(k)$ for all integers $n \geqslant 1$.

\end{lemma}

\medskip
\noindent
{\bf Proof.} The hypothesis implies that the valuation of the different ideal $\mathcal {D}_{E/K}$ is zero, and that
$E/F$ is inert; hence Proposition \ref{residue} yields the desired result.
%\cite{BT}, Corollary 6.2 yields the desired result. 

\section{Hermitian forms and Hasse-Witt invariants}\label{HW}

Assume that  $K$ is a field of characteristic $\not = 2$. The aim of this section is
to give some results relating the invariants of  hermitian forms and those of the
quadratic forms obtained from them via  transfer. Recall that every quadratic
form $q$ over $K$ can be diagonalized, in other words there exist $a_1,\dots,a_n \in K^{\times}$ such
that $q \simeq \langle a_1, \dots ,a_n \rangle$. The {\it determinant} of  $q$ is by definition the product $a_1\dots a_n$,
denoted by ${\rm det}(q)$; 
it is an element of $K^{\times}/K^{\times 2}$. Let us denote by ${\rm Br}(K)$ the Brauer group of $K$,
considered as an additive abelian group, and let ${\rm Br_2}(K)$ be the subgroup of elements
of order $\le 2$ of ${\rm Br}(K)$. The {\it Hasse-Witt invariant} of $q$ is
by definition $w_2(q) = \Sigma_{i < j} (a_i,a_j) \in {\rm Br_2}(k)$, where $(a_i,a_j)$ is the class of the quaternion
algebra over $K$ determined by $a_i,a_j$.

\medskip
Let $E$ be a finite extension of $K$, and let $\sigma : E \to E$
be a non-trivial $K$-linear involution of $E$; 
let $F$ be the fixed field of $\sigma$, and let $d \in F^{\times}$ be such that $E = F(\sqrt d)$. Let 
$\ell : E \to K$ be a non-trivial $K$--linear map such that 
$\ell (x) = \ell (\sigma( x))$ for all $x \in E$.

\medskip

For all $a \in F^{\times}$, let $D_a$ be the determinant of the quadratic form $F \times F \to K$, defined by $(x,y) \mapsto \ell_F(axy)$, where
$\ell_F$ is the restriction of $\ell$ to $F$.

\begin{lemma}\label{det F} For all $a \in F^{\times}$, we have $$D_a = {\rm N}_{F/K}(a) D_1.$$

\end{lemma}

\noindent
{\bf Proof.} Let $L_a : F \to {\rm Hom}_K(F,K)$ be defined by $L_a(x)(y) = \ell_F(axy)$ for all $x,y \in F$, and let $m_a : F \to F$ be the
multiplication by $a$; we have $L_a = L_1 \circ m_a$. Since ${\rm det}(L_a) = D_a$ and ${\rm det}(m_a) = {\rm N}_{F/K}(a)$,
the lemma follows. 

\medskip
If $h : M \times M \to E$ is a non-degenerate hermitian form, composing with $\ell$ gives rise  to a quadratic form $\ell(h) : M \times M \to K$
defined by $\ell(h)(x,y) = \ell(h(x,y))$. 

\begin{prop}\label{det} Let $h : M \times M \to E$ be a non-degenerate hermitian form, and let $n = {\rm dim}(M)$. The determinant of
the quadratic form $\ell(h)$ is $${\rm N}_{F/K}(-d)^n.$$

\end{prop}

\noindent
{\bf Proof.} It suffices to prove the proposition when $n = 1$. Let $\lambda \in F^{\times}$ be such that $h(x,y) = \lambda x \sigma(y)$. Since
$\ell(\sigma(x) = \ell(x)$ for all $x \in E$, we have $\ell(\sqrt d) = 0$, hence $\ell(h)$ is the orthogonal sum of two quadratic forms defined
on the $K$-vector space $F$, namely $(x,y) \mapsto \ell_F(\lambda xy)$ and $(x,y) \mapsto \ell_F(- d \lambda xy)$, where $\ell_F$ is the restriction
of $\ell$ to $F$. The determinant of $\ell(h)$ is
the product of the determinant of these two forms. Lemma \ref{det F} implies that the determinant of the first one is ${\rm N}_{F/K}(\lambda)D_1$, and
of the second one ${\rm N}_{F/K}(-d){\rm N}_{F/K}(\lambda)D_1$; this completes the proof of the proposition. 

\begin{notation}\label{hlambda}
For all $\lambda \in F^{\times}$ and all integers $n \geqslant 1$,  let $h_{n,\lambda}$
be the $n$-dimensional  diagonal hermitian form $h_{n,\lambda} = \langle \lambda, 1, \dots, 1 \rangle$ over $E$, and  set $q _{n,\lambda}= \ell(h_{n,\lambda})$.

\end{notation}

\medskip

\begin{prop}\label{general Hasse-Witt} Let $\lambda \in F^{\times}$ and let  $n \geqslant 1$ be an integer. We have
$$w_2(q_{n,\lambda})= w_2(q_{n,1}) + {\rm cor}_{F/K}(\lambda,d)$$
in ${\rm Br}_2(K)$. 

\end{prop}

To prove this proposition, we recall a theorem of Arason. If $k$ is a field of characteristic $\not = 2$, we denote by
$I(k)$  the ideal of even dimensional forms of the Witt ring $W(k)$, and by $e_2 : I^2(k) \to {\rm Br}_2(k)$ be the homomorphism
sending a $2$-fold Pfister form $\langle \langle a,b \rangle \rangle$ to the class of the quaternion algebra $(-a,-b)$ (see for instance \cite{A2}).

\begin{theo}\label{Arason} {\rm (Arason)} Let $s : F \to K$ be a non-trivial linear form. If $Q \in I^2(F)$, then 
$s(Q) \in I^2(K)$, and $$e_2(s(Q)) = {\rm cor}_{F/K}(e_2(Q))$$
in ${\rm Br}_2(K)$. 

\end{theo}

\noindent {\bf Proof.} See \cite{A2}, Theorem 1.21, or \cite{A1}, Satz 4.1 and Satz 4.18.

\begin{coro}\label{Arason coro}
If $Q \in I^2(F)$, then $$w_2(s(Q)) = {\rm cor}_{F/K}(w_2(Q))$$
in ${\rm Br}_2(K)$.

\end{coro}

\noindent{\bf Proof.}  \cite{Lam}, Proposition 3.20 implies that for $Q \in I^2(F)$, either $w_2(Q) = e_2(Q)$ or 
$w_2(Q) = e_2(Q) + (-1,-1)$. In the first case, there is nothing to prove; assume that $w_2(Q) = e_2(Q) + (-1,-1)$. We have $e_2(Q) = {\rm cor}_{F/K}(e_2(Q))$ by Theorem \ref{Arason}. Since 
$w_2(Q) = e_2(Q) + (-1,-1)$, we  have ${\rm cor}_{F/K}(w_2(Q)) = {\rm cor}_{F/K}(e_2(Q)) + (-1,(-1)^{[F:K]}) = {\rm cor}_{F/K}(e_2(Q)) 
[F:K](-1,-1)$.
On the other hand,  $w_2(s(Q)) = e_2(s(Q)) + [F:k](-1.-1)$ (see \cite{Lam}, Proposition 3.20), and this implies that  $w_2(s(Q)) = {\rm cor}_{F/K}(w_2(Q))$,
as claimed.

\medskip For all $n \in {\bf N}$, set $D_n = {\rm N}_{F/K}(-d)^n$;  Proposition \ref{det} implies that  $${\rm det}(q_{n,\lambda}) = D_n$$ for all $\lambda \in F^{\times}$.

\medskip
\noindent
{\bf Proof of Proposition \ref{general Hasse-Witt}.} 
Assume first that $n = 1$, and let $s : F \to K$ be the restriction of $\ell : E \to K$ to $F$. Note that $\ell(h_{1,1}) = s(\langle 2,-2d \rangle)$, and
$\ell(h_{1,\lambda}) = s(\langle 2,-2\lambda d \rangle)$. Set $q_1 = \langle 2,-2d \rangle$, $q_2 = \langle 2,-2\lambda d \rangle$, and
$Q = q_1 + q_2$.  We have $w_2(Q) = w_2(q_1) + w_2(q_2) + ({\rm det}(q_1), {\rm det}(q_2)) = (\lambda,d) + (-1,d)$, hence
${\rm cor}_{F/K}(w_2(Q)) = {\rm cor}_{F/K}(\lambda,d) + (-1,N_{F/K}(-d))$. On the other hand, we have
$w_2(s(Q)) = w_2(s(q_1) )+w_2(s(q_2)) + ({\rm det}(s(q_1)), {\rm det}(s(q_2))$. By Proposition \ref{det}, we have ${\rm det}(s(q_1)) = {\rm det}(s(q_2)) =
{\rm N}_{F/K}(-d)$, hence $w_2(s(Q)) = w_2(s(q_1) )+w_2(s(q_2)) + (-1,{\rm N}_{F/K}(-d))$. By Corollary \ref{Arason coro}, we have
$w_2(s(Q)) = {\rm cor}_{F/K}(w_2(Q))$, therefore $w_2(s(q_1) )+w_2(s(q_2)) = {\rm cor}_{F/K}(\lambda,d)$.
Since  $q_{1,1} = \ell(h_{1,1})= s(q_1)$ and  $q_{1,\lambda} = \ell(h_{1,\lambda})= s(q_2)$, this proves the proposition when $n = 1$. 
Asssume now that $n \geqslant 2$. We have
$$w_2(q_{n,\lambda}) = w_2(q_{1,\lambda}) + (D_1,D_{n-1}) + w_2(q_{n-1,1}),$$ and
$$w_2(q_{n,1}) = w_2(q_{1,1}) + (D_1,D_{n-1})  + w_2(q_{n-1,1}).$$

\smallskip

By the one-dimensional case, we have $w_2(q_{1,\lambda})  = w_2(q_{1,1}) +
{\rm cor}_{F/K}(\lambda,d)$, hence 
$$w_2(q_{\lambda})= w_2(q_{n,1}) + {\rm cor}_{F/K}(\lambda,d),$$ as claimed.

%\begin{prop}\label{general Hasse-Witt} For all $\lambda \in F^{\times}$, we have
%%$$w_2(q_{n,\lambda})= w_2(q_{n,1}) + {\rm cor}_{F/K}(\lambda,d)$$
%in ${\rm Br}_2(K)$. 

%\end{prop}

\bigskip
{\bf Local fields} 

\medskip Assume that $K$ is a local field.

\medskip
\begin{lemma}\label{herm} Two hermitian forms over $E$
having the same dimension and determinant are isomorphic.

\end{lemma}

\medskip
\noindent
{\bf Proof.} See for instance \cite {Sch}, 10.1.6, (ii).

\medskip
\begin{prop}\label{Hasse-Witt} Let $\lambda \in F^{\times}$, let $h$
be an $n$-dimensional  hermitian form of determinant $\lambda$ over $E$, and  set $q = \ell(h)$. We have
$$w_2(q) = w_2(q_{n,1}) + {\rm cor}_{F/K}(\lambda,d)$$
in ${\rm Br}_2(K)$. 
 
\end{prop}

\medskip
\noindent
{\bf Proof.}  If $K$ is a non-archimedean local field, then
by Lemma  \ref{herm} the
hermitian form $h$ is isomorphic to the diagonal form $h_{n,\lambda} = \langle \lambda, 1, \dots, 1 \rangle$, hence
the statement follows from Proposition \ref{general Hasse-Witt}.  
%[It can also be deduced  from \cite{M}, Theorem 2.7]. 

\medskip
Assume now that $K = {\bf R}$; then $F = K = {\bf R}$, and $E = {\bf C}$, hence $d = -1$.
We have $h = \langle \lambda_1, \dots, \lambda_n \rangle = \langle \lambda_1 \rangle \oplus
\dots \oplus \langle \lambda_n \rangle$ with $\lambda_i \in {\bf R}$. Note that ${\rm det}(\ell (\langle
\lambda_i \rangle) = 1$ for all $i$; hence we have
$w_2(q) = \underset{i \in I} \sum w_2(\ell (\langle \lambda_i \rangle))$ and 
$w_2(q_{n,1}) = \underset {i \in I} \sum w_2 (\ell (\langle 1 \rangle))$. 

\medskip
We have 
$\ell(\lambda_i) = \langle 1,1 \rangle$ if $\lambda_i > 0$, and $\ell(\lambda_i) = \langle -1,-1 \rangle$ if $\lambda_i < 0$; note that $w_2(\langle 1,1 \rangle ) = 0$ and $w_2(\langle -1,-1 \rangle) = 1$. 
This implies that $w_2(q)= w_2(q_{n,1})$ if and only if $\lambda_i$ is negative for an
even number of $i \in I$. On the other hand, ${\rm cor}_{F^v/K_v}(\lambda,d) = 0$ if and
only if  $\lambda_i$ is negative for an
even number of $i \in I$; this concludes the proof of the proposition.

\section{Obstruction}\label{obstruction}

The aim of this section is to describe an obstruction group 
in a general situation. The group depends on a finite set $I$, a set $V$ and,
for all $i,j \in I$, a subset $V_{i,j}$ of $V$; in the applications, the vanishing of the group 
detects the validity of the local-global principle (also called Hasse principle).

\medskip
To obtain a necessary and sufficient condition for the Hasse principle to hold, we need additional data; in the applications, it is 
provided
by  local solutions. For all $v \in V$, we give ourselves a set $\mathcal C^v$
having certain properties (see below for details); this set corresponds to the local data.

\bigskip
{\bf The basic setting}

\medskip
Let $I$ be a finite set, let $\sim$ be an
equivalence relation on $I$, and let $\overline I$ be
the set of equivalence classes. Let $C(I)$ be the set of maps $I \to {\bf Z}/2{\bf Z}$.  
Let $C_{\sim}(I)$ be
the subgroup of $C(I)$ consisting of the maps that are constant on the equivalence classes, and note that  $C_{\sim}(I) = C(\overline I)$; it is a finite
elementary abelian 2-group.

\medskip
Let
$\sha_{\sim}(I)$ be the quotient of $C_{\sim}(I)$ by the constant maps; equivalently, we can regard $\sha_{\sim}(I)$ as the
quotient of $C(\overline I)$ by the constant maps. 

\medskip We start by giving some examples that will be useful in the sequel. 

\begin{example}\label {sha indep}
This example will be used in \S \ref{rational obstruction} and \S  \ref{independent extensions}. We say that two finite extensions
$L_1$ and $L_2$ of a field $K$ are  {\it independent over $K$} if the tensor product $L_1 \otimes_K L_2$ is a field. 
Let $E = \underset{i \in I} \prod E_i$ be a product of finite field extensions $E_i$ of $K$, and let us consider the
equivalence relation $\sim$ generated by the elementary equivalence

\medskip
\centerline {$i \sim_e j \iff$ $E_i$ and $E_j$ are independent over $K$.}

\medskip
We denote by $\sha_{\rm indep}(E) = \sha_{\sim}(I)$ the quotient of $C_{\sim}(I)$ by the constant maps.

\end{example}

\begin{example}\label {Q and Z} For all $i \in I$, let $f_i \in {\bf Z}[X]$ be symmetric, irreducible polynomials of even degree;
if $p$ is a prime number, we denote by $f_i \ {\rm (mod \ p)}$ the polynomial in ${\bf F}_p[X]$ obtained by reduction
modulo $p$. 

\medskip Let $V_{i,j}$ be the set of prime numbers such that $f_i \ {\rm (mod \ p)}$ and $f_j \ {\rm (mod \ p)}$ have
a common irreducible, symmetric factor in  ${\bf F}_p[X]$. Let us consider the equivalence relation $\sim$
generated by the elementary equivalence
$$i \sim_e j \iff V_{i,j} \not = \varnothing,$$ and let $\sha^{\rm int}$ be the group $ \sha_{\sim}(I)$ defined above. This
group will play an important role in the ``integral" Hasse principle of \S \ref{Z}.

\medskip Let $V_i$ be the set of prime numbers such that  $f_i \ {\rm (mod \ p)}$ has an irreducible, symmetric
factor in ${\bf F}_p[X]$. Let us consider the equivalence relation $\approx$
generated by the elementary equivalence
$$i \approx_e j \iff V_i \cap V_j \not = \varnothing,$$ and let $\sha^{\rm rat}$ be the group $ \sha_{\approx}(I)$ defined above.  We have $V_{i,j} \subset V_i \cap V_j$, hence $i \sim j \implies i \approx j$; therefore  $\sha^{\rm rat}$ is
a subgroup of $\sha^{\rm int}$. The group $\sha^{\rm rat}$ plays a role in a ``rational" Hasse principle (see 
\S \ref{rational}). 

\end{example}

\medskip
{\bf Equivalence relation on $C(I)$}

\medskip As above, $I$ is a finite set and $\sim$ an equivalence relation on $I$. 
If $i,j \in I$, let $c_{i,j} \in C(I)$ be such that  $c_{i,j}(i) = c_{i,j}(j) = 1$ and $c_{i,j}(k) = 0$
if $k \not = i,j$. Let $(i,j) : C(I) \to C(I)$ be the map sending $c$ to $c + c_{i,j}$. 
We also denote by $\sim$ the equivalence relation on $C(I)$ generated by the elementary equivalence

$$c \sim_e c' \iff c = (i,j)(c')$$
for some $i, j \in I$ such that $i \sim j.$

\medskip
The following remark will be useful in the sequel
\begin{lemma}\label{indep} Let $a, b \in C(I)$ such that $a \sim b$. Then we have 

$$\underset{i \in I} \sum  \ c(i) a(i) = \underset{i \in I} \sum  \ c(i) b(i)$$ 
for all $c \in C_{\sim}(I)$.

\end{lemma}

\noindent
{\bf Proof.} 
We can assume that
$b = (i,j)(a)$ for some $i,j \in I$ such that $i \sim j$. Since $c \in C_{\sim}(I)$, we have
 $c(i) = c(j)$, hence $\underset{k \in I} \sum  \ c(k) a(k) = \underset{k \in I} \sum  \ c(k) b(k)$, as claimed.

%\bigskip
%Let
%$\sha_{\sim}(I)$ be the quotient of $C_{\sim}(I)$ by the constant maps; equivalently, we can regard $\sha_{\sim}(I)$ as the
%quotient of $C(\overline I)$ by the constant maps.

\bigskip
{\bf The sets $V$, $V_{i,j}$ and the associated equivalence relations}

\medskip

Let $V$ be a set, and for all $i,j \in I$, let $V_{i,j}$ be a subset of $V$. We take for
$\sim$ the 
equivalence relation generated by $$i \sim j \iff V_{i,j} \not = \varnothing,$$
and the equivalence relation $\sim$ on $C(I)$ generated by

$$c \sim c' \iff c = (i,j)(c')$$
for $V_{i,j} \not = \varnothing.$

\medskip
Let $\sha = \sha_{\sim}(I)$ be the group corresponding to the equivalence relation $\sim$ defined as above.

\bigskip
For all $v \in V$, we define equivalence relations $\sim_v$ on $I$ and $C(I)$, generated by
$$i \sim_v j \iff v \in V_{i,j},$$ 
and
 $$c \sim_v c' \iff c = (i,j)(c')$$
for some $i,j \in I$ with $v \in V_{i,j}.$

\bigskip

Let $V = V' \cup V''$, and assume that for all $i,j \in I$, if $V_{i,j} \not = \varnothing$, then $V_{i,j} \cap V' \not = \varnothing$. 

\medskip
For all $v \in V$, let $A^v \in {\bf Z}/2{\bf Z}$ be such that 

\medskip
(i) $A^v = 0$ for almost all $v \in V$, and
 $$\underset {v \in V} \sum A^v = 0.$$

\medskip
For all $v \in V$, we consider subsets $\mathcal C^v$ of $C(I)$ satisfying the conditions (ii) and (iii) below :

\medskip
{\rm (ii}) For all $a^v \in \mathcal C^v$, we have $$\underset{i \in I} \sum a^v(i) = A^v.$$

\medskip

{\rm (iii)} If $v \in V'$, then $\mathcal C^v$ is stable by the maps $(i,j)$ for $i \sim_v j$. 
%If $v \in V'$, then $\mathcal C^v$ is a
%$\sim_v$-equivalence class of $C(I)$. 

\bigskip
Let $\mathcal C$ be the set of the elements $(a^v)$ with $a^v \in \mathcal C^v$ such that
$a^v = 0$ for almost all $v \in V$.

\bigskip
We now prove some results that will be useful in the following sections.

\medskip
For all $a \in \mathcal C$ and $i \in I$, set $\Sigma_i(a) = \underset {v \in V} \sum a^v (i)$. 

\medskip
\noindent
\begin{prop}\label{Si} Let $a \in \mathcal C$, and let $i,j \in I$ be such that $i \sim j$. Then there 
exists $b \in \mathcal C$ such that $\Sigma_i(b) \not = \Sigma_i(a)$, $\Sigma_j(b) \not = \Sigma_j(a)$, and $\Sigma_k(b) = \Sigma_k(a)$ for
all $k \not = i,j$. 

\end{prop}

\medskip
\noindent
{\bf Proof.} Let $i_1,\dots,i_k \in I$ be such that $i = i_1$, $j = i_k$, and that $V_{i_s,i_{s+1}} \not = \varnothing$
for all $s = 1,\dots,k-1$; therefore $V_{i_s,i_{s+1}} \cap V' \not = \varnothing$ for all $s = 1,\dots,k-1$.

\medskip
 If $v \in V_{i_s,i_{s+1}} \cap V'$, then by condition (iii) the map $(i_s,i_{s+1})$ sends $\mathcal C^v$ to $\mathcal C^v$. 
Consider the map $\mathcal C \to \mathcal C$ that is equal to $(i_s,i_{s+1})$ on $C^v$, and is the identity on $\mathcal C^w$ if
$w \not = v$. Applying to $a \in \mathcal C$ the maps induced by $(i_1,i_2),\dots,(i_{k-1},i_k)$ successively
yields $b \in \mathcal C$ with the required properties. 

\bigskip
If $c \in \sha$ and  $a = (a^v) \in \mathcal C$, then by conditions (i) and (ii)  the sum
$$\underset {v \in V} \sum \  \underset{i \in I} \sum  \ c(i) a^v(i) $$ is well-defined. The following result is used in 
\S \ref{rational} and \S \ref{isometrieslattices}  to give necessary and sufficient conditions for some Hasse principles to hold. 

\bigskip

\begin{theo} \label{sum} Let $a = (a^v) \in \mathcal C$ such that 
$$\underset {v \in V} \sum \  \underset{i \in I} \sum  \ c(i) a^v(i) = 0$$ for all 
$c \in \sha$. Then there exists $b =(b^v) \in \mathcal C$ such that $$\underset {v \in V} \sum b^v (i) = 0$$
for all $i \in I$.

\end{theo}

\medskip
\noindent
{\bf Proof.} If $\Sigma_i(a) = 0$ for all $i \in I$, then we are done. Assume that $\Sigma_{i_0}(a) = 1$. 
We claim that there exists $i \in I$ with $i \not = i_0$ and $i \sim i_0$ such that $\Sigma_i (a) = 1$. 
In order to prove this claim,
let $c \in C(I)$ be such that $c(i_0) =1$ and that $c(i) = 0$ if $i \not = i_0$. Then we have 

$$\underset {v \in V} \sum \  \underset{i \in I} \sum  \ c(i) a^v(i) = \Sigma_{i_0}(a) = 1,$$ hence $c \not \in
\sha$. Therefore $c$ is not constant on the equivalence classes. This implies that
there exists $i_1 \in I$ with $i_1 \not = i_0$ such that $i_1 \sim i_0$. If $\Sigma_{i_1}(a) = 1$, we stop.
Otherwise, let $c \in C(I)$ be such that $c(i_0) = c(i_1) = 1$ and $c(i) = 0$ if $i \not = i_0,i_1$. 
We have $$\underset {v \in V} \sum \  \underset{i \in I} \sum  \ c(i) a^v(i) = \Sigma_{i_0}(a) = 1,$$ hence $c \not \in
\sha$. Therefore $c$ is not constant on the equivalence classes. This implies that
there exists $i_2 \in I$, $i_2 \not = i_0,i_1$, such that
either $i_2 \sim i_0$ or $i_2 \sim i_1$. Note that since $i_1 \sim i_0$, we have $i_2 \sim i_0$ in
both cases. Since $I$ is finite and that $\underset {v \in V} \sum \  \underset{i \in I} \sum  \ c(i) a^v(i) = 0$,
we eventually get $i_r \in I$ with $i_r \not = i_0$, $i_r \sim i_0$, and $\Sigma_{i_r}(a) = 1$. By Proposition
\ref{Si}, there exists $b \in \mathcal C$ such that $\Sigma_{i_0}(b) = \Sigma_{i_r} (b) = 0$. Continue inductively until
the theorem is proved. 

\bigskip
Under some additional hypothesis on $\mathcal C^v$ for $v \in V''$, the necessary and sufficient condition
for the Hasse principle can be given by the vanishing of a homomorphism :

\bigskip

{\bf The homomorphism}

\medskip
Let us make the additional assumption that 

\medskip
(iv) For all $v \in V$, the set $\mathcal C^v$ is a subset of an
$\sim_v$-equivalence class of $C(I)$.

\medskip
Let $a \in \mathcal C$. We define a homomorphism $\alpha = \alpha_a : \sha \to {\bf Z}/2{\bf Z}$ as
follows. For all $c \in \sha$, we set

$$\alpha(c) = \underset {v \in V} \sum \  \underset{i \in I} \sum  \ c(i) a^v(i).$$
Note that this is well-defined, since by conditions (i) and (ii) $$\underset {v \in V} \sum \underset{i \in I} \sum a^v(i)  = 
\underset {v \in V} \sum A^v = 0.$$

\begin{prop}\label{independent}

 The homomorphism $\alpha$ is independent of the choice of $a \in \mathcal C$.
 
 \end{prop}
 
 \noindent
 {\bf Proof.} Let $a = (a^v), b = (b^v) \in \mathcal C$, and let us show that $\alpha_a = \alpha_b$.
 Let $v \in V$. Since $a^v, b^v \in \mathcal C^v$, we have $a^v \sim_v b^v$, hence $a^v \sim b^v$.
 Therefore by Lemma \ref{indep},  we have 
 $$\underset{i \in I} \sum  \ c(i) a^v(i) = \underset{i \in I} \sum  \ c(i) b^v(i).$$ This holds for
 all $v \in V$, hence $\alpha_a = \alpha_b$.
 
 \begin{coro} \label{ihomo} Let $a \in \mathcal C$, and assume that $\alpha_{a} = 0$. Then there
exists $b =(b^v) \in \mathcal C$ such that $$\underset {v \in V} \sum b^v (i) = 0$$
for all $i \in I$.

\end{coro}

\medskip
\noindent
{\bf Proof.} This follows from Theorem \ref{sum}.

\section{Obstruction group - the rational case}\label{rational obstruction}

The aim of this section is to associate an ``obstruction group" to certain algebras with involution; this group
will play an important role in the Hasse principle results of the following sections.

\medskip
Assume that $K$ is a global field of characteristic $\not =2$, and let $I$ be a finite set. For all $i \in I$,
let $E_i$ be a finite degree extension of $K$, and let $\sigma_i : E_i \to E_i$ be a non-trivial 
involution; let $F_i$ be the fixed field of $\sigma_i$, and let $d_i \in F_i^{\times}$ be such that $E_i = F_i(\sqrt d_i)$. 
Set $F = \underset{i \in I} \prod F_i$ and
$E = \underset{i \in I} \prod E_i$; let $\sigma : E \to E$ be the involution such that the
restriction of $\sigma$ to $E_i$ is equal to $\sigma_i$.  

%\medskip
%The aim of this section is to associate to $(E,\sigma)$ a finite group $\sha_E = \sha_{E,\sigma}$
%that will be useful in \S \ref{rational}.

\medskip We denote by $V_K$ the set of places of $K$. 

\begin{notation}\label{Vi}
For all $i \in I$, let $V_i$ be set of places $v \in V_K$ such that there exists a place
of $F_i$ above $v$ which is inert or ramified in $E_i$. 

\end{notation}

Let $\sha = \sha_E$
be the group constructed in \S \ref{obstruction}
 using the data $I$
 % $V = V_K$ 
 and $V_{i,j} = V_i \cap V_j$; recall that the equivalence relation on $\sim$ $I$ is generated by the
 elementary equivalence $$i \sim_e j \iff V_i \cap V_j \not = \varnothing,$$
 that $C_{\sim}(I)$ is
the subgroup of $C(I)$ consisting of the maps $I \to {\bf Z}/2{\bf Z}$ that are constant on the equivalence classes, and 
that $\sha_E$ is the quotient of $C_{\sim}(I)$ by
the constant maps. 

\bigskip In the remaining of the section, we give some examples and some results that will be used in the next sections. 
We start by giving some examples in which the obstruction group is trivial.
 
 \medskip
 \begin{example}\label{CM} Assume that there exists a real place $v$ of $K$ such that for all $i \in I$, there exists a real place of $F_i$ above $v$ which
 extends to a complex place of $E_i$. Then $\sha_E = 0$. Indeed, $v \in V_i$ for all $i \in I$, hence $V_i \cap V_j \not =
 \varnothing$ for all $i, j \in I$. Therefore all the elements of $I$ are equivalent, and this implies that $\sha_E = 0$.
 
 \medskip
 In particular, if $K = {\bf Q}$, and if for all $i \in I$, the field $E_i$ is a CM field (that is, $E_i$ is totally complex 
and  $F_i$ is totally real), then $\sha_E = 0$.

  %Assume that $K = {\bf Q}$, and that for all $i \in I$, the field $E_i$ is CM; in other words, $E_i$ is totally complex 
% and $F_i$ is totally real. Then $\sha_E = 0$. Indeed, the infinite place of $\bf Q$ belongs to $V_i$ for all $i \in I$, hence $V_i \cap V_j \not =
% \varnothing$ for all $i, j \in I$. Therefore all the elements of $I$
% are equivalent, and this implies that $\sha_E = 0$.
 
  \end{example}

Recall  that two finite extensions
$K_1$ and $K_2$ of $K$ are  {\it independent} if the tensor product $K_1 \otimes_K K_2$ is a field. If $K_1$ and $K_2$ are
subfields of a field extension $\Omega$ of $K$, then this means that $K_1$ and $K_2$ are linearly disjoint. 

 \begin{prop}\label{2} Assume that $E_i$ and $E_j$ are independent field extensions
of $K$. Then $V_{i,j} \not = \varnothing$.

\end{prop}

\noindent
{\bf Proof.} Let $\Omega/K$ be a Galois extension containing $E_i$
and $E_j$, and set $G = {\rm Gal}(\Omega/L)$. Let $H_i \subset G_i$ and $H_j \subset G_j$  be subgroups of $G$ such that
we have $E_i = \Omega^{H_i}$, $E_j = \Omega^{H_j}$ and $F_i = \Omega^{G_i}$, $F_j = \Omega^{G_j}$. Since $E_i/F_i$ and $E_J/F_J$ are quadratic extensions,  the subgroup
$H_i$ is of index 2 in $G_i$, and $H_j$ is of index 2 in $G_j$. By hypothesis, $E_i$ and $E_j$ are linearly disjoint over $K$, therefore
$[G : H_i \cap H_j ] = [G : H_i] [G : H_j]$. Note that $F_i$ and $F_j$ are also linearly disjoint over $K$, hence
$[G : G_i \cap G_j] = [G : G_i][G : G_j]$.
%We claim that $[G_1 \cap G_2: H_1 \cap H_2] = 4$. Indeed,
%$$[G : G_1 \cap G_2] \geqslant [G : G_1][G : G_2] = {1 \over 4}[G:H_1][G:H_2].$$
%On the other hand, $[G:H_1 \cap H_2] = [G:G_1\cap G_2] [G_1 \cap G_2: H_1 \cap H_2]$, hence
%$[G_1 \cap G_2: H_1 \cap H_2] \geqslant 4$. 
This implies that $[G_i \cap G_j: H_i \cap H_j] = 4$, hence 
the quotient $G_i \cap G_j/H_i \cap H_j$ is an elementary abelian group of order 4. 

\medskip
The field $\Omega$ contains the composite fields $F_i F_j$ and $E_i E_j$. By the above argument, the
extension $E_iE_j/F_i F_j$ is biquadratic. Hence there exists a place of $F_i F_j$ that is inert in both
$E_i F_j$ and $E_j F_i$. Therefore there exists a place $v$ of $K$ and places $w_i$ of $F_i$ and $w_j$ of $F_j$  above $v$ 
such that $w_i$ is inert in $E_i$, and $w_j$ is inert in $E_j$, hence $v \in V_i \cap V_j$; by definition, this
implies that $v \in V_{i,j}$. 

\medskip Recall from \S \ref{obstruction}, Example \ref{sha indep}, the equivalence relation $\approx$ on $I$ generated by the elementary equivalence

\medskip
\centerline {$i \approx_e j \iff$ $E_i$ and $E_j$ are independent over $K$,}

\medskip
\noindent
and that we denote by $\sha_{\rm indep}(E) = \sha_{\approx}(I)$  the associated obstruction group (cf. \S \ref{obstruction}).
 Proposition \ref{2} implies that the identity $C(I) \to C(I)$ induces a surjection $\sha_{\rm indep}(E) \to \sha_E$. 
Hence we have

\begin{coro}\label {0 and 0}
Assume that $\sha_{\rm indep}(E) =0$; then $ \sha_E = 0$.

\end{coro}

On the other hand, one can show that all elementary abelian 2-groups occur as $\sha_E$ for some $(E,\sigma)$. 

\medskip
Finally, we show a result that will be used in \S \ref {odd degree} to prove an odd degree descent property. Let
$K'/K$ be a finite extension of odd degree. We define a homomorphism
$$\sha_E \to  \sha_{E \otimes_K K'}$$ as follows. For all $i \in I$, let $F_i \otimes_K K' = \underset{j \in I(i)} \prod F'_{i,j}$,
where $F'_{i,j}$ is a field extension of $K'$. Set

\medskip

\centerline {$I' = \{ (i,j) \ | \ i \in I, \ j \in I(i)$ and the image of $d_i$ in $F'_{i,j}$ is not a square\}.}

\medskip Let $\pi : I' \to I$ be the map sending $(i,j)$ to $i$; this map induces a homomorphism $\pi ' : C(I) \to C(I')$. 

\begin{prop}\label {odd sha} The map $\pi ' : C(I) \to C(I')$ 
induces an injective homomorphism 
 $$\pi' : \sha_E \to  \sha_{E \otimes_K K'}.$$

\end{prop}

\noindent
{\bf Proof.} Let us show that $\pi'$ sends $C_{\sim}(I)$ to $C_{\sim}(I')$. Let $(i_1,j_1), (i_2,j_2) \in I'$ be
such that $(i_1,j_1) \sim_e (i_2,j_2)$; by definition, this means that $V_{(i_1,j_1), (i_2,j_2)} \not = \varnothing$. Let
$v' \in V_{(i_1,j_1), (i_2,j_2)} $, and let $v \in V_K$ be the restriction of $v'$ to $K$. Then $v \in V_{i_1,i_2}$,
hence $i_1 \sim_e i_2$. This shows that if $c \in C_{\sim}(I)$, then  $\pi'(c) \in C_{\sim}(I')$, hence we have a well-defined
homomorphism  $C_{\sim}(I) \to C_{\sim}(I')$. This in turn induces a homomorphism 
$\pi' : \sha_E \to  \sha_{E \otimes_K K'},$ and it is clearly injective.

\section{Twisting groups and equivalence relations}\label{equiv}

In this section, we introduce some notation that will be used throughout the paper; we also
prove some results concerning equivalence relations defined on the twisting groups.

%\medskip
%Let $K$ be a global field of characteristic $\not =2$, and let $I$ be a finite set. For all $i \in I$,
%let $E_i$ be a finite degree extension of $K$, and let $\sigma_i : E_i \to E_i$ be a non-trivial 
%involution; let $F_i$ be the fixed field of $\sigma_i$. 
%Let $d_i \in F_i$ be such that $E_i =
%%F_i(\sqrt {d_i})$. 
%Set $F = \underset{i \in I} \prod F_i$ and
%$E = \underset{i \in I} \prod E_i$; let $\sigma : E \to E$ be the involution such that the
%restriction of $\sigma$ to $E_i$ is equal to $\sigma_i$.  Set $d = (d_i)$. 

\medskip We keep the notation of \S \ref{rational obstruction}.
 If $v \in V_K$, we denote by $K_v$ the completion of $K$ at $v$, and set $E_i^v = E_i \otimes_K K_v$,
 $F_i^v = F_i \otimes_K K_v$.

\begin{notation}\label{twist}
For all $i \in I$, let $S_i$ be the set of places $w$ of $F_i$ such that $E_i \otimes_{F_i} (F_i)_w$
is a field. If $w \in S_i$, set $E_i^w = E_i \otimes_{F_i} (F_i)_w$, and let

$$\theta_i^w : T(E_i^w,\sigma_i) \to {\bf Z}/2{\bf Z}$$

\medskip
\noindent be the isomorphism defined in (\ref{e1}). If $v \in V_K$, we denote by $S_i^v$ the set of places of $S_i$ above $v$. 
For
all $v \in V_K$, set 

$$T(E_i^v,\sigma_i) = \underset {w \in S_i^v} \prod T(E_i^w,\sigma_i),$$ and
$$T(E^v,\sigma) = \underset {i \in I} \prod T(E_i^v,\sigma_i).$$

\medskip
Recall that $C(I)$ is the set of maps $I \to {\bf Z}/2 {\bf Z}$. For all $\lambda^v = (\lambda_i^w) \in T(E^v,\sigma)$, we define $a(\lambda^v) \in C(I)$ by
setting $$a(\lambda^v)(i) = \underset {w \in S_i^v} \sum \theta_i^w(\lambda_i^w).$$ 

\medskip 
Let $\tilde C(I)$ be the set of maps $I \to {\bf Q}/{\bf Z}$. If $\lambda^v = (\lambda_i^w) \in T(E^v,\sigma)$, 
let  $\tilde a(\lambda^v) \in \tilde C(I)$ be defined by
$$\tilde a(\lambda^v)(i) = {\rm inv}_v({\rm cor}_{F^v/K_v}(\lambda_i^v,d_i)).$$

%\medskip
%Recall from  (\ref{e1}) that the injection ${\bf Z}/2{\bf Z} \to {\bf Q}/{\bf Z}$ sends $\theta_i^w(\lambda_i^w)$
%to ${\rm inv}_w (\lambda_i^w,d_i)$.

\medskip
For all $v \in V_K$, let $C_{T^v}(I)$ be
the subset of $C(I)$ consisting of the maps $a(\lambda^v)$ for $\lambda^v \in T(E^v,\sigma)$. 
Similarly, let us denote by $\tilde C_{T^v}(I)$ be
the subset of $C(I)$ consisting of the maps $\tilde a(\lambda^v)$ for $\lambda^v \in T(E^v,\sigma)$. 

\end{notation}

\begin {lemma}\label{tilde} Let $v \in V_K$. Sending $a(\lambda^v)$ to $\tilde a(\lambda^v)$ gives
rise to a bijection $$C_{T^v}(I) \to  \tilde C_{T^v}(I).$$

\end{lemma} 

\medskip
\noindent
{\bf Proof.} Let $\iota : {\bf Z}/2{\bf Z} \to {\bf Q}/{\bf Z}$ be the canonical injection, and let us
also denote by $\iota$ the induced injection $C(I) \to \tilde C(I)$. Let us show that
for all $\lambda^v \in T(E^v,\sigma)$, we have $\tilde a(\lambda^v) = \iota \circ a(\lambda^v)$.

\medskip
Recall from  (\ref{e1}) that for all $\lambda^v = (\lambda_i^w) \in T(E^v,\sigma)$, the injection $\iota : {\bf Z}/2{\bf Z} \to {\bf Q}/{\bf Z}$ sends $\theta_i^w(\lambda_i^w)$
to ${\rm inv}_w (\lambda_i^w,d_i)$.
Since $a(\lambda^v)(i) = \underset {w \in S_i^v} \sum \theta_i^w(\lambda_i^v)$, the
injection $\iota : C(I) \to \tilde C(I)$ sends $a(\lambda^v)(i)$ to $\underset {w \in S_i^v} \sum
{\rm inv}_w (\lambda_i^w,d_i)$. By (\ref{e2}), this is equal to $\underset {w \in S_i^v} \sum
{\rm inv}_v ( {\rm cor}_{(F_{i,w}/K_v} (\lambda_i^w,d_i))$, which in turn is equal to 
${\rm inv}_v({\rm cor}_{F^v/K_v}(\lambda_i^v,d_i))$; hence we have $\tilde a(\lambda^v) = \iota \circ a(\lambda^v)$, as claimed.

\bigskip
Note
that $C(I)$ is a group, and that $C_{T^v}(I)$ is a subgroup of $C(I)$. 

\medskip
Recall that for all $i \in I$,  we denote by $V_i$ the set of places $v \in V_K$ such that there exists a place
of $F_i$ above $v$ which is inert or ramified in $E_i$. Note that 

$$v \in V_i \iff S_i^v \not = \varnothing.$$

\begin{lemma}\label{CT invariant} The subgroup $C_{T^v}(I)$ of $C(I)$ is stable by the
maps $(i,j)$ for $v \in V_i \cap V_j$. 

\end{lemma}

\medskip
\noindent
{\bf Proof.} Let  $v \in V_i \cap V_j$, and let $\lambda^v = (\lambda_i^w) \in T(E^v,\sigma)$. Let $(+1) : {\bf Z}/2{\bf Z}$ be the
map sending $n$ to $n+1$. Let $w_i \in S_i^v$, and $w_j \in S_j^v$. Set
$$\mu_i^{w_i} = [(\theta_i^{w_i})^{-1} \circ (+1) \circ \theta_j^{w_j}](\lambda_j^{w_j}),$$ 
$$\mu_j^{w_j} = [(\theta_j^{w_j})^{-1} \circ (+1) \circ \theta_i^{w_i}](\lambda_i^{w_i}),$$ 
and $\mu_r^w = \lambda_r^w$ for all $(r,w) \not = (i,w_i), (j,w_j)$. Then 
$\mu^v = (\mu_i^w) \in T(E^v,\sigma)$,
and $$a(\mu^v) = (i,j)a(\lambda^v).$$

\bigskip
Let $v \in V_K$, and let $A^v \in {\bf Z}/2{\bf Z}$.  

\begin{notation}\label{A} Let $\mathcal L_{A^v}$ be the 
set of $\lambda^v = (\lambda_i^w) \in T(E^v,\sigma)$ such that 
$$A^v =  \underset{i \in I} \sum a(\lambda^v)(i),$$  and let $\mathcal C_{A^v}$ be
the set of $a(\lambda^v)$ for $\lambda^v \in \mathcal L_{A^v}$. 

\medskip
Let $\equiv_v$ be the  equivalence relation on  $C(I)$ generated by
$$c \equiv_v c' \iff c = (i,j)(c')$$
for some $i,j \in I$ with $v \in V_i \cap V_j.$  

\end{notation}

\begin{prop}\label{CA} The set $\mathcal C_{A^v}$ is a $\equiv_v$-equivalence class of $C(I)$.

\end{prop}

\medskip
\noindent
{\bf Proof.} 
Let us first show that $\mathcal C_{A^v}$ is stable by the maps $(i,j)$ for $v \in V_i \cap V_j$.
Let $\lambda^v \in \mathcal L_{A^v}$. By Lemma \ref{CT invariant}, there exists 
$\mu^v \in T(E^v,\sigma)$ such that $a(\mu^v) = (i,j)a(\lambda^v).$ Note that this implies
that $$\underset{r \in I} \sum a(\mu^v)(r) = \underset{r \in I} \sum a(\lambda^v)(r) = A^v,$$
hence $\mu^v \in \mathcal L_{A^v}$.

\medskip
Let us show that if $\lambda^v, \mu^v \in \mathcal L_{A^v}$, then $a(\lambda^v) \equiv_v a(\mu^v)$. 
If $a(\lambda^v)(i) = a(\mu^v)(i)$ for all $i \in I$, there is nothing to prove. Suppose that there
exists $i \in I$ such that 
$a(\lambda^v)(i) \not = a(\mu^v)(i)$; then $S_i^v \not = \varnothing$ , hence  $v \in V_i$. By hypothesis, we have
 $$\underset{r \in I} \sum a(\mu^v)(r) = \underset{r \in I} \sum a(\lambda^v)(r) = A^v,$$
therefore there
exists $j \in I$, $j \not = i$, such that $a(\lambda^v)(j) \not = a(\mu^v)(j)$. This implies that $v \in V_j$, hence
$v \in V_i \cap V_j$. The map $(i,j)a(\lambda^v)$  differs from $a(\mu^v)$ in
less elements than $a(\lambda^v)$. 
Since $I$ is a finite set, continuing this way we see that $a(\lambda^v) \equiv_v 
a(\mu^v)$.

\section{Local data - the rational case}\label{quadr}

\medskip
We keep the notation of \S \ref{equiv}; in particular, $K$ is a global field of characteristic $\not = 2$.
Let $M$ be an $A$-module satisfying the hypotheses of \S \ref {factors}
with $\mathcal A_M = (E_i)_{i \in I}$. Recall that  $E_i/K$ is a  finite  field extension of $K$ and that 
 $\sigma(E_i) = E_i$ for all $i \in I$.  In addition, assume that  for all $i \in I$, {\it the restriction of $\sigma$ to $E_i$ is non-trivial}. The fixed
field of this involution is denoted by $F_i$; hence $E_i/F_i$ is a quadratic extension, and $E_i = F_i(\sqrt d_i)$ for some $d_i \in F_i^{\times}$. 
Set $F = \underset{i \in I} \prod F_i$, and $E = \underset {i \in I} \prod  E_i$.

%\medskip

%{\bf Local data}

\medskip Let $v \in V_K$, and let $q$ be a quadratic form over $K_v$ which is compatible with the module $M \otimes_K K_v$. We now
associate a subset of $T(E^v,\sigma)$  to $q$ and $M$.

\medskip We have $M \simeq \underset {i \in I}\oplus M_i$ with $M_i$ isotypic,  such that $M_i$ is a finite dimensional $E_i$-vector space. By Proposition \ref{decomposition} there exist quadratic forms $q_i^v$ over $K_v$ compatible with $M_i\otimes_K K_v$ such that $q 
 \simeq \underset {i \in I} \oplus q_i^v$. Set $n_i = {\rm dim}_{E_i}(M_i)$.
By Proposition \ref{transfer} there exist $n_i$-dimensional hermitian forms $h^v_i$ over $E^v_i$ 
such that $q_i^v \simeq  \ell_i(h_i^v)$.
Set $\lambda_i^v = {\rm det}(h_i^v)$, and let us denote by $\mathcal {LR}_i^v$ the set of $\lambda_i^v
\in T(E_i^v,\sigma_i)$ that are obtained this way. Let $\lambda^v = (\lambda_i^v)$, and let us denote
by $\mathcal {LR}^v$ the set of $\lambda^v$ with $\lambda_i^v \in \mathcal {LR}_i^v$.
Let $\mathcal {CR}^v$ be the set of $a^v = a(\lambda^v)  \in C_{T^v}(I)$ with $\lambda^v \in \mathcal {LR}^v$ (cf. \S \ref{equiv}). 
Let $d = (d_i)$, and set
$q_{n,1} = \underset  {i \in I} \oplus q_{n_i,1}.$
%By  Proposition \ref{Hasse-Witt}. we have $w_2(q_i^v) = w_2(q_{n_i,1}) + {\rm cor}_{F/K}(\lambda_i^v,d_i)$ in ${\rm Br}_2(K_v)$ for all $i \in I$ and
%all $v \in V_K$.

\medskip
\begin{prop}\label{atilda infinite places} Let  $\lambda^v \in \mathcal{LR}^v$. Then we have
$$w_2(q) = w_2(q_{n,1}) + {\rm cor}_{F/K}(\lambda^v,d)$$ in ${\rm Br}_2(K_v).$

\end{prop} 

\noindent
{\bf Proof.} By  Proposition \ref{Hasse-Witt}, we have $w_2(q_i^v) = w_2(q_{n_i,1}) + {\rm cor}_{F/K}(\lambda_i^v,d_i)$ in ${\rm Br}_2(K_v)$ for all $i \in I$ and
all $v \in V_K$. Using \cite{Sch}, Chapter 2, Lemma 12.6, this
 gives the desired result, noting that ${\rm det}(q_i^v) = {\rm det}(q_{n_i,1})$ for all $i \in I$ and $v \in V_K$. 

\bigskip
 Set $\tilde A^v = {\rm inv}_v(w_2(q) + w_2(q_{n,1}))$. Recall that for $\lambda^v = (\lambda_i^w) \in T(E^v,\sigma)$, 
we defined $\tilde a(\lambda^v) \in \tilde C(I)$ by
$\tilde a(\lambda^v)(i) = {\rm inv}_v({\rm cor}_{F^v/K_v}(\lambda_i^v,d_i)).$ Let  $\tilde {\mathcal L}_{\tilde A^v}$ be the 
set of $\lambda^v = (\lambda_i^w) \in T(E^v,\sigma)$ such that 
$$\tilde A^v =  \underset{i \in I} \sum \tilde a(\lambda^v)(i),$$ and let
$\tilde {\mathcal C}_{\tilde A}$ the set of $\tilde a(\lambda^v)$ for $\lambda^v \in \tilde{ \mathcal L}_{\tilde A^v}$.

\medskip
Let $\iota : {\bf Z}/2{\bf Z} \to {\bf Q}/{\bf Z}$ be the canonical injection, and
let $A \in {\bf Z}/2{\bf Z}$ be such that $\iota (A) = \tilde A$. Recall that ${\mathcal L}_{ A^v}$ is the 
set of $\lambda^v = (\lambda_i^v) \in T(E^v,\sigma)$ such that 
$$A^v =  \underset{i \in I} \sum a(\lambda^v)(i),$$  and that $\mathcal C_{A^v}$ is
the set of $a(\lambda^v)$ for $\lambda^v \in \mathcal L_{A^v}$.

\medskip
\begin{prop}\label{A infinite} The set $\mathcal {CR}^v$ is contained in ${\mathcal C}_{A^v}$. 

\end{prop}

\medskip
\noindent
{\bf Proof.} Let us first note that if $\lambda^v \in \mathcal{LR}^v$, then  $\tilde a (\lambda^v)$ belongs to $ \tilde {\mathcal C}_{\tilde A^v}$. Indeed, since
${\rm inv}_v({\rm cor}_{F^v/K_v}(\lambda^v,d)) = \underset {i \in I} \sum {\rm inv}_v({\rm cor}_{F^v/K_v}(\lambda_i^v,d_i)),$ this
is a consequence of 
Proposition \ref{atilda infinite places}; the proposition now follows from Lemma  \ref{tilde}.

\medskip
\begin{coro}\label{rat utile} The set $\mathcal {CR}^v$ is contained in an $\equiv_v$-equivalence class of $C(I)$.

\end{coro}

\medskip
\noindent
{\bf Proof.}  This follows from Propositions \ref{A infinite} and \ref{CA}.

\bigskip

{\bf Local data - finite places}

\medskip
In the case of finite places, we have more precise information : as we will see, the sets 
$\mathcal {CR}^v$ and $\mathcal {C}_{A^v}$
coincide, and hence $\mathcal {CR}^v$ is an $\equiv_v$-equivalence class.

\medskip If $\lambda^v = (\lambda_i^v) \in T(E^v,\sigma)$ for some $v \in V_K$, set $q_{n,\lambda^v} = \underset  {i \in I} \oplus q_{n_i,\lambda_i^v}.$

\begin{prop} \label{finiteplace}  Suppose that $v$ is a finite place. Then ${\mathcal LR}^v$ is equal to the
 set of $\lambda^v \in T(E^v,\sigma)$ such that $q \simeq q_{n,\lambda^v}$.

\end{prop} 

\medskip
\noindent
{\bf Proof.} It is clear that if $\lambda^v \in T(E^v,\sigma)$ such that $q \simeq q_{n,\lambda^v}$, then $\lambda^v \in {\mathcal LR}^v$.
Conversely, let $\lambda^v \in {\mathcal LR}^v$, and let $h^v$ be a hermitan form over $E^v$ such that
$q \simeq \ell(h^v)$ and ${\rm det}(h^v) = \lambda^v$. By Lemma \ref{herm} we have $h^v \simeq h_{n,\lambda^v}$, hence $q \simeq q_{n,\lambda^v}$.

\medskip
\begin{prop}\label{atilda} Suppose that $v$ is a finite place. Let $\lambda^v  \in T(E^v,\sigma)$, Then $\lambda^v \in \mathcal{LR}^v$ if and only if
$$w_2(q) = w_2(q_{n,1}) + {\rm cor}_{F/K}(\lambda^v,d)$$ in ${\rm Br}_2(K_v).$

\end{prop}

\medskip
\noindent
{\bf Proof.} If $\lambda^v \in \mathcal{LR}^v$, then by Proposition \ref{atilda infinite places} we have $w_2(q) = w_2(q_{n,1}) + {\rm cor}_{F/K}(\lambda^v,d)$.  Let us show the converse.   We have
$w_2(q_{n,\lambda^v}) = w_2(q_{n,1}) + {\rm cor}_{F/K}(\lambda^v,d)$ in ${\rm Br}_2(K_v)$, therefore
$w_2(q_{n,\lambda^v}) = w_2(q)$ in ${\rm Br}_2(K_v).$ By Proposition \ref{det}, we have ${\rm det}(q_{n,\lambda^v}) = {\rm det}(q)$ in 
$K_v^{\times}/K_v^{\times 2}$, therefore the quadratic forms $q$ and $q_{n,\lambda}$ have the
same dimension, determinant and Hasse-Witt invariant; hence they are isomorphic over $K_v$. This implies that 
$\lambda^v \in \mathcal{LR}^v$.

\begin{prop}\label{A} Suppose that $v$ is a finite place; then  $\mathcal {CR}^v = {\mathcal C}_{A^v}$. 

\end{prop}

\medskip
\noindent
{\bf Proof.} By Proposition \ref{atilda}, we have $\mathcal {CR}^v = \tilde {\mathcal C}_{\tilde A^v}$; hence the statement follows from
Lemma  \ref{tilde}.

\medskip
\begin{coro}\label{utile} Suppose that $v$ is a finite place. Then the set $\mathcal {CR}^v$ is an $\equiv_v$-equivalence class of $C(I)$.

\end{coro}

\medskip
\noindent
{\bf Proof.}  This follows from Propositions \ref{A} and \ref{CA}.

\bigskip
{\bf Local data - real places}

\medskip
Suppose that  $v \in V_K$ is a real place. We say that a quadratic form $q$ has {\it maximal signature at $v$} (with respect to $M$) if the signature
of $q$  at $v$ is equal to the signature of $q_{n,1}$ or of $- q_{n,1}$  at $v$.

\begin{prop}\label{real ij} Assume that $q$ does not have maximal signature at $v$, and let $a \in {\mathcal CR}^v$. Then there
exist $i, j \in I$ with
$v \in V_i \cap V_j$ such that $(i,j)a$ belongs to  $ {\mathcal CR}^v$.

\end{prop}

\noindent
{\bf Proof.} Let $\lambda^v = (\lambda^w_r)  \in \mathcal {LR}^v$ such that $a = a(\lambda^v)$. Since the signature of $q$ at $v$ is not maximal, there exist $i,j \in I$ with $v \in V_i \cap V_j$ and 
$w_i \in V_{F_i}$,  $w_j \in V_{F_j}$ above $v$ such that $\lambda^{w_i}_i$ and $\lambda^{w_j}_j$ have opposite signs ($\lambda^{w_i}_i > 0$ and
$\lambda^{w_j}_j < 0$, or vice versa). For all $r \in I$ and $w \in V_{F_r}$ above $v$, let $\mu^w_r \in T(E^w_r,\sigma)$ be such that $\mu^w_r = \lambda^w_r$
if $(w,r) \not = (w_i,i)$ and $(w,r) \not = (w_j,j)$, that $\mu^{w_i}_i = - \lambda^{w_i}_i$ and $\mu^{w_j}_j = - \lambda^{w_j}_j$. We have $\mu^v = (\mu^w_r)
\in \mathcal {CR}^v$, and $a(\mu^v) = (i,j) a(\lambda^v)$.

%$\approx$

\section{Local-global problem - the rational case}\label{rational}

 We keep the notation of the previous sections; in particular, $M$ is an $A$-module satisfying the hypotheses
of \S \ref{quadr}, and $\mathcal A_M = (E_i)_{i \in I}$.  If $(V,q$) is a quadratic form over $K$, we consider the following local and global conditions :

\medskip
(L 1) {\it For all $v \in V_K$, the quadratic form $(V,q)\otimes_KK_v$ is compatible with the module $M \otimes_KK_v$.}

\medskip
(G 1) {\it The quadratic form $(V,q$) is compatible with the module $M$.}

\begin{prop}\label{cond} Assume that condition {\rm (L 1)} is satisfied. Then condition {\rm (G 1)}
holds  if and only if there exists $\lambda = (\lambda_i) \in T(E,\sigma)$ such
that $\lambda \in \mathcal {LR}^v$ for all $v \in V_K$. 

\end{prop} 

\medskip
\noindent
{\bf Proof.} We have $M \simeq \underset {i \in I}\oplus M_i$ with $M_i$ isotypic,  such that $M_i$ is a finite dimensional $E_i$-vector space.
%Let $M_i$ be as in \S  \ref{quadr}. 
If condition (G~1) holds,  then $q \simeq \underset{i \in I} \oplus q_i$, where the quadratic form $q_i$ is compatible with $M_i$
(cf. Proposition \ref {decomposition}).
For
all $i \in I$, let $h_i$ be a hermitian form over $E_i$ such that $q_i = \ell_i(h_i)$, 
set $\lambda_i = {\rm det}(h_i)$, and $\lambda = (\lambda_i)$; we have $\lambda \in \mathcal {LR}^v$
for all $v \in V_K$. 

\medskip 
Let us prove the converse. 
Let $\lambda = (\lambda_i) \in T(E,\sigma)$ such
that $\lambda \in \mathcal {LR}^v$ for all $v \in V_K$, and let $h_i$ be a hermitian form of dimension $n_i$ over
$E_i$ such that ${\rm det}(h_i) = \lambda_i$ and that ${\rm sign}_v(h_i) = {\rm sign}_v (h_i^v)$ for
all real places $v \in V_K$; for the existence of such a hermitian form, see for instance \cite {Sch}, 10.6.9.
Let $q_i = \ell_i(h_i)$ and $q' = \underset{i \in I} \oplus q_i$. The quadratic form $q'$ is compatible with $M$ by construction; its
dimension, signatures and determinant are the same as those of $q$ [this clear for the dimension and the signatures; for the determinant,
it follows from Proposition \ref{det}]. By Proposition \ref {Hasse-Witt},
we have $$ {\rm cor}_{F^v/K_v} (\lambda,d) = w_2(q') + w_2(q_{n,1})$$ in ${\rm Br}_2(K_v)$ for
all $v \in V_K$. On the other hand, since $\lambda
\in \mathcal {LR}^v$ for all $v \in V_K$, by Proposition \ref{atilda infinite places} we have ${\rm cor}_{F^v/K_v} (\lambda,d) = w_2(q) + w_2(q_{n,1})$
in ${\rm Br}_2(K_v)$ for
all $v \in V_K$; therefore 
$w_2(q') = w_2(q)$ in ${\rm Br}_2(K)$. This implies that $q' \simeq q$, hence condition
(G 1) holds. 

%\bigskip
%{\bf The obstruction group}
%and recall that $n_i = {\rm dim}_{E_i}(M_i)$. 

\medskip

Recall from
\S \ref {rational obstruction} the definition of the group $\sha_E = \sha_{E,\sigma}$; recall that
for all $v \in V_K$ and $i \in I$, the set $S_i^v$ consists of the places $w$ of $F_i$ such that $E_i \otimes_{F_i} (F_i)_w$
is a field, and that if $i \in I$, we denote by  $V_i$  the set of places $v \in V_K$ such that $S_i^v \not = \varnothing$.
The group $ \sha_E$
is the one constructed in \S \ref{obstruction}
 using the data $I$  and $V_{i,j} = V_i \cap V_j$; see \S \ref{rational obstruction} for details.
 
  \medskip
  Note that $\sha_E$ does not depend on the module $M$, only on the algebra with involution $E$.

\begin{theo}\label{iso field obstruction} {\rm (i)} Assume that 
$\sha_E= 0$, and let $q$ be a quadratic form such that {\rm (L 1)} holds. Then {\rm (G 1)} holds as well.

\medskip
{\rm (ii)} If $\sha_E \not = 0$, there exists a quadratic form satisfying  {\rm (L 1)} but not {\rm (G1)}.

\end{theo}

The proof of this theorem will be given later, after the construction of the obstruction homomorphism.

\medskip
If the obstruction group $\sha_E$  is not trivial, then the validity of the Hasse principle also depends
on the choice of the quadratic form.

\medskip

{\bf Local data}

\medskip
Let $q$ be a non-degenerate quadratic form over $K$, and let us assume that condition (L 1) holds. 
%for all $v \in V$, the quadratic 
%form $q$ has an isometry with characteristic polynomial $f$ and minimal polynomial $g$
%over $K_v$.  
Recall from Section \ref{quadr}  that a local solution gives rise to sets $\mathcal {LR}^v$ and $\mathcal {CR}^v$
for all $v \in V_K$.

\medskip

Let $\mathcal {CR}$ be the set of the elements $(a^v)$, $a^v \in \mathcal {CR}^v$, such that
$a^v = 0$ for almost all $v \in V_K$. Let us show that this set is not empty.

\begin{prop}\label{almost0} Assume that condition {\rm (L 1)} holds.
Then the set $\mathcal {CR}$ is not empty.
  
\end{prop}

\noindent
{\bf Proof.} Let $S$ be the subset of $V_K$ consisting of the places $v \in V_K$ such that $A^v \not = 0$ and the
infinite places; this is a finite set.

\medskip 
Let $v \in V_K$ be such that $v \not \in S$. Let $\lambda^v \in \mathcal {LR}^v$, and assume that there exists $i \in I$ such that $a(\lambda^v)(i) \not = 0$.
 Recall that
$A^v =  \underset{r \in I} \sum a(\lambda^v)(r)$. Since $v \not \in S$, by hypothesis $A^v = 0$; therefore 
there exists
$j \in I$, $j \not = i$, with $a(\lambda^v)(j) \not = 0$. 
Note that since $a(\lambda^v)(i) \not = 0$ and $a(\lambda^v)(j) \not = 0$,
we have $S_i^v \not = \varnothing$ and $S_j^v \not = \varnothing$, hence $v \in V_i \cap V_j$.

\medskip
 Note that $v$ is a finite place, since $v \not \in S$; hence by Proposition \ref{utile},
the map $(i,j)a(\lambda^v)$ belongs to ${CR}^v$. Moreover, this map vanishes at $i$ and $j$. 
If $(i,j)a(\lambda^v) = 0$, we stop; otherwise
we continue, and after a finite number of steps we obtain the zero element of $C(I)$. 
Since this holds for all $v \in V_K$ such that $v \not \in S$, the
proposition is proved.

\bigskip
{\bf The homomorphism and the Hasse principle}

\medskip
We apply the results of \S \ref{obstruction} with the sets $\mathcal C^v = \mathcal {CR}^v$, and the equivalence relation
$\sim_v$ will be the equivalence relation $ \equiv_v$, defined in Section \ref{equiv}. 

\medskip
For all $v \in V_K$, set  $\tilde A^v = {\rm inv}_v(w_2(q) + w_2(q_{n,1}))$. Let $\iota : {\bf Z}/2{\bf Z} \to {\bf Q}/{\bf Z}$ be the canonical injection, and
let $A^v \in {\bf Z}/2{\bf Z}$ be such that $\iota (A^v) = \tilde A^v$. 
Note that $q$ and $q_{n,1}$ are quadratic forms over $K$,  hence $\tilde A^v = 0$ for almost all $v \in V_K$, and
 $\underset {v \in V_K} \sum \tilde A^v = 0$; the same properties hold for $A^v$, therefore condition (i) of \S \ref{obstruction} holds. 

\medskip
 Let $V'$ be the set of finite places, and $V''$ be the set of infinite places of $K$. The sets  $\mathcal {CR}^v$
satisfy conditions (ii)-(iv) of Section \ref{obstruction} : indeed, condition (ii) follows from Proposition \ref{A infinite}, condition (iii) from 
Corollary \ref{utile}, 
and condition (iv) from Corollary \ref{rat utile}.

\medskip
Let $(a(\lambda^v)) \in \mathcal {CR}$. As in \S \ref{obstruction}, we define a homomorphism $\rho  : \sha_E \to {\bf Z}/2{\bf Z}$ as
follows. For all $c \in \sha_E$, set
$$\rho(c) = \underset {v \in V_K} \sum \  \underset{i \in I} \sum  \ c(i) a(\lambda^v)(i).$$
By Proposition \ref{independent}, the homomorphism $\rho$ is independent of the choice of $(a(\lambda^v)) \in \mathcal {CR}$.

 \begin{theo}\label {rat} Let $q$ be a quadratic form, and  assume that condition {\rm (L 1)} holds; then 
 condition {\rm (G 1)} holds if and only if $\rho = 0$.
 
 \end{theo}

 \noindent
 {\bf Proof.} If condition (G 1)  holds, then by Proposition \ref{cond}
 there exists $\lambda = (\lambda_i) \in T(E,\sigma)$ such
that $\lambda \in \mathcal {LR}^v$ for all $v \in V_K$; we have
 $\underset {v \in V_K} \sum a(\lambda)(i) = 0$ for all $i \in I$, hence $\rho = 0$.

 \medskip
  Let us prove the converse. Since $\rho = 0$, 
 Corollary \ref{ihomo} implies that there exists $b = (b^v) \in \mathcal {CR}$ such that 
  $\underset {v \in V_K} \sum b^v(i) = 0$ for all $i \in I$. By definition, there exists $(\lambda^v) \in \mathcal {LR}$ such
  that $b^v = a(\lambda^v)$ for all $v \in V_K$. 
 Recall that $a(\lambda^v)(i) = \underset {w \in S_v} \sum \theta_i^w(\lambda_i^w)$;
 therefore for all $i \in I$, we have  $\underset {w \in V_F} \sum \theta_i^w(\lambda_i^w) = 0$.
 By Theorem \ref{reciprocity} this implies that for all $i \in I$, there exists $\lambda_i \in T(E_i,\sigma_i)$ mapping to
 $\lambda_i^w \in T(E_i^w,\sigma_i)$ for all $w \in V_F$. In particular, we have $(\lambda_i,d_i) =
 (\lambda_i^v,d_i)$ for all $i \in I$ and $v \in V_K$.  Set
 $\lambda = (\lambda_i)$; since $\lambda^v \in \mathcal {LR}^v$ for all $v \in V_K$, Proposition \ref {cond} implies that (G~1) holds. 
 This concludes
 the proof of the theorem.

\medskip
We say that a quadratic form $q$ has {\it maximal signature} (with respect to $M$) if for all real places $v \in V_K$, the signature
of $q$  at $v$ is equal to the signature of $q_{n,1}$ or of $- q_{n,1}$  at $v$. 
For the proof of Theorem  \ref{iso field obstruction}, we need the following proposition. 

 \begin{prop}\label{lemma} For all $i,j \in I$, $i \not = j$, there exists a  quadratic form $q$ over $K$  having maximal
 signature satisfying {\rm (L 1)},  and, for
all $v \in V_K$,
a corresponding local data $\lambda^v \in \mathcal {LR}^v$, such that
for for some $v_1, v_2 \in V_K$, $v_1 \not = v_2$. we have

$$a(\lambda^{v_1})(i) = a(\lambda^{v_2})(j) = 1,$$ and
$$a(\lambda^v)(k) = 0 \ {\rm if} \ (v,k) \not = (v_1,i), (v_2,j).$$

\end{prop}

\begin{lemma}\label{realization} Let $v \in V_K$ be a finite place, and let $Q$ be a quadratic form over $K_v$ such
that ${\rm dim}(Q) = {\rm dim}(q)$ and ${\rm det}(Q) = {\rm N}_{F/K}(-d)$; then $Q$ is compatible with $M \otimes_KK_v$.

\end{lemma}

\noindent {\bf Proof.} Note that $q_{n,1}$ is compatible with $M$ by construction, and that its dimension and determinant coincide
with those of $Q$; if $w_2(Q) = w_2(q_{n,1})$ in ${\rm Br}_2(K_v)$, then $Q \simeq q_{n,1}$ over $K_v$, hence we are done.
Suppose that $w_2(Q) \not = w_2(q_{n,1})$, and let $\lambda$ be a non-trivial element of $T(E_v,\sigma)$. By Proposition
\ref{general Hasse-Witt}, we have $w_2(q_{n,\lambda}) \not = w_2(q_{n,1})$, and therefore $w_2(q_{n,\lambda}) = w_2(Q)$ in ${\rm Br}_2(K_v)$. 
The forms $q_{n,\lambda}$ and $Q$ have the same dimension, determinant and Hasse-Witt invariant, hence they are isomorphic over $K_v$.
Since $q_{n,\lambda}$ is compatible with $M \otimes_K K_v$, so is $Q$.

\medskip
\noindent
{\bf Proof of Proposition \ref{lemma}.} 
For all $k \in I$, let $S(k)$ be the set of places $v$ of $V_K$ such that $w_2(q^v_{n_k,1}) \not = 0$, and let $S$ be
the set of places $v \in V_K$ such that $v \in S(k)$ for some $k \in I$, or that the quaternion algebra $(d_s,d_r)$ is not split 
at $v$ for some $r,s \in I$. 
Let $v_1, v_2 \in V_K$ be two distinct finite places that are not in the finite set $S$. 

\medskip
For all $k \in I$ and $v \in V_K$, let $q^v_k$ be a quadratic form over $K_v$ such that ${\rm dim}(q^v_k) = {\rm dim}(q_{n,1})$, that
${\rm det}(q_k^v) = {\rm N}_{F_k/K}(-d_k)$, that if $v$ is a real place, then the signature of $q_k^v$ is equal to the signature of $q_{n,1}$ at $v$, and
that the Hasse-Witt invariants of $q_k^v$ are as follows 

\smallskip
$\bullet$ ${\rm inv}(w_2(q_i^{v_1})) = {\rm inv}(w_2(q_j^{v_2})) = {1\over 2}$;

\smallskip
$\bullet$ If $ (v,k) \not = (v_1,i), (v_2,j)$, then $w_2(q_k^v) = w_2(q_{n_k,1})$ in ${\rm Br}_2(K_v)$. 

\medskip
For all $v \in V_K$, set  $Q^v = \underset{k \in I} \oplus q_k^v$; the quadratic form $Q^v$ has determinant ${\rm N}_{F/K}(-d)$, hence by
Lemma \ref{realization} it is compatible with $M \otimes_K K_v$ if $v$ is a finite place. If $v$ is an infinite place, then $Q^v$ is isomorphic
to $q_{n,1}$ over $K_v$ by construction, hence it is compatible with $M \otimes_K K_v$.

\medskip
We claim that the number of $v \in V_K$ such that $w_2(Q^v) \not = 0$ is even. Indeed, for all $v \in V_K$ and all $k \in I$, we have

$$w_2(Q^v) = \underset {k \in I} \sum w_2(q_k^v) + \underset {r < s} \sum (d_r,d_s),$$ and
 
$$w_2(q_{n,1}) = \underset {k \in I} \sum w_2(q_{n_k,1}) + \underset {r < s} \sum (d_r,d_s).$$

If $v \not = v_1, v_2$, this implies that $w_2(Q^v) = w_2(q_{n,1})$, and
since 
$q_{n,1}$ is a global form, the number of $v \in V_K$ such that $w_2(q_{n,1}) \not = 0$ in ${\rm Br}(K_v)$ is even. We have $w_2(q^{v_1}) \not = w_2(q_{n,1}^{v_1})$, and $w_2(q^{v_2}) \not = w_2(q_{n,1}^{v_2})$; hence the number of places $v$ such  $w_2(Q^v) \not = 0$ is even. 

\medskip
Let  $q$  be a quadratic form over $K$ such that $q^v \simeq Q^v$ over $K_v$ for all $v \in V_K$; the existence of this form follows from
\cite{Sch}, Theorem 6.6.10. The form $q$ has maximal signature by construction.

\medskip For all $v \in V_K$, let $\lambda^v \in \mathcal {LR}^v$ be the local data corresponding to the quadratic form $Q^v$. 
We claim that $a(\lambda^v) \in \mathcal {CR}^v$ satisfies the required conditions. Indeed, recall that for all $v \in V_K$ and all $k \in I$, we have
$$w_2(q_k^v) = w_2(q^v_{n_k,1}) + {\rm cor}_{F^v/K_v}(\lambda_k^v,d_k),$$ 
cf. Proposition \ref{atilda infinite places}.
Since $v_1$ and $v_2$ are not in $S$, we have $w_2(q_{n_i,1}^{v_1}) = w_2(q_{n_j,1}^{v_2}) = 0$, hence
$${\rm inv}({\rm cor}_{F^v/K_v}(\lambda^{v_1},d_i)) = {\rm inv} ( {\rm cor}_{F^v/K_v}(\lambda_j^{v_2},d_j) )= {1 \over 2},$$
in other words, we have $\tilde a(\lambda^{v_1})(i) = \tilde a(\lambda^{v_2})(j) = 1.$
By Lemma \ref{tilde},
this implies that $$a(\lambda^{v_1})(i) = a(\lambda^{v_2})(j) = 1.$$ 
If $ (v,k) \not = (v_1,i), (v_2,j)$, then $w_2(q_k^v) = w_2(q^v_1(k))$, hence the same argument shows that $\tilde a(\lambda^{v})(k) = 0$,
and therefore 
by Lemma \ref{tilde}, 
we have $a(\lambda^v)(k) = 0$. This completes the proof of the Proposition.

\medskip
\noindent
{\bf Proof of Theorem \ref{iso field obstruction}}. If $\sha_E = 0$, then by Theorem \ref{rat} the Hasse principle holds for
any quadratic form $q$.

\medskip
To prove the converse, assume that $\sha_E \not = 0$; we claim that there exists a quadratic form $q$ satisfying (L 1) but not (G 1). Let
$c \in \sha_E$ be a non-trivial element, and let $i, j \in I$ be such that $c(i) \not = c(j)$. With $q$  and $a \in \mathcal {CR}$
as in  Proposition \ref{lemma}, we have

$$\rho(c) = \underset {v \in V_K} \sum \  \underset{i \in I} \sum  \ c(i) a(\lambda^v)(i) = c(i) + c(j),$$
hence $\rho(c) \not = 0$;  by Theorem \ref{rat} condition (G 1) does not hold.

\begin{example}\label{isometry 3} Let $A = K[\Gamma]$ as in example \ref{isometry 2}, and assume that all the simple $\sigma$-stable
factors in $\mathcal A_M$ are of type (1). In other words, we have $\mathcal A_M = (E_i)_{i \in I}$
%\underset {i \in I} \prod E_i$ 
with $E_i = K[X]/(f_i)$, where $f_i \in K[X]$
are monic, irreducible,  symmetric polynomials of even degree, and $M = \underset{i \in I} \oplus M_i$, with $M_i  = [K[X]/(f_i)]^{n_i}$
for some integers $n_i \geqslant 1$. Let
$q$ be a quadratic form over $K$; then $q$ is compatible with $M$ if and only if $q$ has an isometry with minimal polynomial
$g = \underset {i \in I} \prod  f_i$, and characteristic polynomial $f = \underset {i \in I} \prod f_i^{n_i}$; hence theorem \ref{rat} gives a
necessary and sufficient condition for the Hasse principle to hold for the existence of an isometry with minimal polynomial $g$ and
characteristic polynomial $f$. 

\end{example}

\begin{example}\label{cyclotomic 1} With the notation of example \ref{isometry 3}, assume that $K = {\bf Q}$, and that the polynomials $f_i$
are cyclotomic polynomials for all $i \in I$. Then $E_i = {\bf Q}/(f_i)$ is a cyclotomic field for all $i \in I$, hence a CM field; by example
\ref{CM} this implies that $\sha_E = 0$. By Theorem \ref {iso field obstruction}, the Hasse principle holds for the existence of an isometry with minimal polynomial $g$ and
characteristic polynomial $f$; if $q$ is a quadratic form having an isometry with minimal polynomial $g$ and characteristic polynomial
$f$ locally everywhere, then such an isometry exists over $\bf Q$ as well.

\end{example}

\section{Independent extensions}\label{independent extensions}

Recall  that two finite extensions
$K_1$ and $K_2$ of $K$ are  {\it independent over $K$} if the tensor product $K_1 \otimes_K K_2$ is a field. In this section, we show
that the local-global principle of \S \ref{rational} always holds if the extensions $E_i/K$ are pairwise
independent over $K$. 

\medskip We keep the notation of \S \ref{rational}. Recall from Example \ref{sha indep} the equivalence relation $\approx$ on $I$ generated by the elementary equivalence

\medskip
\centerline {$i \approx_e j \iff$ $E_i$ and $E_j$ are independent over $K$,}

\medskip
\noindent
and that $\sha_{\rm indep}(E) = \sha_{\approx}(I)$  is the associated obstruction group. 

\begin{theo}\label{sha indep 0} 
Assume that $\sha_{\rm indep}(E) = 0$ 
and that  condition {\rm (L 1)} holds. 
Then condition {\rm (G 1)} holds
as well.

\end{theo}

\noindent
{\bf Proof.} By Corollary \ref {0 and 0} the hypothesis implies that $\sha_E =0$; therefore Theorem \ref {iso field obstruction} (i) gives the required result. 

\begin{coro} Assume that there exists $i \in I$ such that for all $j \in I$ with $j \not =i$ the extensions $E_i$ and $E_j$
are independent over $K$. If condition {\rm (L 1)} is satisfied, then
condition {\rm (G 1)} also holds.

\end{coro}

\noindent
{\bf Proof.} This follows immediately from Theorem \ref{sha indep 0}, since the hypothesis implies that $\sha_{\rm indep}(E) = 0$.

\medskip The following corollary is an immediate consequence :

\begin{coro}\label {pairwise independent} Assume that the extensions $E_i/K$ are pairwise
independent over $K$, and that condition {\rm (L 1)} holds. Then
condition {\rm (G 1)} holds as well.

\end{coro}

\begin{example}\label{independent polynomials} Let $f = \underset {i \in I} \prod f_i^{n_i}$ and $g = \underset {i \in I} \prod  f_i$ be as in Example \ref{isometry 3}, with $E_i = K[X]/(f_i)$. Assume that $\sha_{\rm indep}(E) = 0$, and let
$q$ be a quadratic form over $K$. Then the Hasse principle holds;
%for the existence of an isometry with minimal polynomial $g$ and
%characteristic polynomial $f$
if $q$ has an isometry with minimal polynomial $g$ and characteristic polynomial
$f$ locally everywhere, then such an isometry exists over $K$. In particular, this is the case if the extensions
$E_i/K$ are pairwise  independent over $K$.

\end{example}

\section {Local conditions} 

We keep the notation of \S \ref{rational}. Recall that $(V,q)$ is a quadratic form, that
$M$ is an $A$-module satisfying the hypotheses
of \S \ref{quadr}, and $\mathcal A_M = (E_i)_{i \in I}$.  Recall that $E_i/F_i$ is a quadratic extension 
for all $i \in I$, and $d_i \in F_i^{\times}$ is such that $E_i = F_i(\sqrt{d_i})$.

\medskip The aim of this section is to give a necessary and sufficient condition
for $M \otimes_K K_v$ and $(V,q) \otimes_K K_v$ to be compatible for all $v \in V_K$, in other words, for
condition (L 1) to hold. This complements the Hasse principle results of \S \ref{rational} and will also be used
in \S \ref{odd degree}. Since in \S \ref{odd degree} we need this result for several fields, we use the notation
with a suffix $K$ : 

\medskip
${\rm (L \ 1)}_K$  {\it  \  \ For all $v \in V_K$, the quadratic form $(V,q)\otimes_KK_v$ is compatible with the module $M \otimes_KK_v$.}

\medskip The validity of ${\rm (L \ 1)}_K$  will be detected by three conditions, called the determinant condition, the hyperbolicity condition and the signature condition,

\medskip
{\bf Determinant condition}

\medskip
${\rm (det)}_K$  {\it \hskip 2cm   ${\rm det}(q) = \underset{i \in I} \prod {\rm N_{E_i/K}}(-d_i)^{n_i}$ \ \  in $K^{\times}/K^{\times 2}$.}

\begin{remark} \label{determinant} Note that by Proposition \ref{det}, condition ${\rm (det)}_K$ is necessary for
$(V,q)$ and $M$ to be compatible for any field $K$. 

\end{remark}

\medskip
{\bf Hyperbolicity condition} 

\medskip
${\rm (hyp)}_K$  {\it If $v \in V_K$ is such that $E^v = F^v \times F^v$, then $(V,q) \otimes_K K_v$ is hyperbolic.}

\bigskip
{\bf Signature condition}

\medskip If $v \in V_K$ is a real place, we denote by $(r_v,s_v)$ the signature of $(V,q) \otimes_K K_v$, by $(\rho_i)_v$
the number of real places of $F_i$ above $v$ that extend to complex places of $E_i$, and set 
$\sigma_v = {\rm dim}(M) - 2 \underset{i \in I} \sum n_i (\rho_i)_v$. 

\begin{prop}\label{signature condition} Let $v \in V_K$ be a real place. Then $(V,q) \otimes_K K_v$ is
compatible with $M \otimes_K K_v$ if and only if $r_v \geqslant \sigma_v$, $s_v \geqslant \sigma_v$, and
$r_v \equiv s_v \equiv \sigma_v \ {\rm (mod \ 2)}$.

\end{prop}

\noindent
{\bf Proof.} This is straightforward, using for instance the arguments of the proof of \cite {B4}, Proposition 8.1. 

\medskip The signature condition is as follows :

\medskip
${\rm (sign)}_K$  {\it If $v \in V_K$ is a real place, we have 

\medskip \centerline {$r_v \geqslant \sigma_v$, $s_v \geqslant \sigma_v$, and
$r_v \equiv s_v \equiv \sigma_v \ {\rm (mod \ 2)}$.}}

\begin{prop}\label{local cond} The following are equivalent

\smallskip {\rm (a)} Condition ${\rm (L \ 1)}_K$  holds.

\smallskip {\rm (b)} Conditions ${\rm (det)}_K$, ${\rm (hyp)}_K$ and ${\rm (sign)}_K$ hold.

\end{prop}

\noindent {\bf Proof.} Let us show that (a) $\implies$  (b). Condition ${\rm (L \ 1)}_K$ implies ${\rm (det)}_K$
by Proposition \ref{det}, ${\rm (hyp)}_K$ by Proposition \ref{decomposition}, and ${\rm (sign)}_K$ by
Proposition \ref{signature condition}. Conversely, let us show that (b) $\implies$ (a). Let $v \in V_K$. If
$v$ is a real place, then Proposition \ref{signature condition} shows that 
$(V,q) \otimes_K K_v$ is
compatible with $M \otimes_K K_v$. Assume now that $v$ is a finite place If $E^v = F^v \times F^v$, then
by ${\rm (hyp)}_K$ the form $(V,q) \otimes_K K_v$ is hyperbolic, and by Proposition \ref{decomposition} it
is compatible with $M \otimes_K K_v$. Suppose that $E^v \not = F^v \times F^v$; then 
$T(E^v,\sigma)$ is non-trivial. For all $\lambda \in T(E^v,\sigma)$,
we have ${\rm det}(q_{n,\lambda}) = {\rm det}(q)$ by Proposition \ref{det} and condition ${\rm (det)}_K$,
and by Proposition \ref{general Hasse-Witt}  we can choose $\lambda \in T(E^v,\sigma)$ so that $w_w(q_{n,\lambda} ) = w_2(q \otimes_K K_v)$
in ${\rm Br}(K_v)$. Therefore $(V,q) \otimes_K K_v$ is
compatible with $M \otimes_K K_v$, and hence ${\rm (L \ 1)}_K$  holds.

\begin{example}\label{local polynomials} Assume that $M$ and $E$ are as in Example \ref{isometry 3}; recall
that $E_i = K[X]/(f_i)$  and $M_i = [K[X]/(f_i)]^{n_i}$ with $f_i \in K[X]$ irreducible, symmetric polynomials of even
degree. Set $f = \underset {i \in I} \prod f_i^{n_i}$. In this case,
the local conditions can be reformulated as follows :

\medskip
${\rm (det)}_K$  {\it \hskip 2cm \ \ ${\rm det}(q) = f(1) f(-1)$  \ \  in $K^{\times}/K^{\times 2}$.}

\bigskip
We say that a polynomial is {\it hyperbolic} if it is a product of irreducible polynomials of type (2) (see Example
\ref{isometry 2}). 

\medskip
${\rm (hyp)}_K$  {\it If $v \in V_K$ is such that  $f \in K_v[X]$ is hyperbolic, then $(V,q) \otimes_K K_v$ is hyperbolic.}

\bigskip
If $v \in V_K$ is a real place, we denote by $m_v(f)$ the number of roots of $f \in K_v[X]$ with $|z|_v > 1$ (counted
with multiplicity). 

\medskip
${\rm (sign)}_K$  {\it If $v \in V_K$ is a real place, we have 

\medskip \centerline {$r_v \geqslant m_v(f)$, $s_v \geqslant m_v(f)$, and
$r_v \equiv s_v \equiv m_v(f) \ {\rm (mod \ 2)}$.}}

\bigskip Set
$g = \underset {i \in I} \prod  f_i$. We recover a result of \cite {B4}, Theorem 12.1 : the quadratic form $(V,q) \otimes_K K_v$ has an isometry with characteristic polynomial $f$
and minimal polynomial $g$ for all $v \in V_K$  if and only if the above three conditions hold.

\end{example}

\section{Odd degree descent}\label{odd degree}

We keep the notation of the previous sections. The aim of this section is to prove an ``odd degree descent" result :

\begin{theo}\label{odd} If $K'$ is a finite extension of odd degree of $K$ such that
$(V,q) \otimes_K K'$ compatible with $M \otimes_K K'$,
then
$(V,q)$ is compatible with $M$.

\end{theo}

We start with a lemma

\begin{lemma}\label{L 1 K'}  Let $K'/K$ be a finite extension of odd degree. Then $${\rm (L \ 1)}_{K'} \implies {\rm (L \ 1)}_K.$$

\end{lemma}

\noindent
{\bf Proof.} It is clear that ${\rm (det)}_{K'} \implies {\rm (det)}_K$, ${\rm (hyp)}_{K'}  \implies {\rm (hyp)}_K$ and ${\rm (sign)}_{K'} \implies {\rm (sign)}_K$. By Proposition \ref{local cond}, this implies that ${\rm (L \ 1)}_{K'} \implies {\rm (L \ 1)}_K.$

\medskip
\noindent
{\bf Proof of Theorem \ref {odd}.} %By Proposition \ref {odd sha}, we have and injective homomorphism $\sha_E \to \sha_{E \otimes_K K'}$. 
By hypothesis,
condition ${\rm (L \ 1)}_{K'}$ holds; let $$\rho_{K'} : \sha_{E \otimes_K K'} \to {\bf Z}/2{\bf Z}$$ be the associated homomorphism
of \S \ref{rational}. Recall that $\rho_{K'}$ is independent of the chosen local data. By 
Lemma \ref {L 1 K'}, condition ${\rm (L \ 1)}_K$  holds. Let

$$\rho_K(c) = \underset {v \in V_K} \sum \  \underset{i \in I} \sum  \ c(i) a(\lambda^v)(i)$$ be the associated homomorphism.
For all  $v\in V_K$, let us choose a place  $w$ of $K'$ over $v$.
Let us take the extension of $a^w$ of $a^v = a(\lambda^v)$ to $F^v\otimes_{K_v} K'_w$ to define $\rho_{K'}$.

\medskip
If $[K'_w:K_v]$ is odd, then $${\rm inv}_v{\rm cor}_{F_i^v/K_v} (a_i^v,d_i)={\rm inv}_w{\rm cor}_{F_i^v\otimes_{K_v}K'_w/K'_w} (a_i^w,d_i).$$

\medskip
If $[K'_w:K_v]$ is even, then $${\rm inv}_w{\rm cor}_{F_i^v\otimes_{K_v}K'_w/K'_w} (a_i^w,d_i)=0.$$

\medskip
Since $K'$ is an odd degree extension of $K$, the degree  $[K'_w:K_v]$ is odd for an odd number of places $w$ of $K'$ over $v$.
%that $[K'_w:K_v]$ is  odd.
Hence \[{\rm inv}_v{\rm cor}_{F_i^v/K_v} (a_i^v,d_i)=\underset{w|v}{\sum}{\rm inv}_w{\rm cor}_{F_i^v\otimes_{K_v}K'_w/K'_w} (a_i^w,d_i).\]

%\medskip
%and let us take the extension of $a(\lambda^v)$ to $K'$ to define $\rho_{K'}$. Note that for all $c \in \sha_E$, we have 

%$$\underset {v' \in V_{K'}} \sum \  \underset{i \in I} \sum  \ c(i) a(\lambda^{v'})(i) = \underset {v \in V_K} \sum \  \underset{i \in I} \sum  \ c(i) a(\lambda^v)(i).$$ By hypothesis, condition  ${\rm (G  \ 1)}_{K'}$ holds, hence by Theorem \ref{rat}  we have $\rho_{K'} = 0$; by the above argument,  this implies that we have $\rho_K = 0$ as well. Therefore
%by Theorem \ref{rat} the quadratic form $(V,q)$ is compatible with $M$.
Recall from \S \ref{rational obstruction} that  we have an injective homomorphism $$\pi' : \sha_E \to \sha_{E \otimes_K K'}$$
(see Proposition \ref {odd sha}). Let $c \in \sha_E $, and let $c' = \pi'(c)$; in other words, with the notation of \S 
\ref{rational obstruction}, we have $c'(i,j)=c(i)$.
Let $a^{v'}_{i,j}$ be the image of $a^v_i$  in $F'_{i,j}\otimes_{K'} K'_{v'}$; then ($a^{v'}_{i,j}$) is a local
data over $K'$.
Denote by $d_{i,j}$  the image of $d_i$ in $F'_{i,j}$.
Note that $F_i\otimes_K K'\otimes_{K'} K'_w\simeq F_i^v\otimes_{K_v}K'_w$.

\medskip By hypothesis, condition  ${\rm (G  \ 1)}_{K'}$ holds, hence by Theorem \ref{rat}  we have $\rho_{K'} = 0$. Therefore
 $\rho_{K'}(c') = 0$.
%$$\underset {v' \in V_{K'}} \sum \  \underset{(i,j) \in I'} \sum  \ c'(i,j) \ {\rm inv}_{v'} {\rm cor}_{F_{i,j}^{v'}/K'_{v'}} (a_{i,j}^{v'},d_{i,j})\\
%= 0. $$
Using the above observations, we see that $\rho_K(c) = \rho_{K'}(c')$. 
%\begin{align*}
%&\underset {v' \in V_{K'}} \sum \  \underset{(i,j) \in I'} \sum  \ c'(i,j) \ {\rm inv}_{v'} {\rm cor}_{F_{i,j}^{v'}/K'_{v'}} (a_{i,j}^{v'},d_{i,j})\\
% &=\underset {v' \in V_{K'}} \sum \  \underset{i \in I} \sum  \ c(i)\underset{j\in S(i)}{\sum} \ {\rm inv}_{v'} {\rm cor}_{F_{i,j}^{v'}/K'_{v'}} (a_{i,j}^{v'},d_{i,j})\\
 %&=\underset {v' \in V_{K'}} \sum \  \underset{i \in I} \sum  \ c(i) \ {\rm inv}_{v'} {\rm cor}_{F_{i}^{v}\otimes_{K_v} K'_{v'}/K'_{v'}} (a_{i}^{v'},d_{i})\\
 %&= \underset {v \in V_K} \sum \  \underset{i \in I} \sum \ c(i) \ {\rm inv}_v {\rm cor}_{F^v/K_v} (a_i^v,d_i).
%\end{align*}
This implies that $\rho_K = 0$; therefore
by Theorem \ref{rat} the quadratic form $(V,q)$ is compatible with $M$.

\section{Obstruction group - the integral case}\label{integral obstruction}

As in the previous sections, $K$ is a global field; let $O$ be a ring of integers of $K$ with respect to a finite, non-empty set $\Sigma$ of places
of $K$, containing the infinite places if $K$ is a number field. Let $V_{\Sigma}$ be the set of places of $K$ that are not in $\Sigma$. 
If $v \in V_{\Sigma}$, we denote by $O_v$ the ring of integers of $K_v$, and by $k_v$ its residue field. Let
$\Lambda$ be an $O$-algebra, and let $\sigma : \Lambda \to \Lambda$ be an $O$-linear involution;
%we use the notation of the previous sections with 
set $A = \Lambda_K = \Lambda \otimes_OK$. If $v \in V_{\Sigma}$, set  $\Lambda_{K_v} = \Lambda \otimes _O K_v$, and $\Lambda_{k_v} = \Lambda \otimes_Ok_v$. 

\medskip Let $(M_i)_{i \in I}$ be a finite set of $A$-modules,  let $M = \underset{i \in I}\oplus M_i$, and 
%Let $M$ be an $A$-module as in \S \ref{quadr}, and 
let ${\Lambda}_M$ be the image of $\Lambda$ in ${\rm End}(M)$. We assume that the kernel of the homomorphism $\Lambda \to \Lambda_M$ is
stable by the involution $\sigma$, and we also denote by $\sigma : \Lambda_M \to \Lambda_M$ the induced involution.
%; in particular, $\Lambda_M$ is a subring
%of $A_M$. Assume that $\Lambda_M$ is stable by the involution $\sigma$. 
%Recall that $M = \underset{i \in I}\oplus M_i$ with $M_i \simeq E_i^{n_i}$ for all $i \in I$. 

\medskip
The aim of this section is to define a group $\sha_{\Lambda_M,(M_i)_{i \in I}}$  that will be useful in \S \ref{isometrieslattices}. 
%This group maps surjectively to the
%group $\sha_E$, defined in \S \ref{rational obstruction}; 
In general, this group depends on $(M_i)_{i \in I}$ as well as $\Lambda_M$. However, in our main case of interest,
namely when $\Lambda = O[\Gamma]$, it only depends on $\Lambda_M$.

\medskip
%Recall that $E = \underset{i \in I}\prod E_i$; 
The group $\sha_{\Lambda_M,(M_i)_{i \in I}}$ is defined using the general framework of \S \ref {obstruction} using the set $I$
and for all $i,j \in I$, the subsets $V_{i,j}$ of $V_{\Sigma}$ defined as follows.

\begin{notation}
For all $i,j \in I$, we denote by  $V_{i,j}$ the set of places of $V_{\Sigma}$ such that there exists a self-dual
 $\Lambda_{k_v}$-module appearing in the reductions mod $\pi_v$ of both $M_i^v$ and $M_j^v$.
 \end{notation}

 \begin{example}\label{isometry 4} 
 Assume that $\Lambda = O[\Gamma]$,  hence  $A = \Lambda \otimes_O K = K[\Gamma]$.
We keep the notation of Example \ref{isometry 3}, and assume that  $f_i \in O[X]$ for all $i \in I$.
%set $g = \underset {i \in I} \prod  f_i$, and $f = \underset {i \in I} \prod f^{n_i}$.
%If $v \in V$ is a finite place, we denote by $O_v$ the ring of integers of $K_v$, and by $k_v$ the residue field. 
The homomorphism $O_v \to k_v$ induces $p_v : O_v[X] \to k_v[X]$. 
For all $i,j \in I$, the set $V_{i,j}$ defined above is the set of  places
$v \in V_{\Sigma}$ such that
the polynomials $p_v(f_i)$ and $p_v(f_j)$ have a common irreducible and symmetric factor; this follows from Example \ref{isometry 0}. 
%The sets $V_{i,j}$ are finite, and $V_{i,j} \subset V_i \cap V_j$.

 \end{example}

 \medskip
Assume now that the $A$-module $M$ and its decomposition $M = \underset{i \in I}\oplus M_i$ are as in \S \ref {quadr}, with $A_M = \underset{i \in I} \prod E_i$ and 
$M_i \simeq E_i^{n_i}$ for all $i \in I$; 
%The decomposition $M = \underset{i \in I}\oplus M_i$ is then given by the set $\mathcal A_M = (E_i)_{i \in I}$,
set $\sha_{\Lambda_M,(M_i)_{i \in I}} = \sha_{\Lambda_M,M}$.

\medskip
Recall that $S_i$ is the set of places $w$ of $F_i$ that are inert or ramified in $E_i$.

 \begin{notation} If $w_i \in S_i$, we denote by $O^{w_i}$ the ring of integers of $E_i^{w_i}$, and by $\kappa_{w_i}$ its 
 residue field. 
 \end{notation}
 
 Recall that if $N$ is a self-dual $\Lambda_{K_v}$-module, we denote by $\sigma_N : \kappa(N) \to \kappa(N)$ the induced
 involution of $\kappa(N) = (\Lambda_{k_v})_N$. 
 
 \begin{prop}\label{kappa} Let $v \in V_{\Sigma}$, and let $i,j \in I$. Assume that $v \in V_{i,j}$. Then there exists a simple, self-dual 
 $\Lambda_{k_v}$-module
 $N$
 %a $k_v$-linear involution $\sigma_N : \kappa(N) \to \kappa(N)$, 
 and places $w_i \in S_i^v$, 
 $w_j \in S_j^v$ such that 
 
 \medskip
 {\rm (i)} The $\Lambda_{k_v}$-modules $\kappa_{w_i}$, $\kappa_{w_j}$ and $N$ are isomorphic.

 \medskip
 {\rm (ii)} The fields with
 involution $(\kappa_{w_i},\sigma_i)$
 and $(\kappa_{w_j},\sigma_j)$ are isomorphic to $(\kappa(N),\sigma_N)$.
 
 \end{prop}

\noindent 
{\bf Proof.}
(i)  follows from Proposition \ref{reduction}. To prove (ii), note that $({\Lambda_{k_v}})_{\kappa_{w_i}} \simeq \kappa_{w_i}$ and
$({\Lambda_{k_v}})_{\kappa_{w_j}} \simeq \kappa_{w_j}$. Since all the involutions are induced by $\sigma$, this implies (ii).

 \begin{coro}\label{rat int} For all $i,j \in I$, we have $V_{i,j} \subset V_i \cap V_j$. 
 
 \end{coro}
 
 \noindent
 {\bf Proof.} This is an immediate consequence of Proposition \ref{kappa}, recalling that $v \in V_i$ if and only if $S_i^v \not = \varnothing$.
 %$w \in S_i$ above $v$. 
 
 \begin{coro}\label{subgroup} $\sha_E$  is a subgroup of $\sha_{\Lambda_M,M}$.
 
 \end{coro}
 
 \noindent
 {\bf Proof.} Let us denote by $\sim$ the equivalence relation on $I$ generated by $$i \sim j \iff V_{i,j} \not  = \varnothing,$$ and by $\approx$
 the equivalence relation generated by $$i \approx j \iff V_{i} \cap V_j  \not  = \varnothing.$$ By Corollary \ref{rat int}, we have 
 $V_{i,j} \subset V_i \cap V_j$ for all $i, j \in I$, hence $i \sim j \implies i \approx j$. Since $\sha_{\Lambda_M,M}$ is defined by the
 equivalence relation $\sim$ and $\sha_E$ by $\approx$, this implies that $\sha_E$  is a subgroup of $\sha_{\Lambda_M,M}$.

\section{Local data and residue maps}\label{localdata}

We keep the notation of \S \ref{integral obstruction}; in particular, $M = \underset{i \in I}\oplus M_i$ is a module with $M_i \simeq E_i^{n_i}$ for all $i \in I$. 

\medskip Let $v \in V_{\Sigma}$. If $w_i \in S_i^v$ and $\lambda_i^{w_i} \in T(E_i^{w_i},\sigma_i)$, we obtain a bounded $E_i^{w_i}$-quadratic
form $((E_i^{w_i})^{n_i},q_{n_i,\lambda_i^{w_i}})$ (see \S \ref{residuemaps}). Set $M_i^{w_i} = M_i \otimes_{E_i}E_i^{w_i}$. Choosing an isomorphism $M_i \to E_i^{n_i}$, we obtain a bounded
$A$-form on $M_i^{w_i}$, denoted by $(M_i^{w_i},q_{\lambda_i^{w_i}})$.

\medskip
Recall from \S \ref{res} that we have a homomorphism $$\partial_v : W_{\Lambda_{K_v}}^b(K_v) \to W_{\Lambda_{k_v}}(k_v).$$ 

Let $\lambda^v = (\lambda_i^{w_i})$. Set
%If $w_i \in S_i^v$ for some $i \in I$, we have a map
%$$t_{n_i} : T(E_i^{w_i}) \to W_{\Lambda_{k_v}}(k_v)$$
%sending $\lambda_i^{w_i} \in T(E_i^{w_i})$ to the class of the bounded form $((E_i^{w_i})^{n_i},q_{n_i,\lambda_i^{w_i}})$ (see \S \ref{residuemaps}). Set

$$\partial_v(\lambda_i^{w_i}) = \partial_v[(M_i^{w_i},q_{\lambda_i^{w_i}})]$$
%\partial_v[(E_i^{w_i})^{n_i},q_{n_i,\lambda_i^{w_i}}],$$

$$\partial_v(\lambda_i^v) = \underset{w_i \in S_i^v} \oplus \ \partial_v (\lambda_i^{w_i}), \ {\rm and}$$ 

 $$\partial_v(\lambda^v) =  \underset {i \in I} \oplus \ \partial_v(\lambda_i^v).$$

%\begin{notation}  If $w_i \in S_i$, we
%denote by $O_i^{w_i}$ the ring of integers of $E_i^{w_i}$ and by $\kappa_i^{w_i}$ its residue field. 
%For all $i,j \in I$, let $V_{i,j}$ be the set of $v \in V_K$ such that there exist places $w_i \in S_i^v$ and $w_j \in S_j^v$
%such that the fields $\kappa_i^{w_i}$ and $\kappa_j^{w_j}$ are isomorphic. 
%\end{notation}

\begin{prop}\label{nullcobordant} Let $v \in V_{\Sigma}$, and let $\lambda^v, \mu^v \in 
\mathcal {LR}^v$ be such that $a(\mu^v) = (i,j)(a(\lambda^v))$ for some $i,j \in I$ with $v \in V_{i,j}$.
Then we have $\partial_v(\mu^v) = \partial_v(\lambda^v)$. 

\end{prop}

\medskip
\noindent
{\bf Proof.}
Since $v \in V_{i,j}$, by Proposition \ref{kappa} there exists a simple, self-dual 
 $\Lambda_{k_v}$-module
 $N$
 and places $w_i \in V_{F_i}$, $w_j \in V_{F_j}$ such that the $\Lambda_{k_v}$-modules $\kappa_{w_i}$, $\kappa_{w_j}$ and $N$ are isomorphic,
 and that
the fields with
 involution $(\kappa_{w_i},\sigma_i)$
 and $(\kappa_{w_j},\sigma_j)$ are isomorphic to $(\kappa(N),\sigma_N)$.

\medskip
Set 
$$Q = (\partial_v(\lambda_i^{w_i})) \oplus (\partial_v(\lambda_j^{w_j})) \oplus (- \partial_v(\mu_i^{w_i})) \oplus (- \partial_v(\mu_j^{w_j})).$$ 

\medskip
We claim that $Q = 0$ in $W_{\Lambda_{k_v}}(k_v)$; note that $Q$ belongs to the subgroup $W_{\Lambda_{k_v}}(k_v,N)$ of $W_{\Lambda_{k_v}}(k_v)$,
and that $W_{\Lambda_{k_v}}(k_v,N) \simeq W_{\kappa(N)}(k_v)$. 

%\medskip
%Set $E_{i,j} = \kappa_i^{w_i} \simeq \kappa_j^{w_j}$, and let $\sigma_{i,j} : E_{i,j} \to E_{i,j}$ be the $k_v$-linear involution induced by $\sigma$.

\medskip
Let us first assume that $w_i$ and $w_j$ are inert in $E_i^{w_i}$, respectively $E_j^{w_j}$. In
this case, the involution $\sigma_N$ is non-trivial. 
By Proposition \ref{residue} (a), we know that $\partial_v : T(E_i^{w_i},\sigma_i) \to W_{\kappa(N)}(k_v)$ and $\partial_v : T(E_j^{w_j},\sigma_j) \to  W_{\kappa(N)}(k_v)$ are bijective. Therefore
we have 
$\partial_v(\lambda_i^{w_i}) \not = \partial_v(\mu_i^{w_i})$
and
$\partial_v(\lambda_i^{w_i}) \not = \partial_v(\mu_j^{w_j})$. This implies that two of the four
elements  $\partial_v(\lambda_i^{w_i}), \partial_v(\lambda_j^{w_j}),  \partial_v(\mu_i^{w_i}), \partial_v(\mu_j^{w_j})$ are trivial, and two are non-trivial. Since $W_{\kappa(N)}(k_v)$ is
of order $2$, this implies that the class of $Q$ is trivial in $W_{\kappa(N)}(k_v)$, hence in $W_{\Lambda_{k_v}}(k_v)$.

\medskip
Assume  now that  $w_i$ and $w_j$ are ramified in $E_i^{w_i}$, respectively $E_j^{w_j}$, and that
${\rm char}(k_v) \not = 2$. In this case, $\sigma_N$ is the identity, and $W_{\kappa(N)}(k_v) \simeq
W(k_v)$.  Proposition \ref{residue} (b), implies that $\partial_v : T(E_i^{w_i},\sigma_i) \to W_{\kappa(N)}(k_v)$ and $\partial_v : T(E_j^{w_j},\sigma_j) \to W_{\kappa(N)}(k_v)$ are injective,
and their images consist of the classes of forms of dimension equal to $[\kappa(N):k_v]$ $\rm {mod} \ 2$.
The injectivity part of the statement implies that 
$\partial_v(\lambda_i^{w_i}) \not = \partial_v(\mu_i^{w_i})$
and
$\partial_v(\lambda_j^{w_j}) \not = \partial_v(\mu_j^{w_j})$. The forms 
$\partial_v(\lambda_i^{w_i}), \partial_v(\lambda_j^{w_j}),  \partial_v(\mu_i^{w_i}), \partial_v(\mu_j^{w_j})$
all have the same dimension $\rm {mod} \ 2$; therefore ${\rm det}(\partial_v(\lambda_i^{w_i})) \not = {\rm det}( \partial_v(\mu_i^{w_i}))$,
and  
${\rm det}(\partial_v(\lambda_j^{w_j})) \not = 
{\rm det}( \partial_v(\mu_j^{w_j}))$. Hence the forms 
$\partial_v(\lambda_i^{w_i}) \oplus \partial_v(\lambda_j^{w_j})$ and 
$ \partial_v(\mu_i^{w_i}) \oplus  \partial_v(\mu_j^{w_j})$ have the same dimension $\rm {mod} \ 2$ and the same
determinant, therefore they are equal in $W(k_v)$. This implies that $Q = 0$ in $W(k_v)$.

\medskip
Finally, assume that  $w_i$ and $w_j$ are ramified in $E_i^{w_i}$, respectively $E_j^{w_j}$, and that
${\rm char}(k_v) = 2$. Proposition \ref{residue} (c), implies that $\partial_v : T(E_i^{w_i},\sigma_i) \to W_{\kappa(N)}(k_v)$ and $\partial_v : T(E_i^{w_i},\sigma_i) \to W_{\kappa(N)}(k_v)$ are constant,
hence $\partial_v(\lambda_i^{w_i})  = \partial_v(\mu_i^{w_i})$
and
$\partial_v(\lambda_j^{w_j})  = \partial_v(\mu_j^{w_j})$. Since $W_{\kappa(N)}(k_v)$ is
of order $2$, this implies that $Q = 0$ in $W_{\kappa(N)}(k_v)$, hence in $W_{\Lambda_{k_v}}(k_v)$.

\medskip
This holds for any pair $w_i$, $w_j$ with the above properties, therefore the form
$$(\partial_v(\lambda_i^v)) \oplus (\partial_v(\lambda_j^v)) \oplus (- \partial_v(\mu_i^v)) \oplus (- \partial_v(\mu_j^v))$$ is trivial in $W_{\Lambda_{k_v}}(k_v)$; this completes the proof of the proposition.

\begin{notation}\label{cdelta} For all $v \in V_{\Sigma}$, let us fix $\delta = (\delta_v)$, with $\delta_v \in W_{\Lambda_{k_v}}(k_v)$ such that $\delta_v = 0$ for
almost all $v$. 
Let $\mathcal L_{\delta}^v$ be the set of $\lambda^v \in \mathcal {LR}^v$ such that $\partial_v(\lambda^v) = \delta_v$, and let 
$\mathcal C_{\delta}^v$ be the set of $a(\lambda^v) \in C(I)$ such that $\lambda^v \in \mathcal L_{\delta}^v$. 

\end{notation}

Let $\sim_v$  be the equivalence relation on $C(I)$ generated by
$$c \sim_v c' \iff c = (i,j)(c')$$
for some $i,j \in I$ with $v \in V_{i,j}.$

\medskip
\begin{coro}\label{nullcobordant coro} Let $v \in V_{\Sigma}$ and let $a(\lambda^v), a(\mu^v) \in 
\mathcal {CR}^v$ be such that $a(\mu^v) \sim_v a(\lambda^v)$.
If $a(\lambda^v) \in \mathcal C_{\delta}^v$, then $a(\mu^v) \in \mathcal C_{\delta}^v$.

\end{coro}

\noindent
{\bf Proof.} It suffices to show that if $ a(\mu^v) = (i,j)(a(\lambda^v))$ for some $i,j \in I$ with $v \in V_{i,j}$,
then $a(\mu^v) \in \mathcal C_{\delta}^v$; this follows from Proposition \ref{nullcobordant}.

\medskip
\begin{prop}\label{sim}  Let $v \in V_{\Sigma}$.  Then the set $\mathcal C_{\delta}^v$ is a $\sim_v$-equivalence class of $C(I)$.

\end{prop}

\medskip
\noindent
{\bf Proof.} Let us prove that if an element of $C(I)$ is $\sim_v$-equivalent to an element of $\mathcal C_{\delta}^v$,
then it is in $\mathcal C_{\delta}^v$. By Proposition \ref{utile}, this element belongs to $\mathcal {CR}^v$,
and Corollary \ref{nullcobordant coro} implies that it is in $\mathcal C_{\delta}^v$.

\medskip
Let $a(\lambda^v)$, $a(\mu^v) \in \mathcal C_{\delta}^v$, and let us show that $a(\lambda^v) \sim_v a(\mu^v)$. 
Let $J \subset I$ be the set of $i \in I$ such that If $a(\lambda^v)(i) \not = a(\mu^v)(i)$. Since
$$\underset{r \in I} \sum a(\mu^v)(r) = \underset{r \in I} \sum a(\lambda^v)(r) = A^v,$$ the set $J$ has an
even number of elements.

\medskip
Suppose first that $v$ is non-dyadic. This implies that if $i \in J$, then $\partial_v(\lambda_i^v) \not = \partial_v(\mu_i^v)$.
We have  $a(\lambda^v)$, $a(\mu^v) \in \mathcal C_{\delta}^v$ by hypothesis, hence 
$\partial_v(\lambda^v) = \partial_v(\mu^v)$, and therefore there exists $j \in J$ such that 
$\partial_v(W_{\Lambda_{K_v}}(K_v,M_i^v))$ and $\partial_v(W_{\Lambda_{K_v}}(K_v,M_j^v))$ 
have a non-zero intersection, and hence
$v \in V_{i,j}$. The map $(i,j)a(\lambda^v)$  differs from $a(\mu^v)$ in
less elements than $a(\lambda^v)$. 
Since $I$ is a finite set, continuing this way we see that $a(\lambda^v) \sim_v 
a(\mu^v)$. 

\medskip
Assume now that $v$ is dyadic, and let $J'$ be the set of $J$ such that $\partial_v(\lambda_i^v) = \partial_v(\mu_i^v)$.
Since $\partial_v(\lambda^v) = \partial_v(\mu^v)$, the set $J'$ also has an even number of elements. Let us write
$J = J' \cup J''$; then $J''$ has an even number of elements. If $i,j \in J'$, then $v \in V_{i,j}$. 
The map $(i,j)a(\lambda^v)$  differs from $a(\mu^v)$ in
less elements than $a(\lambda^v)$. After applying $(i,j)$ for all $j, j \in J'$ with $i \not = j$, we may assume that $J'$
is empty. Now we have $J = J''$, and the same argument as in the non-dyadic case shows that $a(\lambda^v) \sim_v 
a(\mu^v)$.

\section{Local-global problem - the integral case}\label{isometrieslattices}

We keep the notation of the previous sections.
%; in particular, $\Lambda$ is an $O$-algebra, 
%$A = \Lambda_K = \Lambda \otimes_O K$, and $M$ is an $A$-module such that 

\medskip Let $q$ be a quadratic form over $K$.
For all  $v \in V_{\Sigma}$, let us fix $\delta = (\delta_v)$, with $\delta_v \in W_{\Lambda_{k_v}}(k_v)$ such that $\delta_v = 0$ for
almost all $v$. Recall the following terminology from \S \ref{localglobalproblem} :

\medskip
We say that the {\it local conditions are satisfied} if conditions (L 1) and ${\rm (L  \ 2)}_{\delta}$
 below hold :

\medskip
(L 1) {\it For all $v \in V_K$, the quadratic form $(V,q)\otimes_K K_v$ is compatible with the module $M^v = M \otimes_K K_v$}

%\medskip If (L 1) holds, we obtain an $A^v$-quadratic form $(M^v,q^v)$ over $K_v$ for all $v \in V_K$, giving rise to an element of $W^b_{A}(K_v)$. 

%\medskip
\medskip
${\rm (L  \ 2)}_{\delta}$
 {\it For
all  $v \in V_{\Sigma}$, the quadratic form $(V,q) \otimes_K K_v$ contains an almost unimodular $\Lambda^v$-lattice with discriminant form $\delta_v$.}
%(L 2) For all finite places $v \in V_K$, we have $$\partial_v [M^v,q^v] = \delta_v.$$

\bigskip

We say that the {\it global conditions are satisfied} if conditions (G 1) and ${\rm (G \ 2)}_{\delta}$
 below hold :

\medskip
(G 1) {\it The quadratic form $(V,q)$ is compatible with the module $M$.}

%\medskip If (G 1) holds, we obtain an $A$-quadratic form $(M,q)$ over $K$, giving rise to an element of $W_{A}(K)$
%such 
%that  $(M^v,q^v)$ is a bounded $A^v$-quadratic form over $K_v$ for all finite places $v \in V_K$.

%\medskip
%(G 2) For all finite places $v \in V_K$, we have $$\partial_v [M^v,q^v] = \delta_v.$$

%\bigskip
%We say that the {\it Hasse principle holds}  if the local conditions imply the global conditions. 
\medskip
${\rm (G \ 2)}_{\delta}$
 {\it The quadratic form $(V,q)$ contains an almost unimodular $\Lambda$-lattice with discriminant form $\delta$.}

\begin{prop}\label{condition} Assume that the local conditions are satisfied. Then the
global conditions hold if and only if there exists $\lambda = (\lambda_i) \in T(E,\sigma)$ such
that $\lambda \in \mathcal L_{\delta}^v$ for all $v \in V_K$. 

\end{prop} 

\medskip
\noindent
{\bf Proof.} If  (G 1) holds, then by Proposition \ref{cond} here exists $\lambda = (\lambda_i) \in T(E,\sigma)$ such that 
$\lambda \in \mathcal {LR}^v$ or all $v \in V_K$. Condition ${\rm (G \ 2)}_{\delta}$
 implies that one can choose $\lambda$ such that  $\lambda \in \mathcal L_{\delta}^v$ for all $v \in V_{\Sigma}$.

\medskip 
Let us prove the converse.  
Let $\lambda = (\lambda_i) \in T(E,\sigma)$ such
that $\lambda \in \mathcal L_{\delta}^v$ for all $v \in V_{\Sigma}$. By Proposition \ref{cond}, this implies that (G 1) holds. Moreover, 
since $\lambda \in \mathcal L_{\delta}^v$ for all $v \in V_{\Sigma}$, we 
have $\partial_v(\lambda) = \delta_v$ for all $v \in V_{\Sigma}$, hence condition ${\rm (G \ 2)}_{\delta}$
 is also satisfied.

%$a \approx_v b$

\bigskip
Recall that the equivalence relation $\sim_v$ on  $C(I)$ is generated by
$$c \sim_v c' \iff c = (i,j)(c')$$
for some $i,j \in I$ with $v \in V_{i,j}.$

\bigskip
Let $\mathcal C_{\delta}$ the set of the elements $(a^v)$, $a^v \in \mathcal C_{\delta}^v$, such that
$a^v = 0$ for almost all $v \in V_{\Sigma}$. 

\medskip
\begin{prop}\label{integral almost0}

Assume that the local conditions are satisfied. Then
there exists $(\lambda^v) \in \mathcal L_{\delta}$ such that $(a(\lambda^v))
 \in \mathcal C_{\delta}$.

\end{prop}

\medskip
\noindent
{\bf Proof.} Let $S$ be the subset of $V_K$ consisting of the dyadic places, 
the infinite places,
the places that are ramified in $E_i/K$ for some $i \in I$, the places $v \in V_K$ such that $w_2(q) \not = w_2(q_{n,1})$ in ${\rm Br}_2(k_v)$, 
and the places $v \in V_{\Sigma}$ for which $\delta_v \not = 0$. Let us show that if $v \not \in S$, then there
exists $\lambda^v \in \mathcal L_{\delta}^v$ such that $a(\lambda^v) = 0$. 

\medskip Let $v \in V_K$ such that $v \not \in S$, and let $i \in I$. If $S_i^v = \varnothing$, then $T(E_i^v,\sigma_i) = 0$,
hence there is nothing to prove. Assume that $S_i^v \not= \varnothing$; recall that for all  $w \in S_i^v$ we have an
isomorphism $$\theta_i^w : T(E_i^w,\sigma_i) \to {\bf Z}/2{\bf Z},$$ and
for all $\lambda^v = (\lambda_i^w) \in T(E^v,\sigma)$, we have
 $$a(\lambda^v)(i) = \underset {w \in S_i^v} \sum \theta_i^w(\lambda_i^w).$$ 
 
 \medskip For all $i \in I$ such that $S_i^v \not= \varnothing$ and for all  $w \in S_i^v$, 
 let $\lambda_i^{w} = 1$ in $ T(E_i^w,\sigma_i)$.
 We claim that $a(\lambda_i^v) \in  \mathcal C_{\delta}^v$, and that $a(\lambda^v) = 0$. It is clear 
 that  $a(\lambda^v)(i) = 0$ for all $i \in I$; it remains to show that $(\lambda_i^v)  \in \mathcal L^v_{\delta}$.

 \medskip We first show that $(\lambda_i^v) \in {\mathcal LR}^v$. Since $v \not \in S$, we have 
 $w_2(q)  = w_2(q_{n,1})$ in ${\rm Br}_2(K_v)$. The quadratic forms $q$  is compatible with the module $M \otimes_K K_v$   by hypothesis, 
 hence ${\rm det}(q) = {\rm N}_{F/K}(-d)$ in $K_v^{\times}/K_v^{\times 2}$ (see Proposition
 \ref{det}). Since we also have ${\rm det}(q_{n,1}) = {\rm N}_{F/K}(-d)$, the quadratic forms $q$ and $q_{n,1}$ 
 have the same dimension, determinant and Hasse-Witt invariant over $K_v$, therefore they are isomorphic over $K_v$; this implies that
 $(\lambda_i^v) \in {\mathcal LR}^v$.
 
 \medskip

Since $v \not \in S$, it is unramified in $E_i$; hence Lemma \ref{unramified} implies that
$\partial_v(1) =  \partial_v(\lambda^v) = 0$; as  $\delta_v = 0$ for $v \not \in S$, we have $(\lambda_i^v)  \in \mathcal L^v_{\delta}$.

 \medskip
 
 {\bf Necessary and sufficient conditions} 
 
 \medskip
 Let $V'$ be the set of finite places of $K$.

\begin{prop}\label{homomorphism finite places} Let $(a(\lambda^v)), (a(\mu^v)) \in \mathcal C_{\delta}$, and let $c \in \sha_{\Lambda_M,M}$. Then we have

$$ \underset {v \in V'} \sum \  \underset{i \in I} \sum  \ c(i) a(\lambda^v)(i) = 
 \underset {v \in V'} \sum \  \underset{i \in I} \sum  \ c(i) a(\mu^v)(i).$$

\end{prop}

\medskip
\noindent
{\bf Proof.}
This follows from Corollary \ref{sim} and Proposition \ref{independent}. 

\bigskip

Let

$$\epsilon : \sha_{\Lambda_M,M} \to {\bf Z}/2{\bf Z}$$ 

\medskip
\noindent
be the homomorphism defined by

$$\epsilon(c) =  \underset {v \in V'} \sum \  \underset{i \in I} \sum  \ c(i) a(\lambda^v)(i)$$

\medskip
\noindent
for some $(a(\lambda^v)) \in \mathcal C_{\delta}$. By Proposition \ref{homomorphism finite places}, the homomorphism $\epsilon$
is independent of the choice of $(a(\lambda^v)) \in \mathcal C_{\delta}$.

\bigskip
Let $V''$ be the set of infinite places of $K$, and let $\mathcal  C(V'')$ be the set of $(a(\lambda^v))$ with $v \in V''$
and $a(\lambda^v) \in {\mathcal C}^v$. Note that since $V''$ is finite,
the set $\mathcal  C(V'')$ is also finite. If $V'' = \varnothing$, we set $\mathcal  C(V'') = 0$.

\medskip
For all $a \in \mathcal  C(V'')$ with $a = (a(\lambda^v))$, we define a homomorphism

 $$\epsilon_a :\sha_{\Lambda_M,M} \to {\bf Z}/2{\bf Z}$$ 

by setting 

$$\epsilon_a(c) =  \underset {v \in V''} \sum \  \underset{i \in I} \sum  \ c(i) a(\lambda^v)(i).$$

\medskip
\begin{theo}\label{necessary and sufficient} Assume that the local conditions hold. Then the global conditions are satisfied if
and only if there exists $a \in \mathcal  C(V'')$ such that $\epsilon + \epsilon_a = 0$. 

\end{theo}

\medskip
\noindent
{\bf Proof.} Assume that  the global conditions are satisfied. Then by Proposition \ref{condition} there exists 
$\lambda = (\lambda_i) \in T(E,\sigma)$ such
that $\lambda \in \mathcal L_{\delta}^v$ for all $v \in V_K$. 
We have
 $\underset {v \in V_K} \sum a(\lambda)(i) = 0$ for all $i \in I$. Set $a = (a(\lambda^v))$ for $v \in V''$; then we have  
 $\epsilon + \epsilon_a= 0$. 
 
 \medskip Let us prove the converse. 
Since the local conditions hold, Proposition \ref {integral almost0} implies that 
there exists $(\mu^v) \in \mathcal L_{\delta}$ such that $(a(\mu^v))
 \in \mathcal C_{\delta}$. 
By hypothesis, there exists $a \in \mathcal  C(V'')$ such that $\epsilon + \epsilon_a = 0$. Therefore Theorem \ref{sum}
implies that there exists $b = (b^v) \in \mathcal {C}_{\delta}$ such that 
  $\underset {v \in V_K} \sum b^v(i) = 0$. By definition, there exists $(\lambda^v) \in \mathcal {L}_{\delta}$ such
  that $b^v = a(\lambda^v)$ for all $v \in V_K$. 
 Recall that $a(\lambda^v)(i) = \underset {w \in S_v} \sum \theta_i^w(\lambda_i^w)$;
 therefore for all $i \in I$, we have  $\underset {w \in V_F} \sum \theta_i^w(\lambda_i^w) = 0$.
 By Theorem \ref{reciprocity} this implies that for all $i \in I$, there exists $\lambda_i \in T(E_i,\sigma_i)$ mapping to
 $\lambda_i^w \in T(E_i^w,\sigma_i)$ for all $w \in V_F$. 
 %In particular, we have $(\lambda_i,d_i) = (\lambda_i^v,d_i)$ for all $i \in I$ and $v \in V$.  
 Set
 $\lambda = (\lambda_i)$; we have $\lambda \in \mathcal L_{\delta}^v$ for all $v \in V_K$, hence by Proposition \ref{condition}  the global conditions hold. This completes
 the proof of the theorem. 
 
\begin{example}\label{isometry 5} Suppose  that $\Lambda = O[\Gamma]$, and let $f$ and $g$ be as in example \ref {isometry 4}.
%hence  $A = \Lambda \otimes_O K = K[\Gamma]$.
%With the notation of example \ref{isometry 3}, we have $f_i \in O[X]$; set $g = \underset {i \in I} \prod  f_i$, and $f = \underset {i \in I} \prod f^{n_i}$.
Assume that for all $v \in V_{\Sigma}$, the quadratic form $q$ contains a unimodular $O_v$-lattice having an isometry with minimal polynomial
$g$ and characteristic polynomial $f$; Theorem \ref{necessary and sufficient} gives a necessary and sufficient condition for
such a lattice to exists globally.

\end{example}

\bigskip
Recall that a quadratic form $q$ has maximal signature if for all real places $v \in V_K$, the signature
of $q$  at $v$ is equal to the signature of $q_{1,n}$ or of  $- q_{1,n}$ at $v$. 

\begin{lemma}\label{maximal} Assume that $q$ has maximal signature. Then $\mathcal C (V'')$ has at most one element. 

\end{lemma}

\noindent
{\bf Proof.} If $V''$ does not contain any real places, there is nothing to prove. Let $v$ be a real place of $K$, and let $(r_v,s_v)$ be the signature of $q$ at $v$. Let us assume that the signature of $q$ at $v$ is equal to the signature of $q_{1,n}$; the argument is the same if it is equal to 
the signature of $- q_{1,n}$. For all $i \in I$, let 
$(r_{n_i,1}^v, s_{n_i,1}^v)$ 
be the signature of $q_{n_i,1}$ and 
let $(r_{n,1},s_{n,1})$ be the signature of $q_{n,1}$. By hypothesis, 
we have $s_v = s_{n,1}^v$.
In all splittings of $q$ over $K_v$ into the orthogonal sum of quadratic forms $q_i^v$ over $K_v$ 
with signature $(r_i^v,s_i^v)$ compatible with the module $M_i \otimes_K K_v$, we have $s_i^v \geqslant s_{n_i,1}^v$ for all $i \in I$. 
Note that $s_{n,1}^v = \underset{i \in I} \sum s_{n_i,1}^v$
and $s_v =  \underset{ki\in I} \sum s_i^v$; therefore $s_i^v = s_{n_i,1}^v$ for all $i \in I$, 
and this implies that the local solution is unique. This holds for all real places $v \in V_K$,
hence $\mathcal C (V'')$ has at most one element.

\medskip
Assume that $q$ has maximal signature, and that the local conditions hold; let $a \in C(V'')$, and set $\epsilon' = \epsilon + \epsilon_a$. The following
is an immediate consequence of Theorem \ref{necessary and sufficient} :

\begin{coro}\label{definite} Assume that $q$ has maximal signature and that the local conditions hold. Then the global conditions hold if
and only if $\epsilon' = 0$. 

\end{coro}

%Recall that we say that the Hasse principle holds if and only if the local conditions imply the global conditions. 

 \begin{theo}\label{int} The Hasse principle holds for all quadratic forms $q$ and all $\delta = (\delta_v)$ if and only if $\sha_{\Lambda_M,M} = 0$.
 
 \end{theo}
 
\noindent
{\bf Proof.} If $\sha_{\Lambda_M,M}= 0$, then Theorem \ref{necessary and sufficient} implies that the Hasse principle holds for any $q$ and $\delta$. 
To prove the converse, assume that $\sha_{\Lambda_M,M} \not = 0$; we claim that there exists a quadratic form $q$ and
$\delta = (\delta_v)$ satisfying (L 1) and ${\rm (L \ 2)}_{\delta}$ but not (G 1) and ${\rm (G \ 2)}_{\delta}$.
Let
$c \in \sha_{\Lambda_M,M}$ be a non-trivial element, and let $i, j \in I$ be such that $c(i) \not = c(j)$. Let the quadratic form  $q$  and $a \in \mathcal {CR}$
be as in Proposition \ref{lemma}, and set $\delta_v = \partial_v[q]$ for all $v \in V_{\Sigma}$. Since $q$ is a global form, $\delta_v = 0$ for
almost all $v \in V_{\Sigma}$.
By construction, $q$ and $\delta$ satisfy conditions (L 1) and ${\rm (L  \ 2)}_{\delta}$.
Moreover, the form $q$ has maximal signature. We have

$$\epsilon'(c) = \underset {v \in V_K} \sum \  \underset{i \in I} \sum  \ c(i) a(\lambda^v)(i) = c(i) + c(j),$$
hence $\epsilon'(c) \not = 0$;  Corollary  \ref{definite} implies that  the global conditions are not satisfied.

\section {The integral case - Hasse principle with additional conditions}\label{additional}

Given two integers $r,s \geqslant 0$ and a polynomial $f \in {\bf Z}[X]$, does there exist an even, unimodular
lattice of signature $(r,s)$, and having an isometry with characteristic polynomial $f$ ? This is one of
the motivating questions of the paper. We can be more ambitious, and fix an element of ${\rm SO}_{r,s}({\bf R})$
with characteristic polynomial $f$; 
does it stabilize an even, unimodular lattice ?

\medskip The above question leads to a modified local-global problem : we impose the
local behaviour at the real places. 
As in the previous section, we fix a module $M$ and  a quadratic form  $(V,q)$. Moreover, for all
real places $v \in V_K$ we also fix an $A \otimes_K K_v$-quadratic form $(M \otimes_K K_v, b_v)$. 

\begin{example}\label{additional example} Let $A = K[\Gamma]$, $f \in K[X]$ and  $M$ be as in Example
\ref{isometry 3}. Fixing an $A \otimes_K K_v$-quadratic
form $(M \otimes_K K_v,b_v)$ for all real places $v \in V_K$ amounts to fixing a signature for all irreducible, symmetric
factors of $f \in K_v[X]$, as in \S \ref{signatures} (cf. Example \ref {iso signature}). 

%Example \ref{SO}. 

\end{example}

\medskip
We 
consider the following modified local and global conditions~:

\medskip Condition (L 1) is unchanged :

\medskip
(L 1) {\it For all $v \in V_K$, the quadratic form $(V,q)\otimes_K K_v$ is compatible with the module $M^v = M \otimes_K K_v$.}

\medskip
We have the new condition (L 1)$_{b_v}$ :

\medskip
(L 1)$_{b_v}$ {\it For all real places $v \in V_K$, there exists an isomorphism $$\varphi_v : V \otimes_K K_v \to M \otimes_K K_v$$
such that $(M \otimes_K K_v,b_{\varphi_v}) \simeq (M \otimes_K K_v, b_v)$.}

\medskip We also fix $\delta = (\delta_v)$, with $\delta_v \in W_{\Lambda_{k_v}}(k_v)$ such that $\delta_v = 0$ for
almost all $v$, and we have the (unchanged) condition :

\medskip
${\rm (L  \ 2)}_{\delta}$
 {\it For
all  $v \in V_{\Sigma}$, the quadratic form $(V,q) \otimes_K K_v$ contains an almost unimodular $\Lambda^v$-lattice with discriminant form $\delta_v$.}

\medskip
We say that the {\it local conditions} (L)$_{b_v,\delta}$ {\it are satisfied}  if the three conditions (L 1), (L 1)$_{b_v}$  and
${\rm (L  \ 2)}_{\delta}$ hold. 

\medskip
We say that the {\it global conditions {\rm (G)}$_{b_v,\delta}$ are satisfied if there exists an isomorphism 
of vector spaces $\varphi : V \to M$ such that 

\medskip
$\bullet$ $(M,b_\varphi)$ is an $A$-quadratic form;

\medskip
$\bullet$ For all real places $v \in V_K$, we have $(M \otimes_K K_v,b_{\varphi}) \simeq (M \otimes_K K_v, b_v)$:}

\medskip
$\bullet$ For all $v \in V_{\Sigma}$, we have $\partial_v(M \otimes_K K_v, b_{\varphi_v}) = \delta_v$.

\medskip
(In particular, conditions (G 1) and (G 2)$_{\delta}$ of the previous section hold).

\medskip Assume that the local conditions} (L)$_{b_v,\delta}$ are satisfied.
The obstruction group $\sha_{A_M,M}$ is as in \S \ref{isometrieslattices}, and so is the homomorphism $\epsilon : \sha_{A_M,M}
\to {\bf Z}/2{\bf Z}$; indeed, this homomorphism only depends on the finite places of $K$.

\medskip
If $v \in V_K$ is a real place, then condition (L 1)$_{b_v}$ determines the local data $\lambda^v$ uniquely; therefore
$C(V'')$ has exactly one element, namely $a = (a(\lambda^v))$. We define the homomorphism $\epsilon_a$ as
in \S \ref{isometrieslattices}. Setting  $\epsilon' = \epsilon + \epsilon_a$, we obtain a homomorphism 

$$\epsilon' : \sha_{A_M,M}
\to {\bf Z}/2{\bf Z}.$$

\begin{theo}\label{with additional} Assume that the local conditions 
{\rm (L)$_{b_v,\delta}$ } hold. Then the global conditions  {\rm (G)$_{b_v,\delta}$ }  are satisfied if
and only if $\epsilon'= 0$. 

\end{theo}

To prove this theorem, we need an analog of Proposition \ref {condition}. Let us denote by ${\mathcal L}^v_{b_v,\delta}$
the subset of ${\mathcal L}^v_{\delta}$ in which the components $\lambda^v$ for $v$ real are determined by
$(M \otimes_K K_v,b_v)$. 

\begin{prop}\label{condition additional} Assume that the local conditions {\rm (L)$_{b_v,\delta}$ }  are satisfied. Then the
global conditions {\rm (G)$_{b_v,\delta}$ } hold if and only if there exists $\lambda = (\lambda_i) \in T(E,\sigma)$ such
that $\lambda \in {\mathcal L}^v_{b_v,\delta}$
for all $v \in V_K$. 

\end{prop} 

\medskip
\noindent
{\bf Proof.} If  the global conditions {\rm (G)$_{b_v,\delta}$ } hold, then by Proposition \ref{cond} here exists $\lambda = (\lambda_i) \in T(E,\sigma)$ such that 
$\lambda \in \mathcal {LR}^v$ or all $v \in V_K$. Moreover, the conditions
 impliy that one can choose $\lambda$ such that  $\lambda \in {\mathcal L}^v_{b_v,\delta}$
for all $v \in V_K$. 

\medskip 
Let us prove the converse.  
Let $\lambda = (\lambda_i) \in T(E,\sigma)$ such that $\lambda \in {\mathcal L}^v_{b_v,\delta}$ for all $v \in V_K$.
By Proposition \ref{cond}, this implies that (G 1) holds. Moreover, 
since $\lambda \in {\mathcal L}^v_{b_v,\delta}$ for all $v \in V_K$, we 
have $\partial_v(M \otimes_K K_v, b_{\varphi_v}) = \delta_v$ for all $v \in V_{\Sigma}$, and 
$(M \otimes_K K_v,b_{\varphi}) \simeq (M \otimes_K K_v, b_v)$ for all real places $v \in V_K$, hence
the global conditions {\rm (G)$_{b_v,\delta}$ } hold.

\medskip
\noindent
{\bf Proof of Theorem \ref {with additional}.}  The proof goes along the lines of the one of Theorem \ref {necessary and
sufficient}, applying Proposition \ref {condition additional} instead of Proposition \ref {condition}.

\section{Lattices over $\bf Z$} \label{Z}

Let $f \in {\bf Z}[X]$ be a monic, symmetric polynomial without linear factors; we start by recalling from \cite {GM} some necessary conditions for the existence
of an even, unimodular lattice to have an isometry with characteristic polynomial $f$. We then apply the results of 
\S \ref{isometrieslattices} to give sufficient conditions as well. Set $2n = {\rm deg}(f)$.

\begin{defn}\label{c1}  We say that $f$ satisfies condition (C 1) if the integers  $|f(1)|$, $|f(-1)|$ and $(-1)^nf(1)f(-1)$ are all squares.

\end{defn}

Let $m(f)$ be the number of roots  $z$ of $f$ with $|z| > 1$ (counted with multiplicity). 

\begin{defn}\label{c2} Let $(r,s)$ be a pair of integers, $r, s \geqslant 0$. We say that condition (C 2) holds if $r+s = 2n$,   $r  \equiv s \ {\rm (mod \ 8)}$, 
$r \geqslant m(f)$, $s \geqslant m(f)$, and $m(f) \equiv r \equiv s  \ {\rm (mod \ 2)}$.

\end{defn}

The following lemma is well-known (see for instance \cite{GM}).

\begin{lemma}\label{c1c2}
Assume that there exists an even, unimodular lattice with signature $(r,s)$ having an isometry with
characteristic polynomial $f$. Then conditions  {\rm (C 1)} and  {\rm (C 2)} hold.

\end{lemma}

\medskip
\noindent
{\bf Proof.} It is clear that $r+s = 2n$, and the property $r  \equiv s \ {\rm (mod \ 8)}$ is well-known
(see for instance \cite{S}, Chapitre V, Th\'eor\`eme 2).
For the last part of condition (C 2), see   \cite{B4} Corollary 8.2 (in the case where $f$ is separable,
this is also proved in  \cite{GM} Corollary 2.3).  It is well-known that $(-1)^nf(1)f(-1)$ is a square,  see
for instance  \cite{B4} Corollary 5.2. The fact that $|f(1)|$ and $|f(-1)|$  are squares is proved in \cite {GM}, 
Theorem 6.1 (the hypothesis $f$ separable is not needed in the proof). 

\medskip
Let 

\centerline {$f = \underset {i \in I} \prod f_i^{n_i}$ \ \ \  and \ \ \ $g = \underset {i \in I} \prod f_i$}

\medskip 
\noindent
where $f_i \in {\bf Z}[X]$ are distinct irreducible, symmetric polynomials of even degree. Set $E_i = {\bf Q}[X]/(f_i) = {\bf Q}(\tau_i)$, and let
$\sigma_i : E_i \to E_i $ be the $\bf Q$-linear involution sending $\tau_i$ to $\tau_i^{-1}$; let $F_i$ be the fixed field of $E_i$.

\medskip
The following is essentially contained in \cite{BT}.
%Propositions 7.1 and 9.1.

\begin{theo}\label{c1} Assume that condition  {\rm (C 1)} holds. Then for each prime number $p$ there exists an even, unimodular ${\bf Z}_p$-lattice
having an isometry with characteristic polynomial $f$ and minimal polynomial $g$. 

\end{theo}

\noindent
{\bf Proof.} Let $E_{0,i}$ be an extension of degree $n_i$ of $F_i$, linearly disjoint from $E_i$. Set $\tilde E_i = E_i \otimes E_{0,i}$; then the characteristic
polynomial of the multiplication by $\tau_i$ on $\tilde E_i$ is $f_i^{n_i}$, and its minimal polynomial is $f_i$. Set $\tilde E = \underset{i \in I} \prod \tilde E_i$,  $\tilde E_0 = \underset{i \in I} \prod \tilde E_{0,i}$, 
and let $T : \tilde E \to \tilde E$ be the linear transformation acting on $\tilde E_i$ by multiplication with $\tau_i$; the characteristic polynomial of  $T$ is $f$, and
its minimal polynomial is $g$. The argument of \cite{BT}, proof of Theorem A, shows that for every prime number $p$ there exists a quadratic form on $\tilde E \otimes {\bf Q}_p$ containing an even, unimodular ${\bf Z}_p$-lattice stable by $T$.

\medskip
%\medskip
%Let $2n = {\rm deg}(f)$, and let $(r,s)$ be a pair of natural numbers   with $r + s = 2n$ and $r \ \equiv s \ {\rm (mod \ 8)}$. Set

%$$q = \langle 1,\dots,1,-1,\dots,-1 \rangle$$ considered as a quadratic form over ${\bf Q}$, of signature $(r,s)$. 

The following lemma is well-known :
 
 \begin{lemma}\label{fixed qform} {\rm (i)} Let $(r,s)$ be a pair of integers, $r,s \geqslant 0$,
 with $r \ \equiv s \ {\rm (mod \ 8)}$. Then
 there exists an even, unimodular lattice of signature $(r,s)$.
 
 \medskip
 {\rm (ii)} Two even, unimodular lattices of the same signature become isomorphic over ${\bf Q}$.  
  
 \end{lemma}

\noindent
{\bf Proof.} (i) Let $m =  {{r-s}\over 8}$; the orthogonal sum of $m$ copies of the $E_8$-lattice with $s$ hyperbolic planes has the
required property. 

\medskip
(ii) Let $q$ be a quadratic form over ${\bf Q}$ containing an even, unimodular lattice of signature $(r,s)$. Then
the dimension of $q$ is $r + s$, its signature is $(r,s)$ and its determinant is $(-1)^s$. The Hasse-Witt invariant of $q$ at
a prime $\not = 2$ is trivial (see for instance \cite{OM}, 92:1), and at infinity it is $0$ if $s \ \equiv 0 \ {\rm or}  \ 1 \ {\rm (mod \ 4)}$,
and $1$ if $s \ \equiv 2 \ {\rm or}  \ 3 \ {\rm (mod \ 4)}$. By reciprocity, this also determines the Hasse-Witt invariant of $q$
at the prime $2$ (one can also prove this directly,  see for instance \cite{BT}, Proposition 8.3). Therefore the dimension, determinant, signatures
and Hasse-Witt invariant of $q$ are uniquely determined by $(r,s)$; hence $q$ is unique up to isomorphism. 

\medskip
\noindent
{\bf Notation.}
If $(r,s)$ is a pair of integers  with $r,s \geqslant 0$ and $r \equiv s \ {\rm (mod \ 8)}$, let $(V,q_{r,s})$ be a quadratic form containing an even, unimodular
lattice of signature $(r,s)$; such a form exists by Lemma \ref{fixed qform} (i) and is unique up to isomorphism by Lemma \ref{fixed qform} (ii).

\medskip
Let $\Lambda = {\bf Z}[\Gamma]$, where $\Gamma$ is the infinite cyclic group, and set $A = {\bf Q}[\Gamma]$; let us denote by $\sigma$ the
involution of $\Lambda$ and $A$ sending $\gamma$ to $\gamma^{-1}$ for all $\gamma \in \Gamma$. Set $M_i = [{\bf Q}[X]/(f_i)]^{n_i}$
and $M = \underset {i \in I} \oplus M_i$. Let $(V,q)$ be a quadratic form over $\bf Q$. Recall from \S \ref{isometrieslattices} the local and global conditions (L 1), ${\rm (L  \ 2)}_{\delta}$
and (G 1), ${\rm (G \ 2)}_{\delta}$,
and
note that (for $\Lambda$ and $M$ as above, and for $\delta = 0$) they can be reformulated as follows

\medskip
(L 1) For all $v \in V_{\bf Q}$, the quadratic form $(V,q) \otimes_{\bf Q} {\bf Q}_v$ has an isometry with minimal polynomial $g$ and characteristic polynomial $f$. 

\medskip
(G 1) The quadratic form $(V,q)$ has an isometry with minimal polynomial $g$ and characteristic polynomial $f$. 

\medskip
(L 2) For all finite places $v \in V_{\bf Q}$, the quadratic form $(V,q) \otimes_{\bf Q} {\bf Q}_v$ has an isometry with minimal polynomial $g$ and characteristic polynomial $f$ that stabilizes a unimodular lattice of $V \otimes_{\bf Q} {\bf Q}_v$.

\medskip
(G 2) The quadratic form $(V,q)$ has an isometry with minimal polynomial $g$ and characteristic polynomial $f$ that stabilizes a
unimodular lattice of $V$.

\begin{theo}\label{qrs} Assume that conditions  {\rm (C 1)} and  {\rm (C 2)} hold. Then the local conditions {\rm (L 1)} and {\rm (L 2)} are satisfied
for the quadratic form $q_{r,s}$.

\end{theo}

\noindent
{\bf Proof.} Since condition (C 1) holds, Theorem \ref{c1}  implies that for each prime number $p$ there exists an even, unimodular ${\bf Z}_p$-lattice
having an isometry with characteristic polynomial $f$ and minimal polynomial $g$; we denote by $q^p$ the quadratic form over ${\bf Q}_p$ obtained
by extension of scalars from this lattice.   If $p \not=2$, this is the diagonal form $\langle 1,\dots,(-1)^s \rangle$ (see for instance \cite{OM}, 92:1); if $p = 2$, it is an orthogonal sum of
hyperbolic planes (see \cite{BT}, Proposition 8.3). This implies that $q^p \simeq q_{r,s} \otimes_{\bf Q}{\bf Q}_p$ for all prime numbers $p$. 
In particular, the quadratic form $q_{r,s}$ has an isometry with characteristic polynomial $f$ and minimal polynomial $g$ over ${\bf Q}_p$ for
every prime number $p$. This implies that the local condition (L 1) holds for every prime number $p$. Moreover, since 
$q_{r,s} \otimes_{\bf Q}{\bf Q}_p$ contains a unimodular  ${\bf Z}_p$-lattice stable by this isometry, 
%by Theorem \ref{BT} 
the local
condition ${\rm (L  \ 2)}_{\delta}$
 holds for $\delta = 0$.
%Therefore the local conditions (L 1) and (L 2) hold for every prime number $p$. 
On the other hand, condition (C 2) implies that  $q_{r,s} \otimes_{\bf Q} {\bf R}$ has an isometry with characteristic polynomial $f$ and minimal polynomial $g$ 
over $\bf R$ (cf. \cite{B4}, Corollary 8.2). 
Therefore the local conditions hold for the quadratic form $q_{r,s}$.

%Let $p$ be a prime number; as we have seen
%above,  $q_{r,s} \otimes_{\bf Q}{\bf Q}_p$ 
%contains a unimodular  ${\bf Z}_p$-lattice stable by the isometry.  
%By Theorem \ref{BT} this implies that the local condition (L 2) is satisfied for $\delta = 0$. 

%\medskip
%Let $\sha^{\rm int}_f$ be the group constructed in \S \ref{isometrieslattices}; we also refer to this section for the definitions of 
%the homomorphisms $\epsilon$ and $\epsilon_a$. 

\begin{prop}\label{GM and global} The following properties are equivalent :

\medskip

{\rm (i)} There exists an even, unimodular lattice of signature $(r,s)$ having an isometry with characteristic polynomial $f$ and
minimal polynomial $g$.

\medskip

{\rm (ii)} The global conditions ${\rm (G 1)}$ and ${\rm (G \ 2)}$
 are fulfilled for $(V,q_{r,s})$.

\end{prop}

\noindent
{\bf Proof.} Let us prove that {\rm (i)} implies {\rm (ii)}.  The base change to $\bf Q$ of the lattice is isomorphic to $(V,q_{r,s})$ by
Lemma \ref{fixed qform} (ii); hence $(V,q_{r,s})$ has an isometry with characteristic polynomial $f$ and minimal polynomial $g$, and
it contains a unimodular lattice stable by this isometry. This implies that the global conditions ${\rm (G 1)}$ and ${\rm (G \ 2)}_{\delta}$
are fulfilled for $(V,q_{r,s})$ and $\delta = 0$, and therefore (ii) holds. 

%hence condition (G 1) is satisfied, and gives rise to an element $[q_{r,s}]$ of $W_{A}(K)$. 
%By Theorem \ref{BT}, we have $\partial_p[q_{r,s}] = 0$ for all prime numbers $p$, therefore  condition (G 2) also holds. 

\medskip
Conversely, let us assume that {\rm (ii) holds. Let $t : V \to V$ be an isometry of $q_{r,s}$ with
characteristic polynomial $f$ and minimal polynomial $g$, and such that $(V,q_{r,s})\otimes_{\bf Q}{\bf Q}_p$ contains
%; this is possible by condition (G 1). Condition (G 2) and 
%Theorem \ref{BT} imply that $t : M \to M$ can be chosen in such a way that
%$(M,q_{r,s})\otimes_{\bf Q}{\bf Q}_p$ 
a unimodular  ${\bf Z}_p$-lattice $L_p$ with $t(L_p) = L_p$ for all prime numbers $p$. For $p = 2$, let us chose $L_2$ to be even; we claim that this is possible by \cite{BT}, Theorem 8.1. 
Indeed, condition (i) of Theorem 8.1 is satisfied, since $(V,q_{r,s})\otimes_{\bf Q}{\bf Q}_2$ contains a unimodular ${\bf Z}_2$-lattice; condition (ii) follows from the fact that 
$(V,q_{r,s})\otimes_{\bf Q}{\bf Q}_2$ is an orthogonal sum of hyperbolic planes; and condition (iii) holds by a result of Zassenhaus (see for instance \cite{BT}, Theorem 8.5),
since  $|f(-1)|$ is a square.
%Let $N$ be any lattice in $(M,q_{r,s})$ such that $t(N) = N$; then for almost all prime numbers $p$, we have $L_p = N \otimes_{\bf Z} {\bf Z}_p$. 
Let
$$L = \{x \in V \ | \ x \in L_p \ {\rm for \ all} \ p \}.$$ The lattice $L$ is even, unimodular, and $t(L) = L$; hence $L$ has an isometry with characteristic polynomial $f$ and minimal polynomial $g$. 

\bigskip
Recall from \S \ref{integral obstruction} the construction of the group $\sha_{\Lambda_M,M}$.  By Example \ref{isometry 4} the sets
$V_{i,j}$ consist of the prime numbers $p$ such that $f_i$ mod $p$ and $f_j$ mod $p$ have a common irreducible, symmetric factor. Hence
the group $\sha_{\Lambda_M,M}$ only depends on the polynomal $g$, and will be denoted by $\sha_g$. 

\medskip
We  next give some examples of groups $\sha_g$. For all integers $d \geqslant 1$, we denote by $\Phi_d$ the $d$-th cyclotomic polynomial. 

\begin{example}\label{Salem 10} Let $f_1(X) = X^{10} + X^9 - X^7 - X^6 - X^ 5 - X^4 - X^3 + X +1$, and $f_2 = \Phi_{14}$. Set $f = f_1 f_2^2$ and $g = f_1 f_2$. 

\medskip The resultant of $f_1$ and $f_2$ is 169, and the polynomials
$f_1$ mod 13 and $f_2$ mod 13 have the irreducible, symmetric common factor $X^2 + 7X + 1 \in {\bf F}_{13}[X]$. Therefore $V_{1,2} = \{ 13 \}$,
and $\sha_g = 0$.

\end{example}

\begin{example}\label{Salem 12} Let $f_1(X) = X^{12} - X^{11} + X^{10} - X^9 - X^6 - X^3 + X^2  - X + 1$, $f_2 = \Phi_{14}$ and $f_3 = \Phi_{12}$.
Set $f = g = f_1 f_2 f_3$.

\medskip
The resultant of $f_1$ and $f_2$ is 49, the resultant of $f_1$ and $f_3$ is 169, and the polynomials $f_2$ and $f_3$ are relatively prime.
The polynomials $f_1$ and $f_2$ mod 7 have the irreducible, symmetric common factor $X + 1 \in {\bf F}_7[X]$, and the polynomials
$f_1$ and $f_3$ mod 13 have two such common factors, $X + 2 \in {\bf F}_{13}[X]$ and $X + 7 \in {\bf F}_{13}[X]$. This implies that $V_{1,2} = \{ 7 \}$,
$V_{1,3} = \{ 13 \}$, and $V_{2,3} = \varnothing$, and hence $\sha_g = 0$. 

\end{example}

\begin{example}\label {cyclotomic 2} 
Let $p$ and $q$ be two distinct prime numbers, such that $p \equiv q \equiv 3 \ {\rm (mod \ 4})$. Let $n,m,t \in {\bf Z}$ with $n,m,t \geqslant 1$ and $m \not = t$, and set

$$f_1 = \Phi_{p^nq^m}, \ f_2=  \Phi_{p^nq^t}, \ {\rm and} \ \  f  = g = f_1f_2.$$

\bigskip

If ${\rm ({p \over q})} = 1$, then $V_{1,2} = \{q\}$, and $\sha_g = 0$. 

\medskip

If ${\rm ({p \over q})} =  - 1$, then $V_{1,2} = \varnothing$, and $\sha_g \simeq {\bf Z}/2{\bf Z}$. 

\end{example}

\begin{example}\label {intro 1} Let $p$ be a prime number with $p \equiv 3 \ {\rm (mod \ 4})$, and set $f_1 = \Phi_p$, $f_2 = \Phi_{2p}$,
$g = f_1 f_2$, and $f = f_1^2 f_2^2$.

\medskip

If $p \equiv  3 \ {\rm (mod \ 8})$, then $V_{1,2} = \{2\}$, and $\sha_g = 0$. 

\medskip

If $p \equiv  7 \ {\rm (mod \ 8})$, then $V_{1,2} = \varnothing$, and $\sha_g \simeq {\bf Z}/2{\bf Z}$.

\end{example}

\medskip More generally, the group $\sha_g$ can be  determined for any product of cyclotomic polynomials $\Phi_d$ with $d \geqslant 3$ (recall
that the polynomials we consider here don't have any linear factors, hence $\Phi_1$ and $\Phi_2$ are excluded); see \S \ref{torus}
for more details and examples. 

%\begin{prop} $sha_g =0$, and $sha_g^{\rm int} = 0 \ {\rm or} \ {\bf Z}/2{\bf Z}$. Moreover, $$sha_g^{\rm int} = 0 \iff {\rm ({p \over q})}$$
%\end{prop}

\bigskip
Assume now that conditions  {\rm (C 1)} and  {\rm (C 2)} hold; by \ref{qrs} the local conditions {\rm (L 1)} and {\rm ${\rm (L  \ 2)}_{\delta}$
} are satisfied
for the quadratic form $q_{r,s}$ and $\delta = 0$.
Hence $\mathcal C_0$ is not empty (where  $\mathcal C_0$ is the set  $\mathcal C_{\delta}$ for $\delta = 0$). 

\medskip
Recall from \S \ref{isometrieslattices} that $V'$ is the set of finite primes, and $V''$ the set of infinite primes; in our case,
$V'$ is the set of valuations $v_p$, where $p$ is a prime number, and $V'' = \{ v_{\infty} \}$ , where $v_{\infty}$ is the unique infinite place.
Since $\mathcal C_0$ is not empty, we define as in \S \ref{isometrieslattices} a homomorphism $\epsilon : \sha_{g} \to {\bf Z}/2{\bf Z}$; this homomorphism does
not depend of the choice of the local data $a \in \mathcal C_0$ (cf. Proposition \ref{homomorphism finite places}). 

\medskip
 Moreover, for all
$a \in \mathcal  C(V'')$, we have a homomorphism $\epsilon_a : \sha_g \to {\bf Z}/2{\bf Z}$ (cf. \S  \ref{isometrieslattices}); note that
$C(V'')$ is a finite set, hence we obtain a finite number of such homomorphisms.

\begin{theo}\label{final} Assume that conditions  {\rm (C 1)} and  {\rm (C 2)} hold. Then the following properties are equivalent :

\medskip
{\rm (i)} There exists an even, unimodular lattice of signature $(r,s)$ having an isometry with characteristic polynomial $f$ and
minimal polynomial $g$.

\medskip 
{\rm (ii)} There exists $a \in \mathcal  C(V'')$ such that $\epsilon + \epsilon_a = 0$. 

\end{theo}

\noindent
{\bf Proof.} %Let us prove that {\rm (i)} implies {\rm (ii)}. By Theorem \ref{c1}, {\rm (i)} implies that the global conditions hold for $\delta = 0$,
%hence by Theorem \ref{necessary and sufficient} property 
 %{\rm (ii)} holds.
 Since  conditions  {\rm (C 1)} and  {\rm (C 2)} hold, the local conditions {\rm (L 1)} and {\rm ${\rm (L  \ 2)}_{\delta}$
are satisfied
for the quadratic form $q_{r,s}$ and $\delta = 0$ (see Theorem \ref{qrs}); therefore we can apply
Theorem \ref{necessary and sufficient}, and obtain that (ii) is equivalent to

\medskip
{\rm (ii')} The global conditions ${\rm (G 1)}$ and ${\rm (G \ 2)}_{\delta}$
are fulfilled for $(V,q_{r,s})$ and $\delta = 0$.

\medskip By Proposition \ref{GM and global}, {\rm (ii')} and {\rm (i)} are equivalent. This completes the proof of the theorem.

\bigskip
If moreover the quadratic form $(V,q_{r,s})$ has maximal signature, we can define the homomorphism 
$\epsilon' : \sha_{g} \to {\bf Z}/2{\bf Z}$ and apply Corollary \ref{definite}. Note that $(V,q_{r,s})$ has maximal signature if
and only if $s = m(f)$ or $r = m(f)$. Hence we obtain the following 

\begin{coro}\label{definiteform}
Assume that conditions  {\rm (C 1)} and  {\rm (C 2)} hold, and that moreover $s = m(f)$ or $r = m(f)$. Then
there exists an even, unimodular lattice of signature $(r,s)$ having an isometry with characteristic polynomial $f$ and
minimal polynomial $g$ if and only if $\epsilon' = 0$. 

\end{coro}

\noindent
{\bf Proof.} This follows from Theorem \ref{final} and Corollary \ref{definite}.

\begin{coro}\label{sha 0} Assume that $f$ satisfies condition {\rm (C 1)}, and that $\sha_{g} = 0$. Then for each pair of integers $(r,s)$ such that condition {\rm (C 2)} holds
there exists an even, unimodular lattice with signature $(r,s)$ having an isometry with characteristic polynomial $f$ and
minimal polynomial $g$.

\end{coro}

\noindent
{\bf Proof.} This  follows from Theorem \ref{final}.

\begin{example}\label{cyclotomic 3} Let $f_1 = \Phi_{p^nq^m}, \ f_2=  \Phi_{p^nq^t}, \ {\rm and} \ \  f  = g = f_1f_2$ be as in
Example \ref{cyclotomic 2}. We have $f(1) = f(-1) = 1$ and $m(f) = 0$. Therefore conditions (C 1) and (C 2) are satisfied for all pairs of
integers $(r,s)$ with $r,s \geqslant 0$ such that $r  \equiv s \ {\rm (mod \ 8)}$ and  $r+s = {\rm deg}(f)$.  Since ${\rm deg}(f)$ is divisible by 8, the choice  $r = {\rm deg}(f)$ and $s = 0$ satisfies conditions (C 1) and (C 2). 
%the pair $({\rm deg}(f),0)$ satisfies the condition. 

\medskip
If ${\rm ({p \over q})} = 1$, then $\sha_g = 0$ (cf. Example \ref {cyclotomic 2}). Therefore by Corollary \ref{sha 0}, for all pairs $(r,s)$ as above,
there exists an even, unimodular lattice with signature $(r,s)$ having an isometry with characteristic polynomial $f$ (and
minimal polynomial $g = f$). In particular, there exists a {\it definite} even, unimodular lattice with this property.

\medskip

If ${\rm ({p \over q})} =  - 1$, then  $\sha_g \simeq {\bf Z}/2{\bf Z}$ (cf. Example \ref{cyclotomic 2}). Since conditions (C 1) and (C 2) are
satisfied (for a choice of $(r,s)$ as above), the local conditions are satisfied, and we have $\mathcal C_0 \not = \varnothing$. 

\medskip We denote by $v_p$  the finite place of $\bf Q$ corresponding to the prime number $p$, and by $v_{\infty}$ the unique real place. 
Let us choose $a \in C(V')$ as follows : $a^v(i) =0$ for all $i \in I = \{1,2 \}$  if $v \not = v_p$, and $a^{v_p} (1) = a^{v_p}(2) = 1$. 

\medskip
Let $c : I \to {\bf Z}/2{\bf Z}$ represent the unique non-trivial element of $\sha_g$ : we can take $c$ such that $c(1) = 1$ and $c(2) = 0$;
we have $\epsilon(c) = 1$.

\medskip Let $(r,s) = ({\rm deg}(f),0)$. Then $C(V'')$ has only one element, namely the identically zero one, $a = 0$, hence
 %$a = a^{v_{\infty}} : I \to {\bf Z}/2{\bf Z}$ with $a^{v_{\infty}} (1) = a^{v_{\infty}}(2) = 0$;  
 $\epsilon_a(c) = 0$.
 This implies that $\epsilon + \epsilon_a \not = 0$; therefore by Theorem \ref{final} (or Corollary \ref {definiteform}) there does not
exist any definite, even, unimodular lattice having an isometry with characteristic polynomial $f$.

\end{example}

\medskip The method of Example \ref{cyclotomic 3} can be used to decide which products of cyclotomic polynomials $\Phi_d$ (with $d \geqslant 3$)
occur as characteristic polynomials of isometries of definite even, unimodular lattices; this completes the results of \cite {B1}. 

\begin{example}\label{cyclotomic 4} With the notation of Example \ref{cyclotomic 3}, set $p = 3$, $q = 7$, $n = m = 1$ and $t = 2$. We have
$$f_1 = \Phi_{21}, \ f_2=  \Phi_{147}, \ {\rm and} \ \  f  = g = f_1f_2.$$

Since ${\rm ({3 \over 7})} =  - 1$, we have $\sha_g \simeq {\bf Z}/2{\bf Z}$ (see Example \ref {cyclotomic 2}). We have also seen (cf. Example
\ref{cyclotomic 3})  that
conditions (C 1) and (C 2) are satisfied for all pairs of
integers $(r,s)$ with $r,s \geqslant 0$ such that $r  \equiv s \ {\rm (mod \ 8)}$ and  $r+s = {\rm deg}(f) = 96$.

\medskip This gives rise to 25 possible pairs $(r,s)$. We have already seen that the signatures $(96,0)$ and $(0,96)$ are impossible (see Example \ref {cyclotomic 3}).
Assume that $(r,s) = (92,4)$.

\medskip The homomorphism $\epsilon : \sha_g \to {\bf Z}/2{\bf Z}$  is already computed in Example
\ref {cyclotomic 3}~: namely, if  $c : I \to {\bf Z}/2{\bf Z}$ represents the unique non-trivial element of $\sha_g$,
%(we can take $c$ such that $c(1) = 1$ and $c(2) = 0$)
we have $\epsilon(c) = 1$.

\medskip
The homomorphism $\epsilon_a$ depends on the choice of $a \in \mathcal C(V'') = \mathcal C^{v_\infty}$. There are two possibilities :
$a(1) = a(2) = 0$, and $a'(1) = a'(2) = 1$. Hence $\epsilon + \epsilon_a \not = 0$, but $\epsilon + \epsilon_{a'} = 0$; by Theoren \ref {final}, this
implies that there exists an even, unimodular lattice of signature $(92,4)$ having an isometry with characteristic polynomial~$f$. 

\medskip
It is easy to check that all the other signatures $(r,s)$, with the exception of (96,0) and (0,96), occur as signatures of even,
unimodular lattices having an isometry with characteristic polynomial $f$; this also follows from Proposition \ref {I 2} below.

%\medskip

%The following pairs are possible for $(r,s)$, with $r \geqslant s$  : $(96,0), (92,4), (88,8), (84,12), (80,16), (76,20), (72,24), (68,28), (64,28), 
\end{example}

\medskip
We now give some details concerning the counter-example to the Hasse principle of the introduction :

\medskip
\begin{example}\label{intro 2} Let $f_1 = \Phi_p$, $f_2 = \Phi_{2p}$,
$g = f_1 f_2$, and $f = f_1^2 f_2^2$ be as in Example \ref{intro 1}. We have $f(1) = f(-1) = p^2$, hence Condition (C 1)
holds. 

\medskip Let $p = 7$, and note that  $f$ is the polynomial
$$(X^6 + X^5 + X^4 + X^3 + X^2 + X + 1)^2(X^6 - X^5 + X^4 - X^3 + X^2 - X + 1)^2$$ of the introduction. We have
$m(f) = 0$, ${\rm deg}(f) = 24$, and condition (C 2) holds for $(r,s) = (24,0)$. This implies that the local conditions
are satisfied for $q_{24,0}$ and $\delta = 0$. We denote by $v_2$  the finite place of $\bf Q$ corresponding to the prime number $2$, and by $v_{\infty}$ the unique real place. 
Let us choose $a \in \mathcal C(V')$ as follows : $a^v(i) =0$ for all $i \in I = \{1,2 \}$  if $v \not = v_2$, and $a^{v_2} (1) = a^{v_2}(2) = 1$. 

\medskip By Example \ref{intro 1}, we know that $\sha_g \simeq {\bf Z}/2{\bf Z}$. 
Let $c : I \to {\bf Z}/2{\bf Z}$ represent the unique non-trivial element of $\sha_g$ : we can take $c$ such that $c(1) = 1$ and $c(2) = 0$;
we have $\epsilon(c) = 1$. On the other hand, $\mathcal C(V'')$ has only one element, $a = 0$, hence $\epsilon_a = 0$.
Therefore $\epsilon + \epsilon_a \not = 0$, and by Corollary \ref{definiteform} this implies that there does not
exist any positive definite, unimodular and even lattice having an isometry with characteristic polynomial $f$. 

\medskip
It is well-known that there exist exactly 24 isomorphism classes of even, unimodular, positive definite lattices of
rank 24, including the Leech lattice, and that they are isomorphic over ${\bf Z}_{\ell}$ for all prime
numbers $\ell$.  The above argument shows
that none of them has an isometry with characteristic polynomial $f$, and that they all have such an isometry 
locally everywhere.

\end{example}

\medskip

The following examples involve indefinite forms. Let $(r,s)$ be a pair of integers with $r,s \geqslant 1$ such that $r  \equiv s \ {\rm (mod \ 8)}$, and let $L_{r,s}$ be 
an  even, unimodular lattice  of signature $(r,s)$; it is well-known that such a lattice is unique up to isomorphism (see for instance \cite{S}, Chapitre V, Th\'eor\`eme 5).

\medskip
Recall that a {\it Salem polynomial}  is a monic, irreducible and symmetric polynomial in {\bf Z}[X] with exactly two roots outside the unit circle, both
real and positive.

\begin{example}\label{Salem 6} %As in \ref{GM}, Proposition 5.2
Let $f_1(X) = X^6 - 3 X^5   - X^4 + 5X^3 - X^2 -3X +1$ and $f_2 = \Phi_{12}$; set $f = g = f_1 f_2$. This example is taken from \cite {GM},
Proposition 5.2, where it is shown that $f$ does not arise as the characteristic polynomial of isometry of $L_{9,1}$. 
%Note that $f_1$ is a Salem
%polynomial, hence $m(f_1) = 1$. 
%an even, unimodular lattice with signature (9,1).

\medskip With the point of view of the present paper, this can be shown as follows. The polynomials $f_1$ and $f_2$ are relatively prime,
hence $V_{1,2} = \varnothing$, and this implies that $\sha_g \simeq {\bf Z}/2{\bf Z}$. 

\medskip

We have $f(1) = -1$, $f(-1) = 1$, and $m(f) = 1$; hence conditions (C 1) and (C 2) are fulfilled with $(r,s) = (9,1)$ or (1,9), and $(r,s) = (5,5)$. Let
$c : I \to {\bf Z}/2{\bf Z}$ represent the unique non-trivial element of $\sha_g$; let us chose $c$ so that $c(1) = 1$ and $c(2) = 0$. It is
easy to check that the homomorphism $\epsilon : \sha_g \to {\bf Z}/2{\bf Z}$ associated to the finite places satisfies $\epsilon (c) = 1$.

\medskip
Let $(r,s)$ = (9,1) or (1,9). In this case, the set $\mathcal C^{v_{\infty}}$ has a unique element, namely $a = 0$; hence $\epsilon_a = 0$, and $\epsilon + \epsilon_a \not = 0$. Therefore  by Theorem \ref {final}, $f$ does not arise as the characteristic polynomial of an isometry of $L_{9,1}$ or $L_{1,9}$.
%an even, unimodular lattice with signature (9,1) or (1,9). 

\medskip
Assume now that $(r,s) = (5,5)$. In this case, we can choose $a$ such that $\epsilon_a \not = 0$, and $\epsilon + \epsilon_a = 0$. Therefore
by Theorem \ref {final}, the lattice $L_{5,5}$ has an isometry with characteristic polynomial $f$.

\end{example}

\begin{prop}\label{I 2} Let $f = f_1^{n_1} f_2^{n_2}$ and let $(r,s)$ be a pair of integers with $r, s \geqslant 0$ such that $(r,s)$ is not a maximal signature for $f$. 
Set $g = f_1 f_2$. If conditions  {\rm (C~1)} and {\rm (C 2)} hold, then  $L_{r,s}$ has an isometry
with characteristic polynomial $f$ and minimal polynomial $g$. 

\end{prop}

\noindent
{\bf Proof.}  If $\sha_g = 0$, then this follows from Corollary \ref{sha 0}. Assume that $ \sha_g \not =0$; hence $\sha_g \simeq {\bf Z}/2{\bf Z}$. 
Let $I = \{0,1\}$, and let 
$c : I \to {\bf Z}/2{\bf Z}$ represent the unique non-trivial element of $\sha_g$; let us chose $c$ so that $c(1) = 1$ and $c(2) = 0$.  
Note that $\mathcal C^{v_{\infty}} =  \mathcal {CR}^{v_{\infty}}$.
We  claim that there exists 
$a \in \mathcal {CR}^{v_{\infty}}$ such that $\epsilon + \epsilon_a = 0$. Let $b \in  \mathcal {CR}^{v_{\infty}}$. If $b(1) = \epsilon (c)$,
set $a = b$. If $b(1) \not = \epsilon(c)$, then set $a = (1,2) b$; by Proposition \ref {real ij}, we
have $a \in \mathcal {CR}^{v_{\infty}}$. 
We have $(\epsilon + \epsilon_a)(c) =
\epsilon (c) + a(1)c(1) = 0$, hence $\epsilon + \epsilon_a = 0$. By Theorem \ref{final}, this implies that $L_{r,s}$ has an isometry
with characteristic polynomial $f$ and minimal polynomial $g$.

\begin{prop}\label{Salem times cyclotomic} Let $f = f_1 f_2$ where $f_1$ is a Salem polynomial, and $f_2$ is a power of a cyclotomic polynomial $\Phi_d$
with $d \geqslant 3$. Let $(r,s)$ be a pair of integers with $r,s \geqslant 3$, and assume that conditions {\rm (C 1)} and {\rm (C 2)} hold. 
Set $g = f_1 \Phi_d$.  Then $L_{r,s}$ has an isometry
with characteristic polynomial $f$ and minimal polynomial $g$. 

\end{prop} 

\noindent
{\bf Proof.} This follows from Proposition \ref {I 2}, since for $r,s \geqslant 3$, the signature $(r,s)$ is not maximal for $f$. 

%\medskip Set $E_1 = {\bf Q}[X]/(f_1)$ and $E_2 = {\bf Q}[X]/(\Phi_d)$. Let $\sigma_i : E_i \to E_i$ be the involution induced by $X \mapsto X^{-1}$,
%and let $F_i$ be the fixed field of the involution, for $i = 1,2$. 

%\noindent
%{\bf Remark.} As suggested by the above examples, one can show that when  $f = f_1^{n_1}  f_2^{n_2}$ with 
%$f_1, f_2 \in {\bf Z}[X]$ irreducible, symmetric, such that ${\rm deg}(f_i) > 2m(f_i)$ for $i = 1,2$ and the signature $(r,s)$ is not maximal, 
%then one can choose $a \in \mathcal C(V'')$ such that $\epsilon + \epsilon_a = 0$, hence the Hasse principle holds.

\bigskip

The following two examples might be interesting in the perspective of constructing K3 surfaces with given entropy (cf. McMullen \cite {Mc1}, \cite{Mc2},
\cite{Mc3}). 

\begin{example}\label{Salem 10 bis} As in Example \ref {Salem 10}, let $f_1(X) = X^{10} + X^9 - X^7 - X^6 - X^ 5 - X^4 - X^3 + X +1$, and $f_2 = \Phi_{14}$. Set $f = f_1 f_2^2$ and $g = f_1 f_2$.

%\medskip The resultant of $f_1$ and $f_2$ is 169, and the polynomials
%$f_1$ mod 13 and $f_2$ mod 13 have the irreducible, symmetric common factor $X^2 + 7X + 1 \in {\bf F}_{13}[X]$. Therefore $V_{1,2} = \{ 13 \}$,
%and $\sha_g = 0$.

\medskip Note that $f_1$ is
a Salem polynomial, hence $m(f_1) = 1$. 
We have $f(1) = -1$, $f(-1) = 49$, $m(f) = 1$. Therefore conditions (C 1) and (C 2) hold for all pairs of integers
$(r,s)$ with $r,s \geqslant 1$ such that $r  \equiv s \ {\rm (mod \ 8)}$ and  $r+s = {\rm deg}(f) = 22$.

\medskip
We have seen that  $\sha_g = 0$ (cf. Example \ref {Salem 10}). By Theorem \ref {final}  this implies that for all pairs $(r,s)$ as above, the lattice $L_{r,s}$ has an isometry with
%of signature $(r,s)$ having an isometry with 
characteristic polynomial $f$ and minimal polynomial $g$. For instance, the signature
$(r,s) = (3,19)$ is possible (the other possible signatures are (19,3), (15,7), (11,11) and (7,15)).

\end{example}

%The following examples involve indefinite forms. Let $(r,s)$ be a pair of integers with $r,s \geqslant 1$ such that $r  \equiv s \ {\rm (mod \ 8)}$, and let $L_{r,s}$ be 
%an  even, unimodular lattice  of signature $(r,s)$; it is well-known that such a lattice is unique up to isomorphism (see for instance \cite{S}, \S 5).

\begin{example}\label{Salem 12 bis} Let $f_1(X) = X^{12} - X^{11} + X^{10} - X^9 - X^6 - X^3 + X^2  - X + 1$, $f_2 = \Phi_{14}$ and $f_3 = \Phi_{12}$.
Set $f = g = f_1 f_2 f_3$.

\medskip
Since  $f_1$ is
a Salem polynomial, we have  $m(f_1) = 1$. Moreover,  $f(1) = -1$, $f(-1) = 49$, and  $m(f) = 1$, therefore conditions (C 1) and (C 2) hold for all pairs of integers
$(r,s)$ with $r,s \geqslant 1$ such that $r  \equiv s \ {\rm (mod \ 8)}$ and  $r+s = {\rm deg}(f) = 22$.

\medskip
\medskip
Since $\sha_g = 0$ (see Example \ref {Salem 12}),  Theorem \ref {final}  implies that for all pairs $(r,s)$ as above, the lattice $L_{r,s}$ has an isometry with
characteristic polynomial $f$. 
%there exists an even, unimodular lattice
%of signature $(r,s)$ having an isometry with characteristic polynomial $f$ and minimal polynomial $g$. 
The signature
$(r,s) = (3,19)$ is possible, as well the  signatures (19,3), (15,7), (11,11) and (7,15).

\end{example}

\section{Milnor signatures and Milnor indices}\label{Milnor}

In 1968, Milnor defined a signature invariant for knot cobordism (see \cite{M 68}, \S 5). The aim of this
section is to relate this invariant to the one defined in \S \ref {signatures} (cf. Example \ref {SO}).

\medskip
First, a question of terminology : the word ``signature" has two possible meanings, namely a pair $(r,s)$ or
the difference $r-s$. The convention of the present paper is to call signature the pair $(r,s)$, whereas in
knot theory the difference $r-s$ is used. To avoid confusion, we call {\it index} the difference $r - s$ (see
Notation \ref{index}). 

\medskip
Let $(V,q)$ be a non-degenerate quadratic form over $\bf R$, let $t : V \to V$ an isometry of $q$, and let $f \in {\bf R}[X]$
the characteristic polynomial of $t$. To each irreducible, symmetric factor $\mathcal P$ of $f$, Milnor associates an
index $\tau_{\mathcal P}$ (see \cite {M 68}, \S 5), as follows. Let $V_{\mathcal P(t)}$ be the $\mathcal P(t)$-primary subspace of $V$, consisting
of all $v \in V$ with $\mathcal P(t)^N v = 0$ for $N$ large. The {\it Milnor index} $\tau_\mathcal P (t)$ is by definition the index
of the restriction of $q$ to the subspace $V_{\mathcal P(t)}$. 

\medskip We define the {\it Milnor signature} at $\mathcal P$ as the signature of the restriction of $q$ to $V_{\mathcal P(t)}$. 
If the minimal polynomial of $t$ is square-free, these indices and signatures 
are the same as those defined in \S \ref{signatures} (see Example \ref{SO}). 

\medskip
More generally, we associate to $t : V \to V$ an isometry $t' : V \to V$ of $q$ with characteristic
polynomial $f$ and square-free minimal polynomial (see \cite {M}, \S 3). It follows from \cite {M}, Theorem 3.3 that
$\tau_{\mathcal P}(t') = \tau_{\mathcal P}(t)$ for all irreducible, symmetric factors $\mathcal P$ of $f$.

\section{Lattices and Milnor indices}\label{Z+}

We keep the notation of \S \ref{Z}. Let $f \in {\bf Z}[X]$ be a monic, symmetric polynomial of degree $2n$ 
without linear factors satisfying condition (C 1); let $r \geqslant 0$ and $s \geqslant 0$ be integers
such that condition (C 2) holds for $f$ and $(r,s)$. 

\medskip
Theorem \ref {final} gives a necessary and sufficient condition for the existence of an even, 
unimodular lattice of signature $(r,s)$ having an isometry with characteristic polynomial $f$ and  square-free
minimal polynomial. In this section, we ask a more precise question

\medskip
\noindent
{\bf Question.} Let $t \in {\rm SO}_{r,s}({\bf R})$ be an semi-simple isometry with characteristic polynomial $F$.
Does $t$ preserve an even, unimodular lattice ?

\medskip
Theorem \ref{final SO} below gives a necessary and sufficient condition for this to hold. We consider this
as a Hasse principle problem; Condition (C 1) implies that the local conditions hold at the finite places.
Condition (C 2) ensures the existence of a semi-simple element of ${\rm SO}_{r,s}(\bf R)$ with characteristic polynomial $f$; fixing such an element $t$ determines the local data at $\bf R$.
This is also the point of view of \S \ref{additional}, the results of which will be applied here.

\medskip
Let us write $f = \underset {i \in I} \prod f_i^{n_i}$, where $f_i \in {\bf Z}[X]$ are
distinct monic, irreducible, symmetric polynomials of even degree. Set $2n = {\rm deg}(f)$, and
$g = \underset {i \in I} \prod f_i$. For all $i \in I$, let $M_i = [{\bf Q}[X]/(f_i)]^{n_i}$,  and set $M = \underset{i \in I} \oplus M_i$. 

\medskip Recall that $V'$ is the set of finite places of $\bf Q$, and that $V''= \{v_{\infty} \}$, where $v_{\infty}$ is
the unique infinite place of $\bf Q$. 

\medskip We start by introducing some notation.

\begin{notation}\label{Mil} Let ${\rm Irr}_{\bf R}(f)$ be the set of irreducible, symmetric factors 
of $f \in {\bf R}[X]$. If $\mathcal P \in {\rm Irr}_{\bf R}(f)$, let  $n_{\mathcal P} > 0$ be the integer such that
$\mathcal P^{n_{\mathcal P}}$ is the power of $\mathcal P$ dividing $f$.

\medskip We denote by ${\rm Mil}(f)$ the set of maps $\tau : {\rm Irr}_{\bf R}(f) \to {\bf 2Z}$ such that the image of 
$\mathcal P \in {\rm Irr}_{\bf R}(f)$ belongs to the set $\{ -2n_{\mathcal P},\dots,2n_{\mathcal P} \}$. Let
${\rm Mil}_{r,s}(f)$ be the subset of ${\rm Mil}(f)$ such that $\underset {\mathcal P} \sum \tau_{\mathcal P} = r-s$,
where the sum runs over $\mathcal P \in {\rm Irr}_{\bf R}(f)$.

\end{notation}

\begin{defn}\label{Milnorsignature} An element of ${\rm Mil}(f)$ is called a {\it Milnor index}. 

\end{defn}

\medskip
Recall from the previous sections that the following sets are in bijective correspondance

\medskip
(a) Conjugacy classes of semi-simple elements of ${\rm SO}_{r,s}(\bf R)$ with characteristic polynomial $f$.

\medskip
(b) Isomorphism classes of ${\bf R}[\Gamma]$-quadratic forms on $M \otimes_{\bf Q} {\bf R}$ of signature $(r,s)$.

\medskip
(c) ${\rm Mil}_{r,s}(f)$.

\medskip  To a semi-simple isometry $t \in {\rm SO}_{r,s}(\bf R)$ with characteristic polynomial $f$, the 
above bijection associates 
an ${\bf R}[\Gamma]$-quadratic form $b(t)$ on the module $M \otimes_{\bf Q} {\bf R}$.
 With the notation of
\S \ref{additional} we have 

\begin{prop}\label{SO and global}  Let $t \in {\rm SO}_{r,s}(\bf R)$ be a semi-simple isometry with
characteristic polynomial $f$. The following properties are equivalent :

\medskip

{\rm (i)} The isometry $t$ preserves an even, unimodular lattice.

\medskip

{\rm (ii)} The global condition ${\rm (G)}_{b(t),0}$ is
 fulfilled for $(V,q_{r,s})$.

\end{prop}

\noindent
{\bf Proof.} By Proposition \ref{GM and global}, the existence of an even, unimodular lattice having an
isometry with characteristic polynomial $f$ and square-free minimal polynomial implies conditions
(G 1) and (G 2)$_0$.  Since $t \in {\rm SO}_{r,s}(\bf R)$ preserves the lattice, ${\rm (G)}_{b(t),0}$ holds; therefore
(i) implies (ii). Conversely, let us show that (ii) implies (i). By Proposition \ref{GM and global}, there exists
an even, unimodular lattice having an
isometry with characteristic polynomial $f$ and square-free minimal polynomial; condition ${\rm (G)}_{b(t),0}$ 
implies that it is preserved by $t \in {\rm SO}_{r,s}$.

\bigskip As in \S \ref{Z}, we obtain a homomorphism
$$\epsilon : \sha_g \to {\bf Z}/2{\bf Z},$$ defined in terms of finite places, that is, $C(V')$. Moreover,  to all elements $a \in C(V'')$, we associate a homomorphism 
$$\epsilon_a : \sha_g \to {\bf Z}/2{\bf Z}.$$

Recall from \S \ref{quadr} that  an ${\bf R}[\Gamma]$-quadratic form on $M \otimes_{\bf Q} {\bf R}$ gives rise to an element $\lambda^{v_{\infty}} \in \mathcal {LR}^{v_{\infty}}$, and hence also to $a(\lambda^{v_{\infty}}) \in C(V'')$. 
The above bijection allows us to associate 

\medskip
$\bullet$ to $t \in {\rm SO}_{r,s}(\bf R)$, an element $a_t \in C(V'')$;

\medskip

$\bullet$ to $\tau \in {\rm Mil}_{r,s}(f)$, an element $a_{\tau} \in C(V'')$.

\medskip
\begin{theo}\label{final SO} Let $t \in {\rm SO}_{r,s}(\bf R)$ be a semi-simple isometry with characteristic polynomial $f$.
The following are equivalent

\medskip
{\rm (a)} The isometry $t$ preserves an even, unimodular lattice.

\medskip
{\rm (b)} $\epsilon + \epsilon_{a_t} = 0$. 

\end{theo}

\noindent
{\bf Proof.} By Proposition \ref {SO and global}, property (a) holds if and only if the global condition ${\rm (G)}_{b(t),0}$ is
 fulfilled for $(V,q_{r,s})$. On the other hand, this is equivalent to (b) by Theorem \ref{with additional}. This concludes
 the proof of the theorem.
 
 \begin{coro}\label{all SO}  If $\sha_g = 0$, then all semi-simple elements of ${\rm SO}_{r,s}(\bf R)$ with characteristic
 polynomial $f$ preserve an even, unimodular lattice.
 
 \end{coro}
 
 \noindent
 {\bf Proof.} This is an immediate consequence of Theorem \ref{final SO}.
 
 \medskip Note that this is a generalization of Theorem 1.3 of \cite{GM}; indeed, if $f$ is irreducible, then $\sha_g = 0$. 

\bigskip
If $L$ is a lattice of signature $(r,s)$ having an isometry $t$ with characteristic polynomial $f$ and minimal
polynomial $g$, then $t$ extends to a semi-simple element of ${\rm SO}_{r,s}(\bf R)$, and this
element gives rise to an element of ${\rm Mil}_{r,s}(f)$; this element will be called the {\it Milnor index} of
the pair $(L,t)$. 

\medskip
The following results are  reformulations of Theorem \ref{final SO} and Corollary \ref{all SO}.

\begin{theo}\label{final Mil} Let $\tau \in {\rm Mil}_{r,s}(f)$. 
The following are equivalent

\medskip
{\rm (a)} There exists  an even, unimodular lattice having an isometry with characteristic polynomial $f$ and
Milnor index $\tau$. 

\medskip
{\rm (b)} $\epsilon + \epsilon_{a_\tau} = 0$. 

\end{theo}

\begin{coro}\label{all Milnor} Assume that $\sha_g = 0$. Then every $\tau \in {\rm Mil}_{r,s}(f)$ occurs as the Milnor
index of an even, unimodular lattice with an isometry of characteristic polynomial $f$. 

\end{coro}

As shown by the following example, if $\sha_g \not = 0$, then Theorem   \ref{final Mil} allows us to determine which Milnor indices occur.

\begin{example}\label{cyclotomic Milnor} As in Example \ref{cyclotomic 4}, let 
$$f_1 = \Phi_{21}, \ f_2=  \Phi_{147}, \ {\rm and} \ \  f  = g = f_1f_2.$$
We have already seen that 
conditions (C 1) and (C 2) are satisfied for all pairs of
integers $(r,s)$ with $r,s \geqslant 0$ such that $r  \equiv s \ {\rm (mod \ 8)}$ and  $r+s = {\rm deg}(f) = 96$,
and that $\sha_g \simeq {\bf Z}/2{\bf Z}$; let $c \in \sha_g$ be the only non-trivial element. 

\medskip
The homomorphism $\epsilon : \sha_g \to {\bf Z}/2{\bf Z}$  satisfies $\epsilon (c) =1$. The homomorphism $\epsilon_a$ depends on the choice of $a \in \mathcal C(V'') = \mathcal C^{v_\infty}$. There are two possibilities :
$a(1) = a(2) = 0$, and $a'(1) = a'(2) = 1$. Hence $\epsilon + \epsilon_a \not = 0$, and $\epsilon + \epsilon_{a'} = 0$.

\medskip By Theorem \ref{final Mil}, there exists  an even, unimodular lattice having an isometry with characteristic polynomial $f$ and
Milnor index $\tau$ if and only if $a_{\tau} = a'$. 

\medskip Let $r,s \geqslant 0$ be two integers such that  $r  \equiv s \ {\rm (mod \ 8)}$ and  $r+s = {\rm deg}(f) = 96$.
Since $n_{\mathcal P} = 2$ for all $\mathcal P \in {\rm Irr}_{\bf R}(f)$, the set ${\rm Mil}_{r,s}(f)$  
consists of the maps $$\tau : {\rm Irr}_{\bf R}(f) \to \{ -2, 2 \}.$$ with $\underset {\mathcal P} \sum \tau_{\mathcal P} = r-s$,
where the sum runs over $\mathcal P \in {\rm Irr}_{\bf R}(f)$.

\medskip The Milnor index map $\tau$ is determined by the two maps

$$\tau_1 : {\rm Irr}_{\bf R}(f_1) \to \{ -2, 2 \} \ \ {\rm and} \ \ \tau_2 : {\rm Irr}_{\bf R}(f_2) \to \{ -2, 2 \}.$$

\medskip

For $i = 1,2$,  let $N_i$ be the number of $\mathcal P \in {\rm Irr}_{\bf R}(f_i)$ such that $\tau_i(\mathcal P) = -2$. Since
 $r  \equiv s \ {\rm (mod \ 8)}$ and  $r+s = 96 \equiv 0 \ {\rm (mod \ 8)}$, we have $s \equiv 0 \ {\rm (mod \ 4)}$. Therefore 
$N_1 + N_2$ is even, hence $N_1 \equiv N_2 \ {\rm (mod \ 2)}$. 

\medskip
We have $a_{\tau} = a$ if $N_1$ is even, and $a_{\tau} = a'$ if $N_1$ is odd. Set $N(\tau) = N_1\  \ {\rm (mod \ 2)}$.
In summary, we obtain :

\medskip

There exists  an even, unimodular lattice having an isometry with characteristic polynomial $f$ and
Milnor index $\tau$ if and only $N(\tau) = 1$.

\end{example}

\section{Knots}\label{knots}

A {\it knot} is a smooth, oriented submanifold of $S^3$, homeomorphic to $S^1$. One of the classical knot
invariants is the {\it Alexander polynomial}, another is the {\it signature}. 
%The former is a symmetric polynomial
%$\Delta \in {\bf Z}[X]$ such that $\Delta(1) =  1$; the later is a pair of integers $r, s \geqslant 0$ with
%$r + s = {\rm deg}(\Delta)$. 

\medskip
The results of the previous sections can be applied to decide for which pair $(r,s)$ there exists
a knot with  Alexander polynomial $\Delta$
and signature $(r,s)$. For simplicity, we restrict to monic
polynomials $\Delta$ such that $\Delta(-1) = \pm 1$; these polynomials will be called {\it unramified}. Moreover,
we assume that $\Delta$ is a product of distinct irreducible, symmetric polynomials. The general case will be treated elsewhere.

\medskip
Still assuming that $\Delta$ is unramified and square free, we deal with a
more precise question. To each irreducible, symmetric factor $\mathcal P$ of $\Delta \in {\bf R}[X]$, one associates
a {\it Milnor signature} $\sigma_{\mathcal P} = (2,0)$ or $(0,2)$.  Given $\Delta$ and $\sigma_{\mathcal P} $ as above, we give a necessary
and sufficient criterion for the existence of  a knot with Alexander polynomial $\Delta$
and Milnor signatures $\sigma_{\mathcal P}$.

\medskip
As usual
in knot theory, the results are expressed in terms of the {\it index} $r-s$ rather than the signature $(r,s)$, and the {\it Milnor indices} $\tau_{\mathcal P}$
rather than the Milnor signatures $\sigma_{\mathcal P}$.

\medskip
For background information on the various notions of knot signatures (that is, indices), see the survey
paper of Ghys and Ranicki \cite{GR}, explaining the history of the topic and its connection to other aspects of 
geometry and topology; see also  the more recent survey of Anthony Conway  \cite {Co}.

\medskip
We start by recalling some facts on Seifert forms, and then come back to the applications to knot theory.

%\medskip
%We start by recalling some basic notions of knot theory.

%\medskip
%{\bf Seifert surfaces and Seifert forms}

%\medskip
%Any knot is the boundary of an oriented surface in $S^3$; such a surface is called a {\it Seifert surface} of the knot.

\section{Seifert forms}\label{Seifert}

A {\it Seifert form} is by definition a $\bf Z$-bilinear form $\mathcal A : \mathcal L \times \mathcal L \to {\bf Z}$, where $\mathcal L$ is
a free $\bf Z$-module of finite rank, such that the skew-symmetric form $\mathcal I : \mathcal L \times \mathcal L \to {\bf Z}$
given by $$\mathcal I (x,y) = \mathcal A(x,y) -  \mathcal A (y,x)$$ has determinant 1. The symmetric form
$\mathcal Q : \mathcal L \times \mathcal L \to {\bf Z}$ defined by $$\mathcal Q (x,y) = \mathcal A(x,y) +  \mathcal A (y,x)$$ 
is called the {\it quadratic form} associated to the Seifert form $\mathcal A$.

\medskip
The {\it index} (or signature) of $\mathcal A$ is by definition the index  (or signature) of the quadratic form $\mathcal Q$ (see Trotter 
\cite {T 62}, Proposition 5.1, or Murasugi \cite{Mu}). 

\medskip The {\it Alexander polynomial} of $\mathcal A$, denoted by $\Delta_{\mathcal A}$, is by definition the
determinant of the form $\mathcal L \times \mathcal L \to {\bf Z}[X]$ given by

$$(x,y) \mapsto \mathcal A(x,y) X - \mathcal A(y,x).$$ Note that $\Delta_{\mathcal A} (1) = 1$, that
$\Delta_{\mathcal A} (0) = {\rm det}(\mathcal A)$, and that $\Delta_{\mathcal A} (-1) = {\rm det}(\mathcal Q)$.

\medskip
We say that a polynomial $\Delta \in {\bf Z}[X]$ is {\it unramified} if $\Delta$ is monic, $\Delta(1) = 1$, and
$\Delta(-1) = \pm 1$. In the sequel, we only consider Seifert forms with unramified Alexander polynomials. 
Note that this implies that ${\rm det}(\mathcal A) = 1$ and ${\rm det}(\mathcal Q) = \pm 1$; in particular,
$\mathcal Q$ is a unimodular lattice. Let $\mathcal T : \mathcal L \to \mathcal L$ be
defined by $\mathcal A (\mathcal T (x), y) = \mathcal A (y,x)$. We have $\mathcal Q(\mathcal T (x), \mathcal T(y)) =
\mathcal Q(x,y)$
for all $x,y \in \mathcal L$, in other words, $\mathcal T$ is an isometry of $\mathcal Q$. Note that the
characteristic polynomial of $\mathcal T$ is equal to $\Delta_{\mathcal A}$. 

\medskip
Assume in addition that $\Delta_{\mathcal A}$ is a product of distinct irreducible, symmetric
polynomials of ${\bf Z}[X]$. 

\medskip
For all irreducible, symmetric factors $\mathcal P \in {\bf R}[X]$ of $\Delta_{\mathcal A}$, we define
a {\it Milnor index} (and Milnor signature) of $(\mathcal Q,\mathcal T)$ as in \cite {M 68}, \S 5;  
see also \S \ref {Milnor}. These will be called Milnor indices (and signatures) of $\mathcal A$. 

\medskip
Given a pair $(\mathcal Q,\mathcal T)$ consisting of a unimodular lattice $\mathcal Q$ and an isometry $\mathcal T$
of characteristic polynomial $\Delta$,
we recover a Seifert form $\mathcal A$ with $\mathcal Q (x,y) = \mathcal A(x,y) +  \mathcal A (y,x)$ for all $x, y \in 
\mathcal L$ (see for instance Levine \cite {Le 69}, \S 9, for a similar argument). 

\medskip Recall from  \S \ref{Z+}, Notation \ref{Mil}, that 
${\rm Irr}_{\bf R}(\Delta_{\mathcal A})$ is the set of irreducible, symmetric factors 
of $\Delta_{\mathcal A} \in {\bf R}[X]$. Since $\Delta_{\mathcal A}$ is square-free, we  have $n_{\mathcal P} = 1$ for 
all $\mathcal P$ in ${\rm Irr}_{\bf R}(\Delta_{\mathcal A})$. Moreover,
recall that we denote by ${\rm Mil}(\Delta_{\mathcal A})$ the set of maps $\tau : {\rm Irr}_{\bf R}(\Delta_{\mathcal A}) \to \{-2,2 \}$, and by
${\rm Mil}_{r,s}(\Delta_{\mathcal A})$ the subset of ${\rm Mil}(\Delta_{\mathcal A})$ such that $\underset {\mathcal P} \sum \tau_{\mathcal P} = r-s$,
where the sum runs over $\mathcal P \in {\rm Irr}_{\bf R}(\Delta_{\mathcal A})$.

Hence we have the following

\begin{prop}\label{lattice and Seifert} Let $r \geqslant 0$ and $s \geqslant 0$ be two integers, and let
$\Delta \in {\bf Z}[X]$ be an unramified polynomial that is a product of distinct irreducible, symmetric polynomials
of ${\bf Z}[X]$. Let $\tau \in {\rm Mil}_{r,s}(\Delta)$. The following are equivalent 

\medskip
{\rm (a)} There exists a Seifert form with Alexander polynomial $\Delta$ and Milnor index $\tau$.

\medskip
{\rm (b)} There exists an even, unimodular lattice $L$ having an isometry $t$ of characteristic polynomial $\Delta$
such that the Milnor index of $(L,t)$ is  $\tau$.

\end{prop} 

In particular, there exists a Seifert form with Alexander polynomial $\Delta$ and index $r-s$ if and only
if there exists an even, unimodular lattice of index $r-s$ having an isometry with characteristic polynomial $\Delta$.

\section{Knots and Seifert forms}\label{knots and Seifert}

We keep the notation of Sections \ref{knots} and \ref {Seifert}. To every knot $\Sigma$  in $S^3$, we associate
 a Seifert form $\mathcal A_{\Sigma}$ (see for instance \cite {Ka}, Chapter VII, \cite{Lick}, Definition 6.4, \cite{Liv 93}, Chapter 6, \S 1). It is well-known that the Alexander polynomial
of $\mathcal A_{\Sigma}$
 is an
invariant of the knot; the same thing is true for  the index (see Trotter, \cite {T 62}, Proposition 5.1) and the Milnor indices 
(see Milnor, \cite {M 68}, \S 5). 

\medskip
The following result of Seifert shows that all Seifert forms are realized by knots (see \cite {Se}). 

\begin{theo} \label {Seifert's theorem} Let $\mathcal A$ be a Seifert form; then there exists a
knot $\Sigma$ in $S^3$ such that $\mathcal A$ is the Seifert form of $\Sigma$. 

\end{theo}

Combining this result with Theorem \ref{lattice and Seifert}, we obtain the following

\begin{prop}\label{lattice and knots} Let $r \geqslant 0$ and $s \geqslant 0$ be two integers, and let
$\Delta \in {\bf Z}[X]$ be an unramified polynomial that is a product of distinct irreducible, symmetric polynomials
of ${\bf Z}[X]$. Let $\tau \in {\rm Mil}_{r,s}(\Delta)$. The following are equivalent 

\medskip
{\rm (a)} There exists a knot with Alexander polynomial $\Delta$ and Milnor index $\tau$.

\medskip
{\rm (b)} There exists an even, unimodular lattice $L$ having an isometry $t$ of characteristic polynomial $\Delta$
such that the Milnor index of $(L,t)$ is  $\tau$.

\end{prop} 

In particular, there exists a knot with Alexander polynomial $\Delta$ and index $r-s$ if and only
if there exists an even, unimodular lattice of index $r-s$ having an isometry with characteristic polynomial $\Delta$.
Using this result, we can apply Theorem \ref{final Mil} to answer the questions of \S \ref{knots}. 

\bigskip
We start by recalling some definitions from \S \ref{Z} and \S \ref{Z+}. Let $\Delta \in {\bf Z}[X]$ be an unramified polynomial as above, and set $2n = {\rm deg}(\Delta)$.

\bigskip
{\bf The conditions (C 1) and (C 2)}

\medskip Recall from \S \ref{Z} that the local conditions for the existence of an even, unimodular lattice of signature
$(r,s)$ and characteristic polynomial $\Delta$ can be translated into two conditions (C 1) and (C 2). We now recall
these conditions in our situation. 

\medskip 
$\bullet$ Since $\Delta$ is unramified, condition (C 1) becomes : $\Delta(-1) = (-1)^n$. 

\medskip
Recall that $m(\Delta)$ is the number of roots  $z$ of $\Delta$ with $|z| > 1$.

\medskip
$\bullet$ Let $(r,s)$ be a pair of integers, $r, s \geqslant 0$. Condition (C 2) holds if $r+s = 2n$,   $r  \equiv s \ {\rm (mod \ 8)}$, 
$r \geqslant m(\Delta)$, $s \geqslant m(\Delta)$, and $m(\Delta) \equiv r \equiv s  \ {\rm (mod \ 2)}$.

\medskip
It follows from Lemma \ref{c1c2} that conditions (C 1) and (C 2) are necessary for the existence of a knot with Alexander
polynomial $\Delta$ and index $r - s$. 
%[Note that  $r  \equiv s \ {\rm (mod \ 8)}$ and $r + s = 2n$ imply  $ n  \equiv s \ {\rm (mod \ 2)}$, hence condition (C 2) implies condition (C 1).]

\begin{example}\label {c1c2 torus} Assume moreover that $\Delta$ is a product of cyclotomic polynomials $\phi_m$ with $m \geq 3$. 
Then $m(\Delta) = 0$, hence if  $r+s = 2n, \ r  \equiv s \ {\rm (mod \ 8)}$, then condition (C 2) holds. 
%Moreover, as pointed out above, this implies that condition (C 1) holds as well. 

\end{example} 

\bigskip
{\bf The obstruction group}

\medskip Recall from \S \ref{Z} the definition of the ``obstruction group" $\sha_{\Delta}$.
%Let us write
%$$\Delta = \underset {i \in I} \prod f_i,$$ where $f_i \in {\bf Z}[X]$ are distinct irreducible, symmetric polynomials. 
Let $I$ be the set of irreducible factors of $\Delta$. 
For all $f,g \in I$, let $V_{f,g}$ be the set of prime numbers such that $f$ and $g$ have a common irreducible,
symmetric factor ${\rm (mod \ p)}$. 
%Recall that $C(I)$ is the set of maps $I \to {\bf Z}/2{\bf Z}$, and let
%$C_{\sim}(I)$ 
%be the set of maps $c \in C(I)$ such that $c(i) = c(j)$ if $V_{i,j} \not = \varnothing$; the group $\sha_{\Delta}$ is
%the quotient of $C_{\sim}(I)$ by the constant maps.
Consider the equivalence relation on $I$ generated by the elementary equivalence 
$$f \sim_e g \iff V_{f,g} \not = \varnothing,$$
and let $\overline I$ be the set of equivalence classes; the group $\sha_{\Delta}$ is the quotient of $C(\overline I)$ by
the constant maps.

\bigskip

By Corollary \ref{all Milnor}, we have the following :

\begin{coro}\label{knot trivial sha} Assume that conditions  {\rm (C 1)} and {\rm (C 2)} hold. If $\sha_{\Delta} = 0$, then for
all $\tau \in {\rm Mil}_{r,s}(f)$ there exists a knot with Milnor index $\tau$. In particular, there
exists a knot with index $r - s$. 

\end{coro}

If $\sha_{\Delta} \not = 0$, then we obtain a necessary and sufficient condition for  $\tau \in {\rm Mil}_{r,s}(f)$ to
be the Milnor index of a knot (cf. Theorem \ref{final Mil}).
The aim
of the next section is to illustrate this by some examples.

\section{Alexander polynomials of torus knots and indices }\label{torus}

\medskip Let $u,v > 1$ be two odd integers, and assume that $u$ and $v$ are prime to each other. Set

$$\Delta_{u,v} = {{(X^{uv} - 1)(X-1)} \over {(X^u-1)(X^v-1)}}.$$

\medskip
It is well-known that  $\Delta_{u,v}$ is the Alexander polynomial of the $(u,v)$-torus knot; the indices of torus knots
have been studied in many papers (see for instance \cite{Bo}, \cite {BO}, \cite {C},  \cite{GLM},  \cite{HM}, \cite {Ka} Chapter XII,  \cite{Kearton}, \cite{KM},  \cite{Lith},   \cite {Liv 18}, \cite{Mat}, \cite{Mu 06}). 

\medskip
The aim of this section is to determine which indices occur
for knots with Alexander polynomial $\Delta_{u,v}$. 

\medskip
The polynomial $\Delta_{u,v}$  is a product of cyclotomic polynomials 
$$\Delta_{u,v} = \prod \Phi_{\alpha \beta}$$ where the product is over the set of all pairs  $(\alpha,\beta)$
for which $\alpha$ is a factor of $u$, $\beta$ is a factor of $v$, and both are greater than $1$ (see for instance 
\cite {Liv 18}, Lemma 2.1). Set ${\rm deg}(\Delta_{u,v}) = 2n$, and note that $n$ is even, and that
$\Delta_{u,v}(1) = \Delta_{u,v}(-1) = 1$; hence condition (C 1) holds. 

\medskip Let $(r,s)$ be a pair of integers, $r, s \geqslant 0$, such that  $$r+s = 2n \ {\rm  and} \ r  \equiv s \ {\rm (mod \ 8)}.$$
By Example \ref{c1c2 torus}, condition (C 2) holds for $\Delta_{u,v}$ and $(r,s)$. 
In order to apply the Hasse principle results of \S \ref{Z} and \S \ref{Z+},  the next step is to determine the obstruction group $\sha_{\Delta_{u,v}}$.

\medskip
{\bf The group $\sha_{\Delta_{u,v}}$}

\medskip
If $f$ and $g$ are monic polynomials in ${\bf Z}[X]$, let $V_{f,g}$ be the set of prime numbers $p$ such that 
$f$ ${\rm (mod \ p)}$ and $g$ ${\rm (mod \ p)}$ have a common irreducible and symmetric factor in ${\bf F}_p[X]$. 
Let $I$ be the set of irreducible factors of $\Delta_{u,v}$, and recall that the definition of the group $\sha_{\Delta_{u,v}}$
involves the equivalence relation on $I$ generated by the elementary equivalence 
$$f \sim_e g \iff V_{f,g} \not = \varnothing.$$
Therefore in order to determine the group $\sha_{\Delta_{u,v}}$, we need to decide when $V_{f,g} \not = \varnothing$ for two cyclotomic
polynomials $f$ and $g$. We start with a lemma

\begin{lemma}\label {cyclotomic resultant} Let $m$ and $m'$ be two odd integers with $m \geqslant 3$ and $m' > m$. 
Let $p$ be a prime number. The cyclotomic polynomials $\Phi_m$ and $\Phi_{m'}$ have a common factor ${\rm (mod \ p)}$ 
if and only if $m' = mp^e$ for some integer $e \geqslant 1$. 

\end{lemma}

\noindent
{\bf Proof.} The resultant of the cyclotomic polynomials  $\Phi_m$ and $\Phi_{m'}$ is divisible by $p$ if and only
if $m' = mp^e$ for some integer $e \geqslant 1$ (see for instance \cite{Ap}, Theorem 4).
%\cite {St}, Proposition 3.4). 
By a well-known property
of the resultant, this is equivalent for $\Phi_m$ and $\Phi_{m'}$ to have a common factor ${\rm (mod \ p)}$ .

\medskip
We need a criterion for the existence of {\it symmetric}  irreducible factors ${\rm (mod \ p)}$; this will be done in the
following proposition, relying on well-known properties of cyclotomic polynomials and cyclotomic fields. For all integers $m \geqslant 2$, let
$\zeta_m$ be a primitive $mth$ root of unity, and let ${\bf Q}(\zeta_m)$ be the corresponding cyclotomic field.

\begin{prop} \label{m and p} Let $m \geqslant 3$ be an odd integer, and let $p$ be an odd prime number. 
The following properties are equivalent

\medskip
{\rm (a)} The polynomial $\Phi_m$ has a symmetric irreducible factor ${\rm (mod \ p)}$.

\medskip
{\rm (b)} The prime ideals above $p$ in ${\bf Q}(\zeta_m + \zeta_m^{-1})$ are inert in ${\bf Q}(\zeta_m)$. 

\medskip
{\rm (c)} The subgroup of $({\bf Z}/m{\bf Z})^{\times}$ generated by $p$ contains $-1$.

\end{prop} 

\noindent
{\bf Proof.} The equivalence (a) $ \iff $ (b)  is well-known, see for instance \cite {W}, Proposition 2.14.  Let us show that (a) $ \iff $ (c).  Let $A$ be
the group $({\bf Z}/m{\bf Z})^{\times}$  written additively, let $m_p : A \to A$ be the multiplication by $p$, and let
$\epsilon : A \to A$ be the map sending $a$ to $-a$. Both (a) and (c)  are equivalent to

\medskip
(d) There exists an orbit of $m_p$  stable by $\epsilon$. 

\medskip
This concludes the proof of the proposition. 

\bigskip Proposition \ref{m and p} and Lemma \ref{cyclotomic resultant} suffice to determine the obstruction group
of a product of cyclotomic polynomials, in particular the group $\sha_{\Delta_{u,v}}$. A full description would be
rather heavy, so we only give some simple special cases and examples. The first remark is that 
if  moreover $m$ is a prime and $m  \equiv 3 \ {\rm (mod \ 4)}$, property (c) of Proposition \ref{m and p} takes a
very simple form.

\begin{coro}\label {q and p}  Let $m$ be a prime number with $m  \equiv 3 \ {\rm (mod \ 4)}$, and let $p$ be an odd prime number. Then 

\smallskip
\centerline {$\Phi_m$ has a symmetric irreducible factor ${\rm (mod \ p)}$ $\iff$ ${\rm ({p \over m})} =  - 1$.}

\end{coro} 

\noindent
{\bf Proof.} Indeed, it is easy to see that subgroup of $({\bf Z}/m{\bf Z})^{\times}$ generated by $p$ contains $-1$
if and only if $p$ is not a square ${\rm (mod \ m)}$, in other words, ${\rm ({p \over m})} =  - 1$. Therefore the corollary
follows from Proposition \ref{m and p}. 

\begin{example}\label{sha pq} Let $p$ and $q$ be two distinct odd prime numbers, and assume that $q  \equiv 3 \ {\rm (mod \ 4)}$.
Let $e \geqslant 1$ be an integer, and set $\Delta = \Delta_{p^e,q}$. Then 

\bigskip
\centerline {$\sha_{\Delta} = 0$ \  if \ ${\rm ({p \over q})} =  - 1$, and  $\sha_{\Delta} \simeq ({\bf Z}/2{\bf Z})^{e-1}$ \  if  \ ${\rm ({p \over q})} =  1$.}

\bigskip Indeed, let $I$ be the set of irreducible factors of $\Delta$; the set $I$ consists of the cyclotomic
polynomials $\Phi_{p^kq}$ for $1 \leqslant k \leqslant e$. If ${\rm ({p \over q})} =  - 1$, then by Corollary \ref{q and p} and
Lemma \ref{cyclotomic resultant} all these polynomials have common symmetric, irreducible factors ${\rm (mod \ p)}$, hence  $V_{f,g} = \{p \}$ for all $f,g \in I$.  Therefore all the polynomials are equivalent, and the set $\overline I$ of equivalence
classes has one element; $C(\overline I)$ modulo the constant maps is trivial, hence $\sha_{\Delta} = 0$.

\medskip
 On the other hand, if ${\rm ({p \over q})} =  1$, then by the above results 
$V_{f,g} = \varnothing$ for all $f,g \in I$. All the polynomials are in different equivalence classes, hence $\overline I$ has
$e$ elements, $C(\overline I) \simeq ({\bf Z}/2{\bf Z})^{e}$,
and $\sha_{\Delta} \simeq ({\bf Z}/2{\bf Z})^{e-1}$.

\end{example} 

\begin{example}\label{three torus} Let $p$, $p_1$ and $p_2$ be distinct prime numbers with  $p \equiv p_1 \equiv p_2   \equiv 3 \ {\rm (mod \ 4)}$,
and set $\Delta = \Delta_{p,p_1p_2}$. Since $p_1$ and $p_2$ play symmetric roles and  $ ({p_2 \over p_1}) =   - ({p_1 \over p_2})$, we may
assume that  $({p_1 \over p_2})=  1$. We have 

$$\sha_{\Delta} \simeq ({\bf Z}/2{\bf Z})^2 \iff ({p_2 \over p}) = 1,$$ and 

$$\sha_{\Delta} \simeq {\bf Z}/2{\bf Z}) \iff ({p_2 \over p}) = - 1.$$  In particular, $\sha_{\Delta}$ cannot be trivial. 

\medskip
Indeed, let $I$ be the set of irreducible factors of $\Delta$; the set $I$ consists of the cyclotomic polynomials $\Phi_{pp_1}$, $\Phi_{pp_2}$ and
$\Phi_{pp_1p_2}$. The first two are not equivalent (cf. Lemma \ref{cyclotomic resultant}). Note that since $({p_1 \over p_2})=  1$, the polynomials
$\Phi_{pp_2}$ and $\Phi_{pp_1p_2}$ are not equivalent; this follows from Proposition \ref{m and p}.  Hence $\overline I$ has at least two elements, and therefore  $\sha_{\Delta} \not = 0$. On the
other hand, since $({p_2 \over p_1}) = -1$, Proposition \ref{m and p} shows that 
$\Phi_{pp_1}$ and $\Phi_{pp_1p_2}$ are equivalent if and only if $({p_2 \over p}) = - 1.$

\end{example}

\bigskip We already know that if $\sha_{\Delta_{u,v}} = 0$, then all signatures $(r,s)$ as above (hence all indices $r-s$) 
are possible, as well as all Milnor indices (see Corollary \ref{knot trivial sha}).  If however $\sha_{\Delta_{u,v}} \not = 0$,
then we need further information, in particular the homomorphism $\epsilon$.

\medskip
{\bf The homomorphism $\epsilon : \sha_{\Delta_{u,v}} \to {\bf Z}/2{\bf Z}$}

\medskip Recall that $V'$ is the set of finite places of $\bf Q$; that is, the set of prime numbers, and that
$\epsilon : \sha_{\Delta_{u,v}} \to {\bf Z}/2{\bf Z}$ is defined in terms of the local data associated to $V'$. Since $\Delta_{u,v}$
is a product of cyclotomic polynomials, $\epsilon$ can be described explicitly, as follows.

\medskip As above, we denote by $I$ the set of irreducible factors of $\Delta_{u,v}$. If $f \in I$ with $f = \Phi_m$ and 
if $p$ is a prime number, set

%$\delta^p_f = 0$ if $p$ does not divide $m$, or if $p$ divides $m$ and $p  \equiv 1 \ {\rm (mod \ 4)}$; 

\medskip
$\delta^p_f = 1$ if $p$ divides $m$, $p  \equiv 3 \ {\rm (mod \ 4)}$, and the polynomial $\Phi_m$ has a symmetric irreducible factor ${\rm (mod \ p)}$; 

\medskip
$\delta^p_f = 0$ otherwise. 

\medskip For all prime numbers $p$, let $a^p : I \to {\bf Z}/2{\bf Z}$ be defined by $a^p(f) = \delta^p_f$.

\begin{prop}\label {epsilon tori} The homomorphism $\epsilon : \sha_{\Delta_{u,v}} \to {\bf Z}/2{\bf Z}$ is given by

$$\epsilon(c) = \underset{p \in V'} \sum  \ \underset{f \in I} \sum c(f) a^p(f).$$

\end{prop} 

\noindent
{\bf Proof.} For all $f \in I$, set $E_f = {\bf Q}[X]/(f)$ and let $\sigma_f : E_f \to E_f$ be the involution induced by
$X \mapsto X^{-1}$; note that $E_f$ is a cyclotomic field, and $\sigma_f$ is the complex conjugation. Let
$F_f$ be the fixed field of $\sigma_f$.  Let $d_f$ be the discriminant of $E_f$, and
let $p^{e_p(f)}$ be the power of $p$ dividing $d_f$. Then $e_p(f) \equiv \delta^p_f \ {\rm (mod \ 2)}$ (see for
instance \cite {W}, Propositions 2.1 and 2.7). 

\medskip Since $f = \Phi_m$ for some integer $m$ that is not a prime power, no
finite places of $F_f$ ramify in $E_f$ (see for instance \cite{W}, Proposition 2.15 (b)).  Moreover, the prime ideals of $F_f$ above $p$ are 
inert in $E_f$ if and only if 
$\Phi_m$ has a symmetric irreducible factor ${\rm (mod \ p)}$; cf. Proposition \ref{m and p}. 
Hence Corollary 6.2 and Proposition 6.4 of \cite {BT} imply that we can take $a^p(f) = \delta^p_f$. The proposition
now follows from the definition of $\epsilon : \sha_{\Delta_{u,v}} \to {\bf Z}/2{\bf Z}$  in \S \ref{isometrieslattices}, where it is also shown that the homomorphism
is independent of the choice of local data.

\bigskip
{\bf Milnor indices and the homomorphism $\epsilon_{\tau} :  \sha_{\Delta_{u,v}} \to {\bf Z}/2{\bf Z}$}

\medskip
Let $\tau \in {\rm Mil}_{r,s}(\Delta_{u,v})$; the homomorphism $\epsilon_{\tau} :  \sha_{\Delta_{u,v}} \to {\bf Z}/2{\bf Z}$
is as follows. For all $f \in I$, let $n(f)$ be the number of $\mathcal P \in {\rm Irr}(\bf R)$ dividing $f$ such that 
$\tau(\mathcal P) = -2$. Let $a_{\tau} : I \to {\bf Z}/2{\bf Z}$ be defined by $a_{\tau}(f) = n(f)$ ${\rm (mod \ 2)}$. Then
$\epsilon_{\tau} :  \sha_{\Delta_{u,v}} \to {\bf Z}/2{\bf Z}$ is by definition 

$$\epsilon_{\tau}(c) =   \underset{f \in I} \sum c(f) \ a_{\tau}(f).$$

\bigskip
{\bf A necessary and sufficient criterion}

\medskip Applying Theorem \ref{final Mil}, we get the following result

\begin{theo}\label{final torus} Let $\tau \in {\rm Mil}_{r,s}(\Delta_{u,v})$. There exists a knot with Alexander polynomial $\Delta_{u,v}$ and Milnor index $\tau$ if and only if $\epsilon + \epsilon_{\tau} = 0$. 

\end{theo} 

\noindent {\bf Proof.} This is an immediate consequence of Theorem \ref{final Mil}, noting that the local
conditions always hold.

\bigskip We illustrate this result by the following example :

\begin{example}\label{final pq} We keep the notation of Example \ref{sha pq}; in particular,
$p$ and $q$ are two distinct odd prime numbers with $q  \equiv 3 \ {\rm (mod \ 4)}$,
$e \geqslant 1$ is  an integer, and  $\Delta = \Delta_{p^e,q}$. 

\medskip $\bullet$ If  ${\rm ({p \over q})} =  - 1$, then $\sha_{\Delta} = 0$. This implies that all $\tau \in {\rm Mil}_{r,s}(\Delta)$
occur as Milnor indices of knots. An integer $\iota$ is the index of a knot 
with Alexander polynomial $\Delta$ if and only if $\iota  \equiv 0 \ {\rm (mod \ 8)}$,
and $|\iota| \leqslant 2n$ (see Example \ref{sha pq} and Corollary \ref{knot trivial sha}). 

\medskip $\bullet$ If  ${\rm ({p \over q})} =  1$, then $\sha_{\Delta} \simeq ({\bf Z}/2{\bf Z})^{e-1}$ (cf Example \ref{sha pq}). 
Let $\tau \in {\rm Mil}_{r,s}(\Delta)$. In order to decide whether $\tau$ occurs as the Milnor index of a knot, we determine
the homomorphisms $\epsilon$ and $\epsilon_{\tau}$.

\medskip Recall that $I$ is the set of irreductible factors of $\Delta$, and that $I$ consists of the cyclotomic
polynomials $\Phi_{p^kq}$ for $1 \leqslant k \leqslant e$.

\medskip Still assuming   ${\rm ({p \over q})} =  1$,  the result depends on the congruence class of $p$ ${\rm (mod \ 4)}$ :

\medskip

$\bullet$ Assume first that $p  \equiv 3 \ {\rm (mod \ 4)}$, and note that ${\rm ({q \over p})} = - {\rm ({p \over q})} = -1$.

\medskip If $f \in I$, we have $\delta^q_f = 1$ and
$\delta^{\ell}_f = 0$ for all prime numbers $\ell \not = q$. 

\medskip
For all prime numbers $\ell$, let $a^{\ell} : I \to {\bf Z}/2{\bf Z}$ be defined by $a^{\ell}(f) = \delta^{\ell}_f$, hence
$a^{q}(f) = 1$  and $a^{\ell}(f) = 0$ if $\ell \not = q$. The 
associated homomorphism $\epsilon : \sha_{\Delta} \to {\bf Z}/2{\bf Z}$ is given by
$$\epsilon(c) = \underset{\ell \in V'} \sum  \ \underset{f \in I} \sum c(f) a^\ell(f) =  \underset{f \in I} \sum c(f).$$
On the other hand, we have
$$\epsilon_{\tau}(c) =   \underset{f \in I} \sum c(f) \ a_{\tau}(f).$$ Therefore 

\medskip \centerline {$\epsilon + \epsilon_{\tau} = 0 \iff
 a_{\tau}(f) = 1$ for all $f \in I$.}

\medskip By definition, this means that $n(f)$ is odd for all $f \in I$, in other words, that each $f \in I$ is divisible by
an odd number of $\mathcal P \in {\rm Irr}({\bf R})$ with $\tau (\mathcal P) = -2$. 
In summary, we have :

\medskip An element $\tau \in {\rm Mil}_{r,s}(\Delta)$ occurs as the Milnor index of a knot 
if and only if for all $f \in I$, the number
of $\mathcal P \in {\rm Irr}({\bf R})$ dividing $f$  such that  $\tau (\mathcal P) = -2$ is odd.

\medskip
In particular, an  integer $\iota$ is the index of a knot with Alexander polynomial $\Delta$ if and only if $\iota  \equiv 0 \ {\rm (mod \ 8)}$, and
$|\iota| \leqslant 2n - 4(e-1)$.

\medskip
$\bullet$ Assume now that $p  \equiv 1 \ {\rm (mod \ 4)}$, and note that ${\rm ({q \over p})} =  {\rm ({p \over q})} = 1$. In this case, 
$\delta^{\ell}_f = 0$ for all prime numbers $\ell$ and for all
$f \in I$. This implies that $\epsilon (c) = 0$ for all $c \in  \sha_{\Delta}$. Recall that $\epsilon_{\tau} : \sha_{\Delta} \to {\bf Z}/2{\bf Z}$ is given by
$$\epsilon_{\tau}(c) =   \underset{f \in I} \sum c(f) \ a_{\tau}(f).$$ Therefore

 \medskip \centerline {$\epsilon + \epsilon_{\tau} = 0 \iff
 a_{\tau}(f) = 0$ for all $f \in I$.}
 
 \medskip This means that $n(f)$ is even for all $f \in I$, in other words, that each $f \in I$ is divisible by
an even number of $\mathcal P \in {\rm Irr}({\bf R})$ with $\tau (\mathcal P) = -2$. 
This imposes a condition on the possible Milnor indices. However, this does not gives rise to  additional conditions on the index; 
an integer $\iota$ is the index of a knot 
with Alexander polynomial $\Delta$ if and only if $\iota  \equiv 0 \ {\rm (mod \ 8)}$,
and $|\iota| \leqslant 2n$.

%hence $\epsilon(c) = \underset{f \in I} \sum c(f) a^q(f) = \underset{f \in I} \sum c(f)$. 

\end{example}

\

\bigskip

\bigskip

Eva Bayer--Fluckiger 

EPFL-FSB-MATH

Station 8

1015 Lausanne, Switzerland

\medskip

eva.bayer@epfl.ch

\end{document}